\documentclass[11pt]{article}

\usepackage{fullpage}
\usepackage{natbib}
\usepackage{algorithm}
\usepackage{algorithmic}

\usepackage{comment} 
\usepackage{lipsum} 
\usepackage{amsmath}
\usepackage[utf8]{inputenc}
\usepackage[T2A]{fontenc}
\usepackage{amssymb}  
\usepackage{algorithm}
\usepackage{algorithmic}
\usepackage{graphicx}
\usepackage[dvipsnames]{xcolor}
\usepackage{nicefrac}
\usepackage{bm}

\usepackage{threeparttable}
\usepackage{makecell}

\usepackage{multirow}
\usepackage{colortbl}
\definecolor{bgcolor}{rgb}{0.8,1,1}
\definecolor{bgcolor2}{rgb}{0.8,1,0.8}

\usepackage{graphicx} 

\newcommand{\cR}{{\cal R}}

\newcommand{\dataset}[1]{{\tt#1\xspace}}

\newcommand{\R}{\mathbb{R}}

\newcommand{\prox}{\operatorname{prox}}

\newcommand{\eqdef}{\; { := }\;}

\usepackage{xspace}

\def\<{\left\langle}
\def\>{\right\rangle}
\def\[{\left[}
\def\]{\right]}
\def\({\left(}
\def\){\right)}

\usepackage{tcolorbox}
\usepackage{pifont}
\definecolor{mydarkgreen}{RGB}{39,130,67}
\definecolor{mydarkred}{RGB}{192,47,25}

\usepackage{amsmath,amsfonts,amssymb,amsthm,array}

\usepackage{mdframed} 
\usepackage{thmtools}
\usepackage{textcomp}

\declaretheorem[within=section]{definition}
\declaretheorem[sibling=definition]{theorem}

\declaretheorem[sibling=definition]{assumption}

\declaretheorem[sibling=definition]{lemma}

\declaretheorem[sibling=definition]{remark}

\usepackage[colorinlistoftodos,bordercolor=orange,backgroundcolor=orange!20,linecolor=orange,textsize=scriptsize]{todonotes}

\usepackage{microtype}
\usepackage{subfigure}
\usepackage{booktabs} 
\usepackage{grffile}
\usepackage{hyperref}

\title{\bf Catalyst Acceleration of Error Compensated Methods Leads to Better Communication Complexity}

\author{Xun Qian\thanks{Shanghai Artificial Intelligence Lab, Shanghai, China.} \and Hanze Dong \thanks{The Hong Kong University of Science and Technology, Hong Kong.} \and Tong Zhang\thanks{The Hong Kong University of Science and Technology, Hong Kong.}   \and Peter Richt\'{a}rik\thanks{King Abdullah University of Science and Technology, Thuwal, Saudi Arabia.}}
\date{January 24, 2023}

\begin{document}

\maketitle

\begin{abstract}
 Communication overhead is well known to be a key bottleneck in large scale distributed learning, and a particularly successful class of methods which help to overcome this bottleneck is based on the idea of communication compression. Some of the most practically effective gradient compressors, such as TopK, are biased, which causes convergence issues unless one employs a well designed {\em error compensation/feedback} mechanism. Error compensation is therefore a fundamental technique in the distributed learning literature. In a recent development, Qian et al (NeurIPS 2021) showed that the error-compensation mechanism can be combined with acceleration/momentum, which is another key and highly successful optimization technique. In particular, they developed the error-compensated loop-less Katyusha (ECLK) method, and proved an accelerated linear rate in the strongly convex case. However, the dependence of their rate on the compressor parameter does not match the best dependence obtainable in the non-accelerated error-compensated methods. Our work addresses this problem. We propose several new accelerated error-compensated methods using the {\em catalyst acceleration} technique, and obtain results that match the best dependence on the compressor parameter in non-accelerated error-compensated methods up to logarithmic terms.
\end{abstract}

\tableofcontents
 
\clearpage 

\section{INTRODUCTION}

In large scale machine learning optimization problems, the data and training need to be distributed among many machines \citep{DMLsurvey}. Also in federated learning \citep{FEDLEARN, FEDOPT, FL2017-AISTATS,  FL_survey_2019}, training occurs on edge devices such as mobile phones and  smart home devices, where the data is originally captured. In these applications, the distributed machine learning can be characterized as the following composite finite-sum problem 
\begin{equation}\label{primal-LSVRG}
	\min \limits_{x\in \mathbb{R}^d} P(x) \eqdef \left\{\frac{1}{n} \sum\limits_{\tau=1}^n  f^{(\tau)}(x) + \psi(x)\right\},
\end{equation}
where $\{f^{(\tau)}(x)\}_{\tau=1}^n$ are smooth convex functions distributed over $n$ nodes, and $\psi:\R^d \to \R\cup \{+\infty\}$ is a regularizer, which is a proper closed convex but possibly non-smooth function. On each node $\tau$, $f^{(\tau)}(x) \eqdef  \frac{1}{m} \sum \limits_{i=1}^m f^{(\tau)}_i(x)$ is the average loss over the training data stored on this node and each $f_i^{(\tau)}$ is smooth and convex.  


In distributed and especially federated settings, communication is generally much slower than the local training, which makes the communication overhead become a key bottleneck. In order to overcome this bottleneck, several methods were proposed in the literature, such as using large mini-batches \citep{Goyal17, You17}, asynchronous learning \citep{tsitsiklis1986distributed, Agarwal11, Lian15, Recht11},  and gradient compression \citep{Seide14, Alistarh17, Bernstein18, Wen17, Mish19}. In this work, we focus on the error-compensated method, which is a gradient compression method and is capable to deal with some effective but biased compressors, such as the TopK compressor.

{\bf Related Work.} The {\em error compensation/feedback} mechanism was first introduced in 1-bit SGD \citep{Seide14}. Then the error-compensated SGD (ECSGD) was proved to have the same convergence rate as vanilla SGD in the strongly convex case \citep{Stich18} and non-convex case \citep{karimireddy2019error,Tang19} when $P$ is smooth. ECSGD was further studied in \citep{stich2020error} under weaker assumptions. When $P$ is non-smooth, it was shown that ECSGD converges at the rate of ${\cal O}(\nicefrac{1}{\sqrt{\delta T}})$ in \citep{karimireddy2019error}, where $T$ denotes the iteration number and $\delta$ is the compressor parameter defined in (\ref{eq:contractor}). The non-accelerated linear convergence can be obtained in EC-LSVRG-DIANA \citep{gorbunov2020linearly} in the smooth case, and in the error-compensated loop-less SVRG, Quartz, and SDCA  \citep{ecsdca} in the composite case. In a recent development, the error-compensated loop-less Katyusha was proposed in \citep{eclk}, and the accelerated linear rate was achieved.

{\bf Compressor.} In error-compensated methods, contraction compressors are generally used. A randomized map $Q:\R^d\to \R^d$ is called a {\em  contraction compressor} if there exists a constant $\delta \in (0, 1]$ such that  
\begin{equation}\label{eq:contractor}
	\mathbb{E} \left[\|x - Q(x) \|^2\right] \leq (1-\delta)\|x\|^2, \qquad \forall x\in \R^d. 
\end{equation}
\noindent Some frequently used contraction compressors include TopK \citep{Alistarh18} and RandK \citep{Stich18}. Let $1\leq K \leq d$. The TopK compressor  is defined as
$$
({\rm TopK}(x))_{\pi(i)} = \left\{ \begin{array}{rl}
	(x)_{\pi(i)}  &\mbox{ if $i\leq K$, } \\
	0  \quad \quad &\mbox{ otherwise, }
\end{array} \right.
$$
where $\pi$ is a permutation of $\{1, 2, ..., d\}$ such that $(|x|)_{\pi(i)} \geq (|x|)_{\pi(i+1)}$ for $i = 1, ..., d-1$. For TopK and RandK compressors, we have $\delta \geq \nicefrac{K}{d}$ \citep{Stich18}.

The {\em unbiased compressor} is also frequently used in compression algorithms, which is defined as a randomized map $\tilde{Q}:\R^d\to \R^d$, where there exists a constant $\omega \geq 0$ such that $\mathbb{E}[{\tilde Q}(x)] = x$, and 
\begin{equation} \label{eq:unbiased}
	\mathbb{E} \left[\|{\tilde Q}(x)\|^2 \right] \leq (\omega + 1)\|x\|^2, \qquad \forall x\in \R^d.
\end{equation}
Some frequently used unbiased compressors include random dithering \citep{Alistarh17}, random sparsification \citep{Stich18}, and natural compression \citep{horvath2019natural}. For any ${\tilde Q}$ satisfying \eqref{eq:unbiased}, $\frac{1}{\omega+1}{\tilde Q}$ is a contraction compressor satisfying \eqref{eq:contractor} with $\delta = \nicefrac{1}{(\omega+1)}$\citep{biased2020}. Furthermore, unbiased compressors and contraction compressors can be composed to generate new contraction compressors \citep{ecsdca}.

\begin{table*}[h]
	\caption{Communication Complexity Results for Different Error-Compensated Algorithms ($r_Q$ represents the communication cost of the compressed vector $Q(x)$ for $x\in \R^d$. For simplicity, we choose $Q=Q_1$, and assume $L_f \geq \lambda$, $\nicefrac{R^2}{\gamma} \geq \lambda$, where $R$ is defined in Algorithm \ref{alg:ec-spdc}, hence the term $\nicefrac{1}{\delta}$ is omitted.) }
	\label{tab:summary}
	\begin{center}
		{\footnotesize
			\begin{tabular}{|c|c|c|}
				\hline
				Algorithm & \begin{tabular}{c} Communication complexity  \\ when $\delta \leq \nicefrac{1}{m}$	\end{tabular}&    \begin{tabular}{c}  Communication complexity under \\ Assumption    \ref{as:expcompressor} when $\delta \leq \nicefrac{1}{m}$	\end{tabular}\\
				\hline 
				\begin{tabular}{c} 
					EC-LSVRG \\
					Smooth Case  \citep{ecsdca}
				\end{tabular}  & 
				${\cal O} \left(  \tfrac{r_Q}{\delta} \tfrac{\sqrt{L_f {\bar L}}}{ \lambda}  \log \tfrac{1}{\epsilon} \right)$ & ${\cal O} \left(  \tfrac{r_Q}{\delta} \tfrac{L_f}{ \lambda}  \log \tfrac{1}{\epsilon} \right)$  \\
				\hline
				\begin{tabular}{c} 
					EC-SDCA \\
					\citep{ecsdca}
				\end{tabular}  & 
				${\cal O} \left(  \tfrac{r_Q}{\delta} \tfrac{R{\bar R}}{ \lambda \gamma}  \log \tfrac{1}{\epsilon} \right)$  & ${\cal O} \left(  \tfrac{r_Q}{\delta} \tfrac{R^2}{ \lambda \gamma}  \log \tfrac{1}{\epsilon} \right)$   \\
				\hline
				\begin{tabular}{c} 
					ECLK \\
					\citep{eclk}
				\end{tabular}  & 
				${\cal O} \left(  \tfrac{r_Q}{\delta \sqrt{\delta}} \sqrt{\tfrac{\bar L}{\lambda}} \log \tfrac{1}{\epsilon} \right)$ &  ${\cal O} \left(   \tfrac{r_Q}{\delta \sqrt{\delta}} \sqrt{\tfrac{L_f}{\lambda}} \log \tfrac{1}{\epsilon} \right)$  \\
				\hline
				\begin{tabular}{c} 
					ECSPDC\\
					{\bf This work}
				\end{tabular}  & 
				${\cal O} \left(   \tfrac{r_Q}{\delta^2 \sqrt{m}} \sqrt{\tfrac{{\bar R}^2}{\lambda \gamma}}  \log \tfrac{1}{\epsilon} \right)$ & ${\cal O} \left(  \tfrac{r_Q}{\delta^2 \sqrt{m}} \sqrt{\tfrac{R^2}{\lambda \gamma}}  \log \tfrac{1}{\epsilon} \right)$  \\
				\hline
				\begin{tabular}{c} 
					EC-LSVRG + Catalyst\\
					Smooth Case \  \ 	{\bf This work}
				\end{tabular}  & 
				${\tilde {\cal O}} \left(  \tfrac{r_Q}{\delta} \sqrt{\tfrac{\bar L}{\lambda}}  \log \tfrac{1}{\epsilon}  \right)$ &  ${\tilde {\cal O}} \left(  \tfrac{r_Q}{\delta} \sqrt{\tfrac{L_f}{\lambda}}  \log \tfrac{1}{\epsilon}  \right)$ \\
				\hline
				\begin{tabular}{c} 
					EC-SDCA + Catalyst\\
					{\bf This work}
				\end{tabular}  & 
				${\tilde {\cal O}} \left(  \tfrac{r_Q}{\delta} \sqrt{\tfrac{{\bar R}^2}{\lambda \gamma}} \log \tfrac{1}{\epsilon}  \right)$  & ${\tilde {\cal O}} \left(  \tfrac{r_Q}{\delta} \sqrt{\tfrac{R^2}{\lambda \gamma}} \log \tfrac{1}{\epsilon}  \right)$ \\
				\hline
			\end{tabular}
		}
	\end{center} 
\end{table*}

\subsection{Motivation}

{\bf Communication Complexity of ECLK.} There are two contraction compressors $Q$ and $Q_1$ in ECLK \citep{eclk} with parameter $\delta$ and $\delta_1$ respectively. We first claim that when $Q$ and $Q_1$ in ECLK are the same type of contraction compressor, but with possibly different compressor parameters (for example, $Q$ and $Q_1$ are both TopK, but with different values of $K$), we could always choose the same compressor parameters for $Q$ and $Q_1$ such that the total communication complexity is less than before or remains the same order as before.

First, from the iteration complexity results for ECLK, it is easy to verify that the iteration complexity will decrease as $\delta$ or $\delta_1$ increases. Without less of generality, we assume the communication cost of $Q(x)$ is higher than that of $Q_1(x)$. Since $Q$ and $Q_1$ are the same type of compressor, we will have $\delta_1 \leq \delta$. Then we can change $Q_1$ to be $Q$. In this way, the total communication cost of $Q(x)$ and $Q_1(x)$ at each iteration is at most twice as before, but $\delta_1$ will increase to $\delta$, which implies that the iteration complexity will decrease and the communication complexity is at most twice as before. Thus, for simplicity, we consider $Q=Q_1$ for ECLK.

{\bf Dependence on $\delta$ for the Iteration Complexity of ECLK.} We introduce the following assumption for Problem (\ref{primal-LSVRG}). 

\begin{assumption}\label{as:primal}
	$\tfrac{1}{n} \sum_{\tau=1}^n f^{(\tau)}$ is $L_f$-smooth, $f^{(\tau)}$ is ${\bar L}$-smooth, $f_i^{(\tau)}$ is $L$-smooth, and $\psi$ is $\lambda$-strongly convex. 
\end{assumption}

Under Assumption \ref{as:primal}, from Theorem 3.8 in \citep{eclk}, the iteration complexity is 
$$
{\cal O} \left(  \left(  \tfrac{1}{\delta} + \tfrac{1}{\delta_1} + \tfrac{1}{p} + \sqrt{\tfrac{L_f}{\lambda}} + \sqrt{\tfrac{{\cal L}_2}{\lambda p}}  \right)  \log\tfrac{1}{\epsilon}  \right),
$$ 
where $${\cal L}_2 = \tfrac{6L}{n} +  \tfrac{112(1-\delta) {\bar L}}{3\delta^2} + \tfrac{28(1-\delta) L}{3\delta} + \tfrac{224(1-\delta) {\bar L}p}{\delta^2 \delta_1} \left(  1 + \tfrac{2p}{\delta_1}  \right)$$ and $p\in(0, 1]$ is the update frequency of the check point. Considering $\delta_1=\delta$, it is easy to see that the iteration complexity of ECLK is at least ${\cal O} \left(  \tfrac{\sqrt{1-\delta}}{\delta\sqrt{\delta}} \sqrt{\tfrac{\bar L}{\lambda}}   \log\tfrac{1}{\epsilon}   \right)$. Hence, when $1-\delta = \Theta(1)$, the dependence on $\delta$ of the communication complexity of ECLK would be $\nicefrac{1}{\delta^{\frac{3}{2}}}$, which is worse than EC-LSVRG in the smooth case and EC-SDCA in the composite case \citep{ecsdca}, where the dependence on $\delta$ is $\nicefrac{1}{\delta}$ only. This leads to the following question:

\begin{quote} {\em  Can we design provably accelerated gradient-type methods that work with contractive compressors and the dependence on the compressor parameter $\delta$ is $\nicefrac{1}{\delta}$.  }
\end{quote}

Let us first recall the results for the error-compensated non-accelerated methods. In the composite case, the dependence on the compressor parameter $\delta$ of EC-SDCA is better than that of EC-LSVRG \citep{ecsdca}. Noticing that L-SVRG \citep{hofmann2015variance,LSVRG} is a primal method and SDCA \citep{prox-SDCA} is a primal-dual method, the better dependence on $\delta$ of EC-SDCA than EC-LSVRG indicates that primal-dual methods may be more suitable for the error feedback mechanism. Therefore, it is natural to apply error feedback to SPDC \citep{SPDC}, which is an accelerated primal-dual algorithm, and expect better dependence on the compressor parameter than ECLK. We first propose error-compensated SPDC (Algorithm \ref{alg:ec-spdc}), but unfortunately, we show that the dependence on $\delta$ of error-compensated SPDC is at least $\nicefrac{1}{\delta^{\frac{3}{2}}}$. This fact makes us to consider the indirect accelerated methods.

In this work, we give a confirmed answer to the above question by applying Catalyst \citep{catalyst}, which is a generic method for accelerating first-order algorithms in the sense of Nesterov, to non-accelerated error-compensated methods, where the dependence on $\delta$ of the communication complexity could be ${\tilde {\cal O}}\left(  \nicefrac{1}{\delta}  \right)$. Here ${\tilde {\cal O}}$ hides some logarithmic terms.

\subsection{Contributions}

1, First, we propose the error-compensated SPDC (ECSPDC), which is a combination of the error feedback machanism and SPDC \citep{SPDC}, and achieve the accelerated linear convergence rate. In the special case where $\delta=1$, ECSPDC is also an extension of SPDC in the sense that $A_{i\tau}$ in problem (\ref{primal-spdc}) is a matrix rather than a vector, and the convergence rate is actually better than SPDC. Specifically, the convergence result in \citep{SPDC} does not achieve linear speed up with respect to the number of nodes, while ours can obtain linear speed up when the number of nodes is in a certain range. 

\noindent 2, We apply Catalyst \citep{catalyst} to EC-LSVRG in the smooth case and EC-SDCA in the composite case \citep{ecsdca}, respectively. The accelerated linear convergence rates are obtained for both cases, and the dependence on $\delta$ of the communication complexities is ${\tilde {\cal O}}(\nicefrac{1}{\delta})$, which matches the best dependence on the compressor parameter in non-accelerated error-compensated methods up to logarithmic terms. The communication complexities of them are summarized in Table \ref{tab:ECLSVRG+C} and Table \ref{tab:ECSDCA+C} in the Appendix, and the comparison of the communication complexity results of different error-compensated algorithms when $\delta \leq \nicefrac{1}{m}$ are summarized in Table \ref{tab:summary}.

\section{ERROR COMPENSATED SPDC}

For primal-dual methods, the following problem is usually studied: 
\begin{equation}\label{primal-spdc}
	\min_{x\in \mathbb{R}^d} P(x) \eqdef \tfrac{1}{N} \sum_{\tau=1}^n \sum_{i=1}^m \phi_{i \tau}(A_{i\tau}^\top x) + g(x),
\end{equation}
where $N = mn$ and $A_{i\tau} \in \R^{d\times t}$. Problem (\ref{primal-spdc}) is actually equivalent to Problem (\ref{primal-LSVRG}). First, by choosing $f_i^{(\tau)}(x) = \phi_{i \tau}(A_{i\tau}^\top x)$ and $\psi=g$, Problem (\ref{primal-spdc}) is a special case of Problem (\ref{primal-LSVRG}). On the other hand, by choosing $A_{i\tau}$ to be the identity matrix, $\phi_{i\tau} = f_i^{(\tau)}$ , and $g=\psi$, Problem (\ref{primal-spdc}) becomes Problem (\ref{primal-LSVRG}). For simplicity, we assume $L_f = \nicefrac{R^2}{\gamma}$, ${\bar L} = \nicefrac{{\bar R}^2}{\gamma}$, and $L = \nicefrac{R_m^2}{\gamma}$, where $R^2$, ${\bar R}^2$, and $R_m^2$ are defined in Algorithm \ref{alg:ec-spdc}. To save space, we only list the assumptions and main results here. The rest can be found in the Appendix.

\begin{assumption}\label{as:contracQQ1-ecSPDC}
	The two compressors $Q$ and $Q_1$ are contraction compressors with parameters $\delta$ and $\delta_1$, respectively. 
\end{assumption}

\begin{assumption}\label{as:ecSPDC}
	Each $\phi_{i\tau} : \R^t \to \R$ is convex and $\nicefrac{1}{\gamma}$-smooth. The regularizer $g : \R^d \to \R$ is $\lambda$-strongly convex. 
\end{assumption}
\noindent Sometimes, we will use the following assumption on the contraction compressor to get better results. 

\begin{assumption}\label{as:expcompressor}
	$\mathbb{E}[Q(x)] = \delta x$ and $\mathbb{E}[Q_1(x)] = \delta_1 x$. 
\end{assumption}

Under Assumption \ref{as:contracQQ1-ecSPDC} and Assumption \ref{as:ecSPDC} , the iteration complexity of ECSPDC is 
$$
{\cal O} \left(  \left(  \tfrac{1}{\delta} + \tfrac{1}{\delta_1} + m + \cR_2 \sqrt{\tfrac{m}{\lambda \gamma}}   \right) \log \tfrac{1}{\epsilon}  \right), 
$$
where 
\begin{align*} 
	\cR_2^2  = 2R^2 + \tfrac{2R_m^2}{n}  \quad + \tfrac{3(1-\delta)}{4}  \left(  \tfrac{14{\bar R}^2}{\delta^2}  + \tfrac{7R_m^2}{2\delta}   + \tfrac{84(1-\delta_1){\bar R}^2}{\delta^2 \delta_1^2 m^2} + \tfrac{42R_m^2}{\delta^2 \delta_1m^2} \right). 
\end{align*} 
If Assumption \ref{as:expcompressor} is further invoked, the iteration complexity is improved to 
$
{\cal O} \left(  \left(  \tfrac{1}{\delta} + \tfrac{1}{\delta_1} + m + \cR_3 \sqrt{\tfrac{m}{\lambda \gamma}}   \right) \log \tfrac{1}{\epsilon}  \right), 
$
where 
\begin{align*}
	\cR_3^2   = 2R^2 + \tfrac{2R_m^2}{n} + \tfrac{21(1-\delta)}{4} \left(  \tfrac{2R^2}{\delta^2} + \tfrac{11R_m^2}{2\delta n}  + \tfrac{12(1-\delta){\bar R}^2}{\delta^2 n}  \tfrac{12R^2}{5\delta^2 \delta_1^2 m^2}  + \tfrac{228R_m^2}{5\delta^2 \delta_1m^2 n} + \tfrac{432(1-\delta_1){\bar R}^2}{5\delta^2 \delta_1^2m^2 n} \right). 
\end{align*}

{\bf Comparison to SPDC.} If there is no compression in ECSPDC, i.e., $\delta=\delta_1=1$, the iteration complexity becomes 
$$
{\cal O} \left(  \left(  \tfrac{N}{n} + \tfrac{(\sqrt{n}R + R_m)}{n} \sqrt{\tfrac{N}{\lambda \gamma}}  \right) \log \tfrac{1}{\epsilon}  \right), 
$$
which is better than that of SPDC obtained in \citep{SPDC}: ${\cal O} \left(  \left(  \tfrac{N}{n} + R_m \sqrt{\tfrac{N}{n \lambda \gamma}}  \right) \log \tfrac{1}{\epsilon}  \right)$. Moreover, our result achieves linear speed up with repect to $n$ when $n \leq \nicefrac{R_m^2}{R^2}$. 

{\bf Dependence on $\delta$.} Consider $Q=Q_1$ in ECSPDC. When $\nicefrac{1}{m} \leq \delta$, the iteration complexity is at least ${\cal O} \left(  \tfrac{\bar R}{\delta} \sqrt{\tfrac{m}{\lambda \gamma}} \log \tfrac{1}{\epsilon}  \right) \geq {\cal O} \left(  \tfrac{1}{\delta\sqrt{\delta}} \tfrac{\bar R}{\sqrt{\lambda \gamma}} \log \tfrac{1}{\epsilon} \right)$. When $\delta \leq \nicefrac{1}{m}$, we have $\delta \leq \nicefrac{{\bar R}^2}{R_m^2}$. Then the iteration complexity becomes $${\cal O}\left(  \left(  \tfrac{1}{\delta} + \tfrac{1}{\delta^2 \sqrt{m}} \tfrac{\bar R}{\sqrt{\lambda\gamma}}  \right) \log \tfrac{1}{\epsilon}  \right) \geq {\cal O}\left(    \tfrac{1}{\delta \sqrt{\delta}} \tfrac{\bar R}{\sqrt{\lambda\gamma}}   \log \tfrac{1}{\epsilon}  \right). $$
Hence, the dependence of ECSPDC on $\delta$ is at least $\nicefrac{1}{\delta^{\frac{3}{2}}}$.

\section{EC-LSVRG + CATALYST IN THE SMOOTH CASE}

EC-LSVRG \citep{ecsdca} is a combination of L-SVRG and error feedback, and the iteration complexity has the better dependence on the compressor parameter in the smooth case than that in the non-smooth case. In this section, we apply Catalyst to EC-LSVRG in the smooth case. First, we restate the Catalyst algorithm and convergence result as follows.

\begin{algorithm}[h!]
	\caption{Catalyst}
	\label{alg:catalyst}
	\begin{algorithmic}
		\STATE {\bfseries Parameters:} $\kappa\geq 0$, $\alpha_0$, sequence $\{ \epsilon_k \}_{k\geq 0}$
		\STATE {\bfseries Initialization:} 
		$y^0 = x^0\in \R^d$; $q = \lambda/(\lambda+\kappa)$
		\STATE {\bf for} {$k = 1, 2, 3, ...$} {\bf do}
		\STATE ~~~ Find an approximate solution of the following problem 
		\begin{align*}
			x^k \approx & \arg\min_{x\in \R^d}  \left\{  G_k(x) \eqdef P(x) + \tfrac{\kappa}{2}\|x-y^{k-1}\|^2  \right\} \\
			& {\rm such \  that} \  \   G_k(x^k) - G_k^* \leq \epsilon_k
		\end{align*}
		
		\STATE ~~~ Compute $\alpha_k \in (0,1)$ from equation $\alpha_k^2 = (1-\alpha_k)\alpha_{k-1}^2 + q\alpha_k$
		\STATE ~~~ Compute 
		$$
		y^k = x^k + \beta_k(x^k - x^{k-1})  \  \ 
		{\rm with} \  \   \beta_k = \tfrac{\alpha_{k-1}(1-\alpha_{k-1})}{\alpha_{k-1}^2 + \alpha_k}$$
		\STATE {\bf end for}
	\end{algorithmic}
\end{algorithm}

\begin{theorem}\label{th:catalyst}[\citealp{catalyst}]
	Choose $\alpha_0 = \sqrt{q}$ with $q = \nicefrac{\lambda}{(\lambda+\kappa)}$ and 
	$$
	\epsilon_k = \tfrac{2}{9} (P(x^0) - P^*) (1 - \rho_0)^k \ \ {\rm with} \ \ \rho_0 < \sqrt{q}. 
	$$
	Then, Algorithm \ref{alg:catalyst} generates iterates $\{  x^k  \}_{k\geq 0}$ such that 
	\begin{equation}\label{eq:Catalyst}
		P(x^k) - P^* \leq C (1-\rho_0)^{k+1} (P(x^0) - P^*). 
	\end{equation}
	with $C = \tfrac{8}{(\sqrt{q} - \rho_0)^2}$. 
\end{theorem}

\noindent In Catalyst (Algorithm \ref{alg:catalyst}), $G_k^*$ represents the minimum of $G_k$. In Theorem \ref{th:catalyst}, $P^*$ is the minimum of $P$, and as discussed in \citep{catalyst}, the term $P(x^0)-P^*$ in $\epsilon_k$ can be replaced by its upper bound, which only affects the corresponding constant in (\ref{eq:Catalyst}). 

We use EC-LSVRG to solve the subproblem in Catalyst for the smooth case where $\psi$ is smooth in Problem (\ref{primal-LSVRG}). The main challenge is proposing suitable initial conditions for the subproblem and estimate the corresponding expected inner iteration number. 

To save space, we restate EC-LSVRG (and also EC-SDCA) in the Appendix. It should be noticed that EC-LSVRG in the smooth case is applied to the problem without the regularizer term. Thus, to minimize $G_k$, we move $\psi$ and the quadratic term $\frac{\kappa}{2} \|x-y^{k-1}\|^2$ to each $f_i^{(\tau)}$. We use subscript $(k)$ and superscript $K$ to denote the variables at the $k$-th outer iteration and $K$-th inner iteration (for example, $x_{(k)}^K$, ${\bar x}_{(k)}^K$, $x_{(k)}^*$, $e_{\tau, (k)}^K$, and $h_{\tau, (k)}^K$).

In \citep{catalyst}, the Catalyst acceleration was applied to the first-order methods whose convergence rate has the following form 
\begin{equation}\label{eq:Catalystfirst}
	G_k(z_t) - G_k^* \leq A (1 - \theta)^t (G_k(z^0) - G_k^*), 
\end{equation}
where $A$ is some constant. If we initial $h_{\tau, (k)}^0$ by the gradient of $f_i^{(\tau)} + \psi + \frac{\kappa}{2}\|\cdot - y^{k-1}\|^2$ at $x_{(k)}^0$. Then the form of the convergence rate of EC-LSVRG becomes form (\ref{eq:Catalystfirst}), and we can get the following lemma.

\begin{lemma}\label{lm:eclsvrg-1}
	Under Assumptions \ref{as:primal}, \ref{as:contracQQ1-ecSPDC} and the premise of Theorem \ref{th:catalyst}, let us run EC-LSVRG (Algorithm \ref{alg:ec-lsvrg}) to minimize $G_k$ and output $x^k \eqdef {\bar x}_{(k)}^{T_k}$, where 
	$
	T_k \eqdef \inf \{  K\geq 1, G_k({\bar x}_{(k)}^{K}) - G_k^* \leq \epsilon_k   \}
	$. 
	For the initialization of EC-LSVRG at the $k$-th outer iteration, we choose $p=\Theta(\delta_1)$, $x_{(k)}^0 = x^{k-1}$, $e_{\tau, (k)}^0 = 0$ and $h_{\tau, (k)}^0 = \nabla f^{(\tau)}(x^0_{(k)}) + \nabla \psi(x^0_{(k)}) + \kappa (x^0_{(k)} - y^{k-1})$. Then 
	\begin{align*}
		\mathbb{E}[T_k]  \leq {\tilde {\cal O}}\left( \tfrac{1}{\delta} + \tfrac{1}{\delta_1} + \tfrac{\sqrt{(1-\delta)(L_f+\lambda+\kappa)({\bar L} + \lambda+\kappa)}}{\delta (\lambda+\kappa)}  + \tfrac{L_f}{\lambda+\kappa} + \tfrac{L}{n(\lambda+\kappa)}   +  \tfrac{\sqrt{(1-\delta)(L_f+\lambda+\kappa)(L+\lambda+\kappa)}}{\sqrt{\delta}(\lambda+\kappa)} \right),
	\end{align*} 
	where the notation ${\tilde {\cal O}}$ hides some universal constants and some logarithmic dependencies in $\delta$, $\delta_1$, $\lambda$, $\kappa$, $L_f$, and $N$. 
\end{lemma}

\begin{remark}
	1, It is easy to verify that an optimal choice of $p$ in EC-LSVRG is $\Theta(\delta_1)$. Hence, we choose $p=\Theta(\delta_1)$ in Lemma \ref{lm:eclsvrg-1} (and also in Lemma \ref{lm:eclsvrg-2}) for simplicity. 
	
	2, As discussed in \citep{catalyst}, the stopping criteria in the inner loop can be checked by calculating some upper bound of $G_k({\bar x}_{(k)}^{K}) - G_k^*$, such as the duality gap. However, this would cause additional computation and also communication cost. Hence, we can actually view the inner iteration number as a parameter and use Lemma \ref{lm:eclsvrg-1} as the guidance. 
\end{remark}

If we further invoke Assumption \ref{as:expcompressor}, we can get the following lemma. Since the proof is similar to that of Lemma \ref{lm:eclsvrg-1}, we omit it.

\begin{lemma}\label{lm:eclsvrg-2}
	Under Assumptions \ref{as:primal}, \ref{as:contracQQ1-ecSPDC}, \ref{as:expcompressor}, and the premise of Theorem \ref{th:catalyst}, let us run EC-LSVRG to minimize $G_k$. Choose the output $x^k$, $T_k$, and the initialization of EC-LSVRG at the $k$-th outer iteration be the same as that in Lemma \ref{lm:eclsvrg-1}. Then 
	$$
	\mathbb{E}[T_k] \leq  {\tilde {\cal O}}\left( \tfrac{1}{\delta} + \tfrac{1}{\delta_1} + \tfrac{L_f}{\lambda+\kappa} + \tfrac{L}{n(\lambda+\kappa)} + \tfrac{\sqrt{1-\delta}(L_f + \lambda+\kappa)}{\delta (\lambda+\kappa)}  \right).
	$$
\end{lemma}

\subsection{Communication Complexity}\label{sec:cm-eclsvrg-ca}

In this subsection, we discuss the total communication cost by using EC-LSVRG + Catalyst. Same as the claim in the discussion of the communication complexity of ECLK, for simplicity, we choose $Q=Q_1$ in EC-LSVRG. 

Denote the communication cost of an vector in $\R^d$ as $U_d$ and the communication cost of  the compressed vector in $\R^d$ by using the compressor $Q$ as $r_Q$. From Theorem \ref{th:catalyst}, to achieve $P(x^k) - P^* \leq \epsilon$, the outer iteration number is ${ \tilde {\cal O}} \left( \frac{\sqrt{\lambda+\kappa}}{\sqrt{\lambda}} \log\frac{1}{\epsilon} \right)$, and from Lemma \ref{lm:eclsvrg-1}, the expected inner iteration number is 
\begin{align*}
	& \quad  {\tilde {\cal O}} \left(   \tfrac{1}{\delta} + \tfrac{L_f+L/n}{\lambda+\kappa} + \tfrac{\sqrt{(1-\delta)(L_f+\lambda+\kappa)({\bar L} + \lambda+\kappa)}}{\delta (\lambda+\kappa)}   + \tfrac{\sqrt{(1-\delta)(L_f+\lambda+\kappa)(L+\lambda+\kappa)}}{\sqrt{\delta}(\lambda+\kappa)} \right) \\ 
	& = {\tilde {\cal O}} \left(   \tfrac{1}{\delta} + \tfrac{ a_1 }{\lambda+\kappa} + \tfrac{b_1}{ \sqrt{\lambda+\kappa}}  \right), 
\end{align*}
where we denote $a_1 \eqdef L_f + \tfrac{L}{n} + \tfrac{\sqrt{1-\delta}(\sqrt{L_f {\bar L}} + \sqrt{\delta L_f L})}{\delta}$ and $b_1 \eqdef \tfrac{\sqrt{1-\delta} \left(  \sqrt{{\bar L}} + \sqrt{\delta L}  \right)}{\delta}$. 
Noticing that at each outer iteration, we need to communicate the uncompressed vector $h_{\tau, (k)}^0$, the expected total communication cost becomes 
\begin{align*}
	&  {\tilde {\cal O}} \left(  \left( \tfrac{\sqrt{\lambda+\kappa}}{\sqrt{\lambda}} \left(   \tfrac{1}{\delta} + \tfrac{ a_1 }{\lambda+\kappa} + \tfrac{b_1}{ \sqrt{\lambda+\kappa}}    \right) r_Q  +  \tfrac{\sqrt{\lambda+\kappa}}{\sqrt{\lambda}} U_d   \right)   \log \tfrac{1}{\epsilon}  \right) \\ 
	& = {\tilde {\cal O}} \left(   \tfrac{r_Q}{\sqrt{\lambda}} \log \tfrac{1}{\epsilon}  \left(  \left(  \tfrac{1}{\delta} + \tfrac{U_d}{r_Q}  \right) \sqrt{\lambda + \kappa}    + \tfrac{a_1}{\sqrt{\lambda+\kappa}} + b_1   \right) \right). 
\end{align*}

{\bf Optimal $\kappa$.} Since $\kappa\geq 0$ in Catalyst, it is easy to get the optimal $\kappa$ for minimizing the expected total communication cost. Let $\lambda_1 \eqdef a_1/\left(  \frac{1}{\delta} + \frac{U_d}{r_Q}  \right)$. If $\lambda \leq \lambda_1$, then the optimal $\kappa$ is $\lambda_1-\lambda$. If $\lambda>\lambda_1$, then the optimal $\kappa$ is $0$. Or equivalently, the optimal $\kappa = \max\{ \lambda_1, \lambda \} - \lambda$. 

Similarly, under the additional Assumption \ref{as:expcompressor}, from Theorem \ref{th:catalyst} and Lemma \ref{lm:eclsvrg-2}, the expected total communication cost is 
$$
{\tilde {\cal O}} \left(   \tfrac{r_Q}{\sqrt{\lambda}} \log \tfrac{1}{\epsilon}  \left(  \left(  \tfrac{1}{\delta} + \tfrac{U_d}{r_Q}  \right) \sqrt{\lambda + \kappa}    + \tfrac{a_2}{\sqrt{\lambda+\kappa}}  \right) \right),
$$
where $a_2 \eqdef L_f + \tfrac{L}{n} + \tfrac{\sqrt{1-\delta}L_f}{\delta}$. Let $\lambda_2 \eqdef a_2/\left(  \frac{1}{\delta} + \frac{U_d}{r_Q}  \right)$. Then the optimal $\kappa = \max\{ \lambda_2, \lambda \} - \lambda$. 

\noindent For TopK, if we use 64 bits for each element in $\R^d$,  
$
\tfrac{U_d}{r_Q} = \tfrac{64d}{(64+\log d)K} = \Theta \left(  \tfrac{d}{K \log d}  \right).
$
Even though the theoretical $\delta$ for TopK is $\nicefrac{K}{d}$, the actual value could be much larger than $\nicefrac{K}{d}$ in practice. Then $\nicefrac{U_d}{r_Q}$ may not be able to be bounded by ${\cal O}(\nicefrac{1}{\delta})$, and thus the communication complexity may be even worse than ECLK and ECSPDC.

\subsection{Remove the Dependence on $\nicefrac{U_d}{r_Q}$}
Due to the communication of uncompressed vectors at each outer iteration of the stratergies in Lemmas \ref{lm:eclsvrg-1} and \ref{lm:eclsvrg-2}, the expected total communication complexities depend on $\nicefrac{U_d}{r_Q}$, which may be much larger than $\nicefrac{1}{\delta}$. In this subsection, we show that we can actually remove the dependence on $\nicefrac{U_d}{r_Q}$ by communicating the compressed vector only. The initialization procedures and estimations of the expected inner iteration number are states in the following two lemmas. 

\begin{lemma}\label{lm:eclsvrg-1-re}
	Under Assumptions \ref{as:primal}, \ref{as:contracQQ1-ecSPDC}, and the premise of Theorem \ref{th:catalyst}, let us run EC-LSVRG to minimize $G_k$ and output $x^k \eqdef {x}_{(k)}^{T_k}$, $h_{\tau, (k)}^{T_k}$, and $e_{\tau, (k)}^{T_k}$, where 
	$
	T_k \eqdef \inf \{  K\geq 1, \Phi_{3, (k)}^K + G_k({x}_{(k)}^{K}) - G_k^* \leq \epsilon_k   \}
	$. 
	For the initialization of EC-LSVRG at the $k$-th outer iteration, we choose $p=\Theta(\delta_1)$, $x_{(k)}^0 = x^{k-1}$, $e_{\tau, (k)}^0 = 0$ or $e_{\tau, (k-1)}^{T_{k-1}}$, and $h_{\tau, (k)}^0 = h_{\tau, (k-1)}^{T_{k-1}}$ or $h_{\tau, (k-1)}^{T_{k-1}} + \kappa(y^{k-2}-y^{k-1})$ ($y^{-1}=y^0$). Then 
	\begin{align*}
		\mathbb{E}[T_k] & \leq {\tilde {\cal O}}\left( \tfrac{1}{\delta} + \tfrac{1}{\delta_1}  + \tfrac{\sqrt{(1-\delta)(L_f+\lambda+\kappa)({\bar L} + \lambda+\kappa)}}{\delta (\lambda+\kappa)}    +  \tfrac{L_f}{\lambda+\kappa} + \tfrac{L}{n(\lambda+\kappa)}  + \tfrac{\sqrt{(1-\delta)(L_f+\lambda+\kappa)(L+\lambda+\kappa)}}{\sqrt{\delta}(\lambda+\kappa)} \right),
	\end{align*} 
	where the notation ${\tilde {\cal O}}$ hides some universal constants and some logarithmic dependencies in $\delta$, $\delta_1$, $\lambda$, $\kappa$, $L_f$, and $N$. 
\end{lemma}

\begin{lemma}\label{lm:eclsvrg-2-re}
	Under Assumptions \ref{as:primal}, \ref{as:contracQQ1-ecSPDC}, \ref{as:expcompressor}, and the premise of Theorem \ref{th:catalyst}, let us run EC-LSVRG to minimize $G_k$ and output $x^k \eqdef {x}_{(k)}^{T_k}$, $h_{\tau, (k)}^{T_k}$, and $e_{\tau, (k)}^{T_k}$, where 
	$
	T_k \eqdef \inf \{  K\geq 1, \Phi_{4, (k)}^K + G_k({x}_{(k)}^{K}) - G_k^* \leq \epsilon_k   \}
	$. 
	Choose the initialization of EC-LSVRG at the $k$-th outer iteration be the same as that in Lemma \ref{lm:eclsvrg-1-re}. Then 
	$
	\mathbb{E}[T_k] \leq {\tilde {\cal O}}\left( \tfrac{1}{\delta} + \tfrac{1}{\delta_1} + \tfrac{L_f}{\lambda+\kappa} + \tfrac{L}{n(\lambda+\kappa)} + \tfrac{\sqrt{(1-\delta)}(L_f + \lambda+\kappa)}{\delta (\lambda+\kappa)}  \right).
	$
\end{lemma}

{\bf Communication Complexity.} Same as the analysis in Section \ref{sec:cm-eclsvrg-ca}, the expected total communication cost of EC-LSVRG + Catalyst with the output and initialization precedures in Lemmas \ref{lm:eclsvrg-1-re} and \ref{lm:eclsvrg-2-re} can be obtained by simply replacing $\nicefrac{U_d}{r_Q}$ with $0$. It is evident that the communication complexity depends on $\nicefrac{1}{\delta}$ only up to logarithmic terms. In particular, if $1-\delta = \Theta(1)$, $\delta \leq \min \{  \nicefrac{\bar L}{L}, \nicefrac{n^2L_f}{L} \}$ and $L_f \geq \lambda$, then an optimal $\kappa$ is $\sqrt{L_f {\bar L}} - \lambda$, and the corresponding communication complexity is ${\tilde {\cal O}} \left(  \tfrac{r_Q}{\delta} \sqrt{\tfrac{\bar L}{\lambda}} \log \tfrac{1}{\epsilon}  \right)$. If Assumption \ref{as:expcompressor} is further invoked, when $1-\delta = \Theta(1)$, $\delta \leq \nicefrac{nL_f}{L}$, and $L_f \geq \lambda$, an optimal $\kappa$ is $L_f-\lambda$, and the corresponding communication complexity is ${\tilde {\cal O}} \left(  \tfrac{r_Q}{\delta} \sqrt{\tfrac{L_f}{\lambda}} \log \tfrac{1}{\epsilon}  \right)$.

\section{EC-SDCA + CATALYST}

In this section, we consider Problem (\ref{primal-spdc}). Let $\xi \eqdef \frac{1}{\lambda}g$. Then $\xi$ is 1-strongly convex if $g$ is $\lambda$-strongly convex. We apply the catalyst to problem (\ref{primal-spdc}), and for the subproblem, we use the error-compensated SDCA (Algorithm \ref{alg:ec-sdca}) in \citep{ecsdca} to solve it. At the $k$-th outer iteration, we use EC-SDCA to minimize $G_k(x) \eqdef P(x) + \frac{\kappa}{2}\|x-y^{k-1}\|^2$, and we also use subscript $(k)$ and superscript $K$ to denote the variables at the $k$-th outer iteration and $K$-th inner iteration (for instance, $x^K_{(k)}$, $\alpha^K_{(k)}$, $e^K_{\tau, (k)}$, $e^K_{(k)}$, and $u^K_{(k)}$). 

To apply EC-SDCA at the $k$-th outer iteration in Algorithm \ref{alg:catalyst}, we need to initialize $\alpha_{i \tau, (k)}^0$. It is natural to use the values of $\alpha_{i\tau}$ in the last inner loop to initialize $\alpha_{i \tau, (k)}^0$, and this is indeed the case in \citep{accSDCA}, where the accelerated SDCA was studied. Then in order to initialize $u^0_{(k)} = \tfrac{1}{(\lambda+\kappa)N} \sum_{\tau=1}^n \sum_{i=1}^m A_{i\tau} \alpha_{i\tau, (k)}^0$, the uncompressed vector $A_{i\tau} \alpha_{i\tau, (k)}^0$ need to be communicated. We state the initialization procedures formally and estimate the expected inner iteration number in the next two lemmas.

\begin{lemma}\label{lm:ecsdca-1}
	Assume $\delta<1$. Under Assumptions \ref{as:contracQQ1-ecSPDC}, \ref{as:ecSPDC}, and the premise of Theorem \ref{th:catalyst},  let us run EC-SDCA (Algorithm \ref{alg:ec-sdca}) to minimize $G_k$ and output $(x^k, \alpha^k) \eqdef (x_{(k)}^{T_k+1}, \alpha_{(k)}^{T_k})$, where 
	$
	T_k \eqdef \inf \{  K\geq 1, \sqrt{4n+\delta mn} \Psi^K_{3,(k)} + 2(G_k(x_{(k)}^{K+1}) - G_k^*) \leq \epsilon_k   \}
	$. 
	For the initialization of EC-SDCA at the $k$-th iteration, we choose $\alpha_{(k)}^0 = \alpha^{k-1}$ ($\alpha^0=0$) and $e_{\tau, (k)}^0 = 0$. Then 
	$$
	\mathbb{E}[T_k] \leq {\tilde {\cal O}}\left( \tfrac{1}{\delta} + m + \tfrac{a_3}{\lambda+\kappa} + \tfrac{b_3}{\sqrt{\lambda+\kappa}} \right),
	$$ 
	where $a_3 \eqdef \frac{R_m^2}{n \gamma} + \frac{R^2}{\gamma} + \frac{\sqrt{1-\delta}R{\bar R}}{\delta \gamma} + \frac{\sqrt{1-\delta}RR_m}{\sqrt{\delta}\gamma}$, $b_3 \eqdef \frac{1}{\delta} \sqrt{\frac{(1-\delta)({\bar R}^2 + \delta R_m^2)}{\gamma}}$ and the notation ${\tilde {\cal O}}$ hides some universal constants and some logarithmic dependencies in $\delta$, $\lambda$, $\kappa$, $R$, and $N$. 
\end{lemma}

\begin{remark}
	In EC-SDCA, $\nicefrac{R^2}{\gamma} \geq \lambda+\kappa$ is assumed. However, by adding the term $ \frac{b_3}{\sqrt{\lambda+\kappa}} \log \frac{1}{\epsilon}$ to the iteration complexity, the assumption $\nicefrac{R^2}{\gamma} \geq \lambda+\kappa$ is no longer needed, which can be seen easily from the proof of Theorem 3.3 in \citep{ecsdca}. 
\end{remark}

If we further invoke Assumption \ref{as:expcompressor} on the compressors in EC-SDCA, we can get the following better result. The proof is similar to that of Lemma \ref{lm:ecsdca-1}, thus we omit it.

\begin{lemma}\label{lm:ecsdca-2}
	Assume $\delta<1$. Under Assumptions \ref{as:contracQQ1-ecSPDC}, \ref{as:ecSPDC}, \ref{as:expcompressor}, and the premise of Theorem \ref{th:catalyst}, let us run EC-SDCA to minimize $G_k$ and output $(x^k, \alpha^k) \eqdef (x_{(k)}^{T_k+1}, \alpha_{(k)}^{T_k})$, where 
	$
	T_k \eqdef \inf \{  K\geq 1, 3\sqrt{2+\delta m} \Psi^K_{4,(k)} + 2(G_k(x_{(k)}^{K+1}) - G_k^*) \leq \epsilon_k   \}
	$. 
	For the initialization of EC-SDCA at the $k$-th iteration, we choose $\alpha_{(k)}^0 = \alpha^{k-1}$ ($\alpha^0=0$) and $e_{\tau, (k)}^0 = 0$. Then $\mathbb{E}[T_k] \leq {\tilde {\cal O}}\left( \tfrac{1}{\delta} + m + \tfrac{a_4}{\lambda+\kappa} \right), $ where $a_4 \eqdef \frac{R_m^2}{n \gamma} + \frac{R^2}{\gamma} + \frac{\sqrt{1-\delta}R^2}{\delta \gamma} $. 
\end{lemma}

\subsection{Communication Complexity}\label{sec:cm-ecsdca-ca}

In this subsection, we discuss the total communication cost by using EC-SDCA + Catalyst. From Theorem \ref{th:catalyst}, to get $P(x^k) - P^* \leq \epsilon$, the outer iteration number is ${\tilde {\cal O}} \left( \frac{\sqrt{\lambda+\kappa}}{\sqrt{\lambda}} \log\frac{1}{\epsilon} \right)$, and from Lemma \ref{lm:ecsdca-1}, the expected inner iteration number is ${\tilde {\cal O}} \left(   \tfrac{1}{\delta} + m + \tfrac{a_3}{\lambda+\kappa}  + \tfrac{b_3}{\sqrt{\lambda+\kappa}}   \right)$. Noticing that at each outer iteration, we need to communicate the uncompressed vector to initialize $u^0_{(k)}$, the expected total communication cost is 
\begin{align*}
	& \quad  {\tilde {\cal O}} \left(  \left( \tfrac{\sqrt{\lambda+\kappa}}{\sqrt{\lambda}} \left(   \tfrac{1}{\delta} + m + \tfrac{a_3}{\lambda+\kappa}  + \tfrac{b_3}{\sqrt{\lambda+\kappa}}  \right) r_Q   +  \tfrac{\sqrt{\lambda+\kappa}}{\sqrt{\lambda}} U_d   \right)   \log \tfrac{1}{\epsilon}  \right) \\ 
	&= {\tilde {\cal O}} \left(   \tfrac{r_Q}{\sqrt{\lambda}} \log \tfrac{1}{\epsilon}  \left(  \left(  \tfrac{1+\delta m}{\delta}  + \tfrac{U_d}{r_Q}  \right) \sqrt{\lambda + \kappa}    + \tfrac{a_3}{\sqrt{\lambda+\kappa}} + b_3  \right) \right).
\end{align*}

{\bf Optimal $\kappa$.} Since $\lambda+\kappa \geq \lambda$, it is easy to obtain the optimal $\kappa$ for minimizing the expected total communication cost. Let $\lambda_3 \eqdef {a_3}/ ({\tfrac{1}{\delta} + m + \tfrac{U_d}{r_Q}})$. Then the optimal $\kappa$ is $\max\{\lambda, \lambda_3\} - \lambda$. 

Similarly, under the additional Assumption \ref{as:expcompressor}, from Theorem \ref{th:catalyst} and Lemma \ref{lm:ecsdca-2}, the expected total communication cost is 
$$
{\tilde {\cal O}} \left(   \tfrac{r_Q}{\sqrt{\lambda}} \log \tfrac{1}{\epsilon}  \left(  \left(  \tfrac{1}{\delta} + m  + \tfrac{U_d}{r_Q}  \right) \sqrt{\lambda + \kappa}    + \tfrac{a_4}{\sqrt{\lambda+\kappa}}   \right) \right).
$$
Let $\lambda_4 \eqdef {a_4}/ ({\tfrac{1}{\delta} + m + \tfrac{U_d}{r_Q}})$. Then the optimal $\kappa$ is $\max\{\lambda, \lambda_4\} - \lambda$. 

The term $\nicefrac{U_d}{r_Q}$ also shows up in the expected total communication cost of EC-SDCA + Catalyst. As we analyzed in Section \ref{sec:cm-eclsvrg-ca}, the presence of $\nicefrac{U_d}{r_Q}$ may make the communication complexity worse than ECLK and ECSPDC. In next subsection, we try to remove the dependence on $\nicefrac{U_d}{r_Q}$.

\subsection{Remove the Dependence on $\nicefrac{U_d}{r_Q}$}

As we can see from the analysis of the communication complexity, the term $\nicefrac{U_d}{r_Q}$ shows up because of the communication of uncompressed vectors. Hence, in order to remove the dependence on $\nicefrac{U_d}{r_Q}$, we need to find initialization procedures that do not need the communication of uncompressed vectors. Fortunately, by investigating the proofs of EC-SDCA, we find out that the relation $u^0_{(k)} = \tfrac{1}{(\lambda+\kappa)N} \sum_{\tau=1}^n \sum_{i=1}^m A_{i\tau} \alpha_{i\tau, (k)}^0$ in the initialization is not necessary, and the relation ${\tilde u}^K_{(k)} = \tfrac{1}{(\lambda+\kappa)N} \sum_{\tau=1}^n \sum_{i=1}^m A_{i\tau} \alpha_{i\tau, (k)}^K$ is actually essential in the proofs, and need to be maintained. This leads to the initialization procedures in the next two lemmas, and the communication of uncompressed vectors is actually not needed for the initialization at each outer iteration.

\begin{lemma}\label{lm:ecsdca-1-re}
	Assume $\delta<1$. Under Assumptions \ref{as:contracQQ1-ecSPDC}, \ref{as:ecSPDC}, and the premise of Theorem \ref{th:catalyst}, let us run EC-SDCA to minimize $G_k$ and output $x^k \eqdef x_{(k)}^{T_k+1}$, $\alpha^k \eqdef \alpha_{(k)}^{T_k}$, $u_{(k)}^{T_k}$, and $e_{\tau, (k)}^{T_k}$, where 
	$
	T_k \eqdef \inf \{  K\geq 1, \sqrt{4n+\delta mn} \Psi^K_{3,(k)} + 2(G_k(x_{(k)}^{K+1}) - G_k^*) \leq \epsilon_k   \}
	$. 
	For the initialization of EC-SDCA at the $k$-th iteration, we choose $\alpha_{(k)}^0 = \alpha^{k-1}$ ($\alpha^0=0$), $u_{(k)}^0 = u_{(k-1)}^{T_{k-1}}$ ($u_{(1)}^0=0$), and $e_{\tau, (k)}^0 = e_{\tau, (k-1)}^{T_{k-1}}$ ($e_{\tau, (1)}^0=0$). Then 
	$
	\mathbb{E}[T_k] \leq {\tilde {\cal O}}\left( \tfrac{1}{\delta} + m + \tfrac{a_3}{\lambda+\kappa} + \tfrac{b_3}{\sqrt{\lambda+\kappa}}  \right). 
	$
\end{lemma}


\begin{lemma}\label{lm:ecsdca-2-re}
	Assume $\delta<1$. Under Assumptions \ref{as:contracQQ1-ecSPDC}, \ref{as:ecSPDC}, \ref{as:expcompressor}, and the premise of Theorem \ref{th:catalyst}, let us run EC-SDCA to minimize $G_k$ and output $x^k \eqdef x_{(k)}^{T_k+1}$, $\alpha^k \eqdef \alpha_{(k)}^{T_k}$, $u_{(k)}^{T_k}$, and $e_{\tau, (k)}^{T_k}$, where 
	$
	T_k \eqdef \inf \{  K\geq 1, 3\sqrt{2+\delta m} \Psi^K_{4,(k)} + 2(G_k(x_{(k)}^{K+1}) - G_k^*) \leq \epsilon_k   \}
	$. 
	Choose the initialization of EC-SDCA at the $k$-th iteration be the same as that in Lemma \ref{lm:ecsdca-1-re}. Then 
	$
	\mathbb{E}[T_k] \leq {\tilde {\cal O}}\left( \tfrac{1}{\delta} + m + \tfrac{a_4}{\lambda+\kappa} \right). 
	$
\end{lemma}

{\bf Communication Complexity.} Same as the analysis in Section \ref{sec:cm-ecsdca-ca}, the expected total communication cost of EC-SDCA + Catalyst with the output and initialization precedures in Lemmas \ref{lm:ecsdca-1-re} and \ref{lm:ecsdca-2-re} can be obtained by simply replacing $\nicefrac{U_d}{r_Q}$ with $0$, and only depends on $\nicefrac{1}{\delta}$ up tp logarithmic terms. In particular, if $\delta \leq \nicefrac{1}{m}$ and $\nicefrac{R{\bar R}}{\gamma} \geq \lambda$, then an optimal $\kappa$ is $\nicefrac{R{\bar R}}{\gamma} - \lambda$, and the corresponding communication complexity is ${\tilde {\cal O}} \left(  \tfrac{r_Q}{\delta} \sqrt{\tfrac{{\bar R}^2}{\lambda \gamma}} \log \tfrac{1}{\epsilon}  \right)$. If Assumption \ref{as:expcompressor} is further invoked, when $\delta \leq \nicefrac{1}{m}$ and $\nicefrac{R^2}{\gamma} \geq \lambda$ an optimal $\kappa$ is $\nicefrac{R^2}{\gamma} - \lambda$, and the corresponding communication complexity is ${\tilde {\cal O}} \left(  \tfrac{r_Q}{\delta} \sqrt{\tfrac{R^2}{\lambda \gamma}} \log \tfrac{1}{\epsilon}  \right)$.

\begin{figure*}[ht]
	\vspace{.05in}
	\centerline{\begin{tabular}{ccc}
			\includegraphics[width=0.3\linewidth]{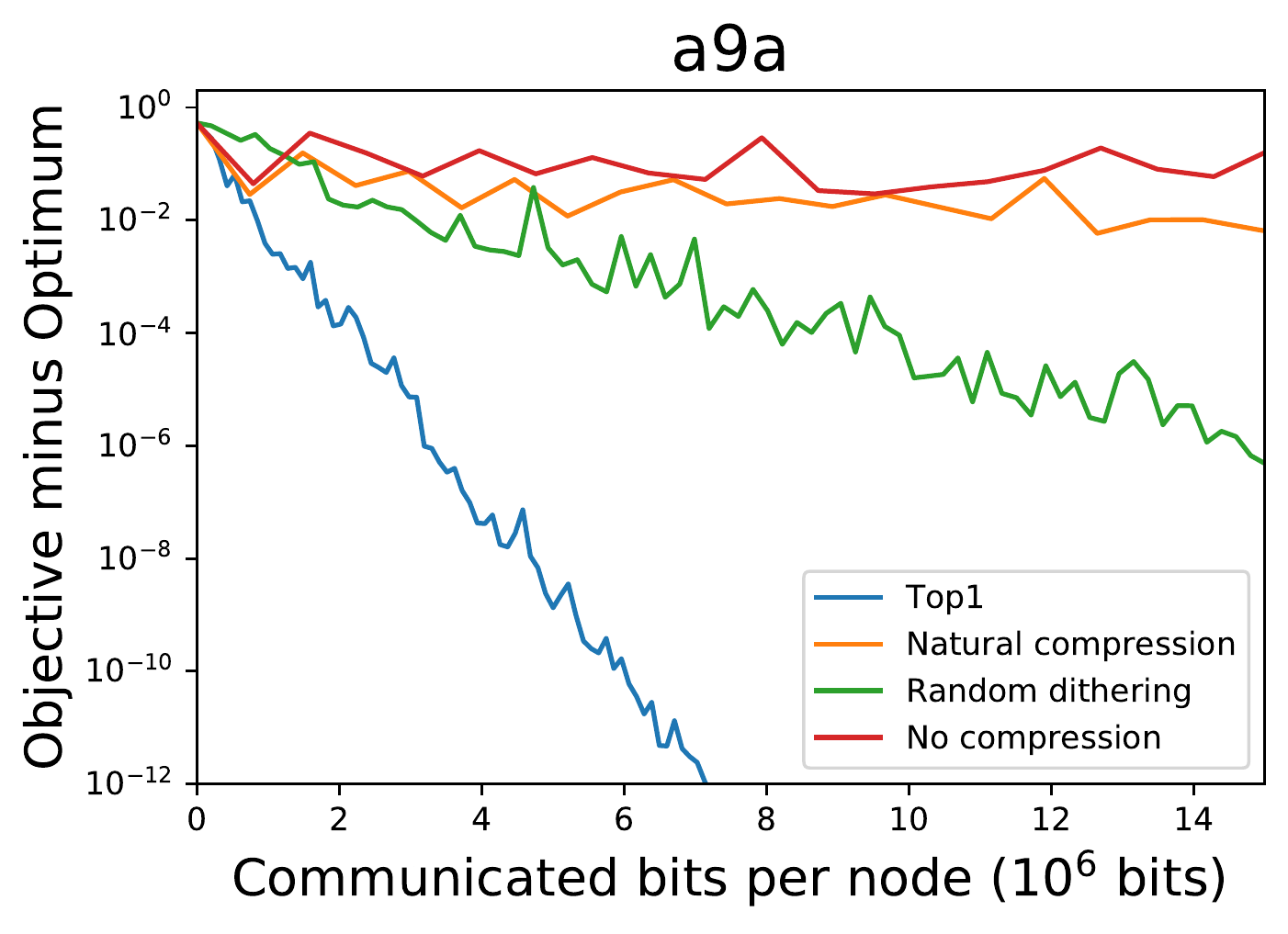}
			&\includegraphics[width=0.3\linewidth]{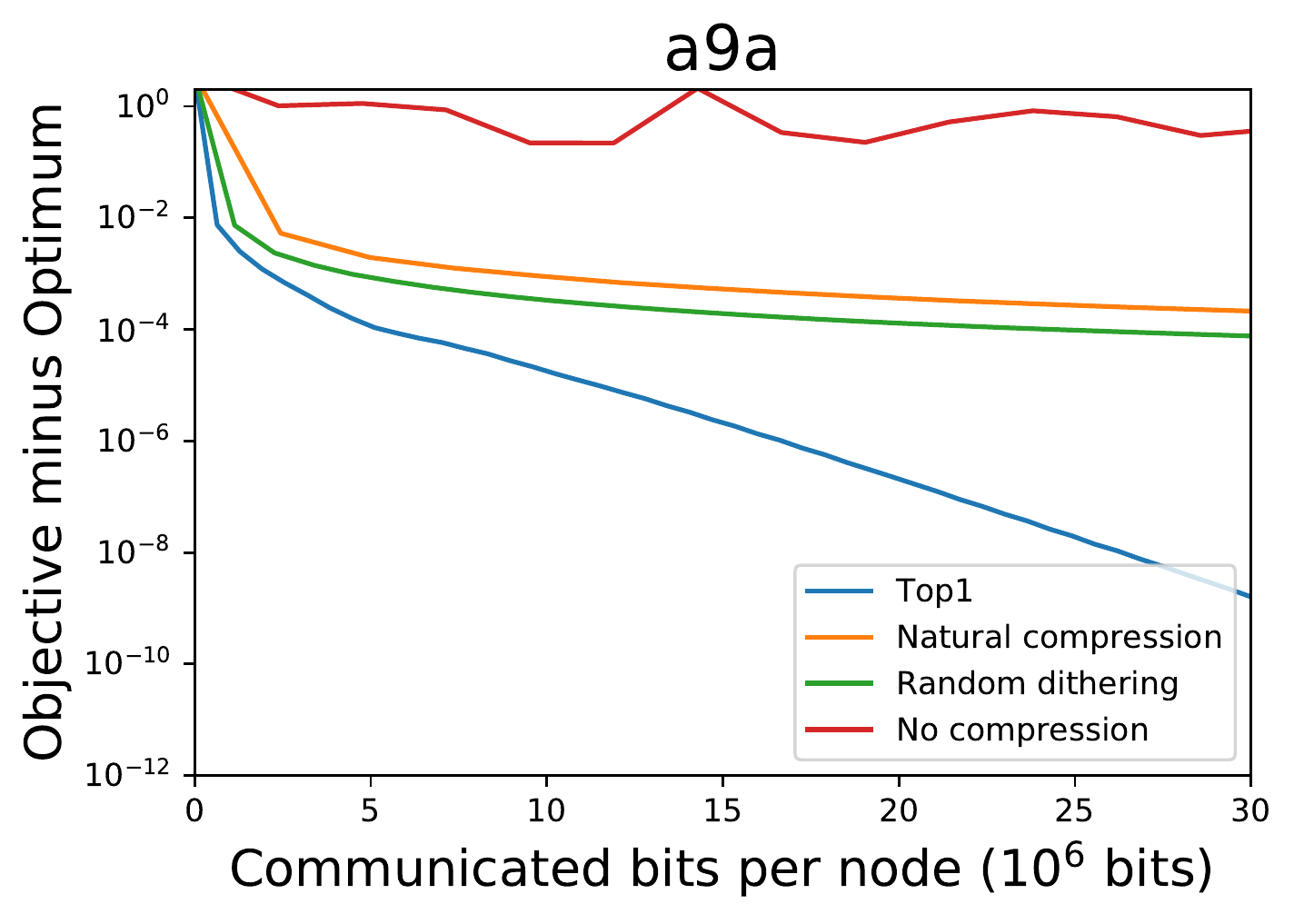}  
			&\includegraphics[width=0.3\linewidth]{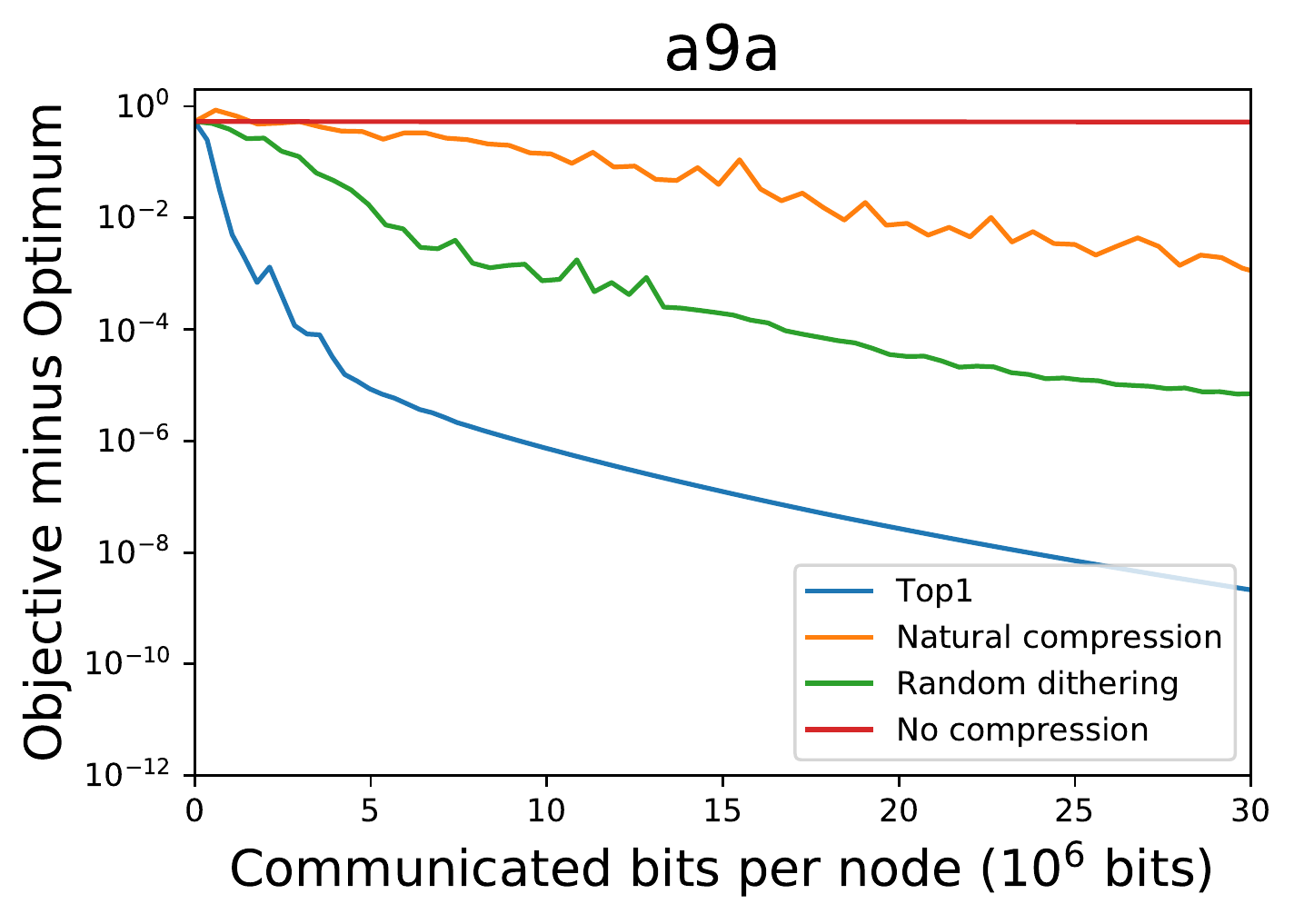} \\
		\end{tabular}
	}
	\vspace{.3in}
	\caption{The Communication Complexity Performance of ECSDCA-Catalyst, ECLSVRG-Catalyst, and ECSPDC Used with Compressors: Top1 VS Random Dithering VS Natural Compression VS No Compression on \dataset{a9a} Data Set}
	\label{fig:diff-a9a}
\end{figure*}

\begin{figure*}[ht]
	\vspace{.1in}
	\centerline{\begin{tabular}{ccc}
			\includegraphics[width=0.3\linewidth]{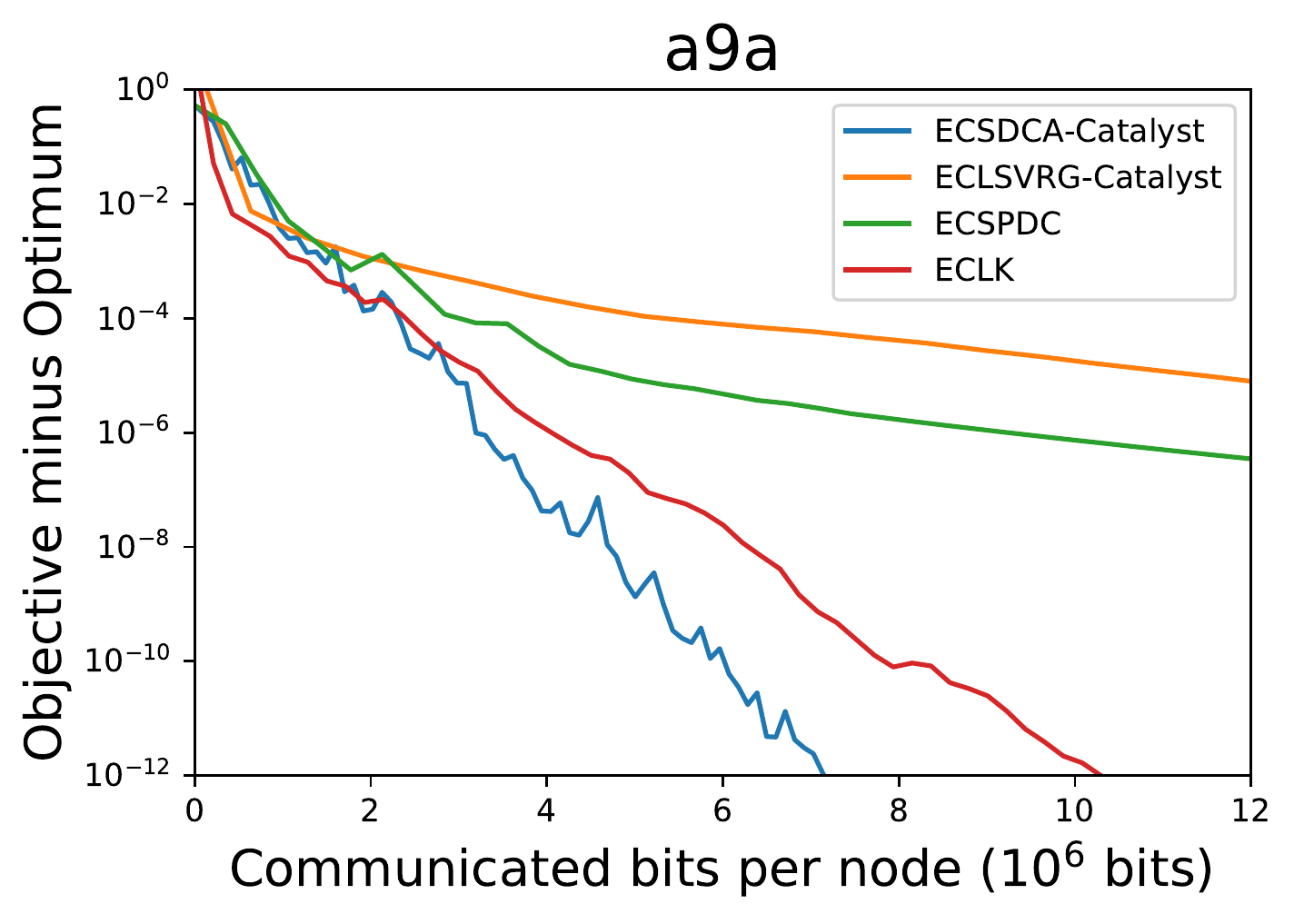}   
			&\includegraphics[width=0.3\linewidth]{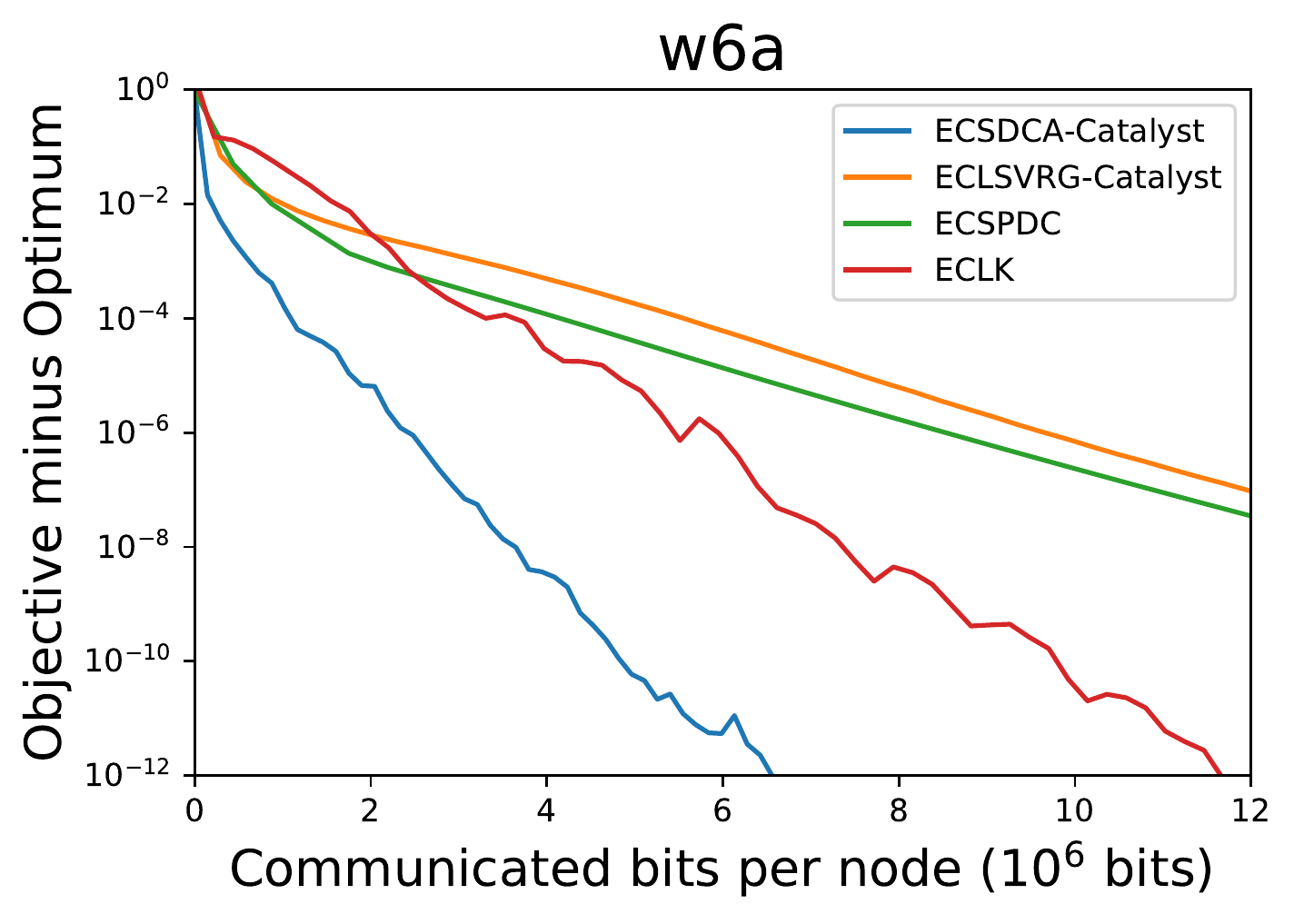}
			&	\includegraphics[width=0.3\linewidth]{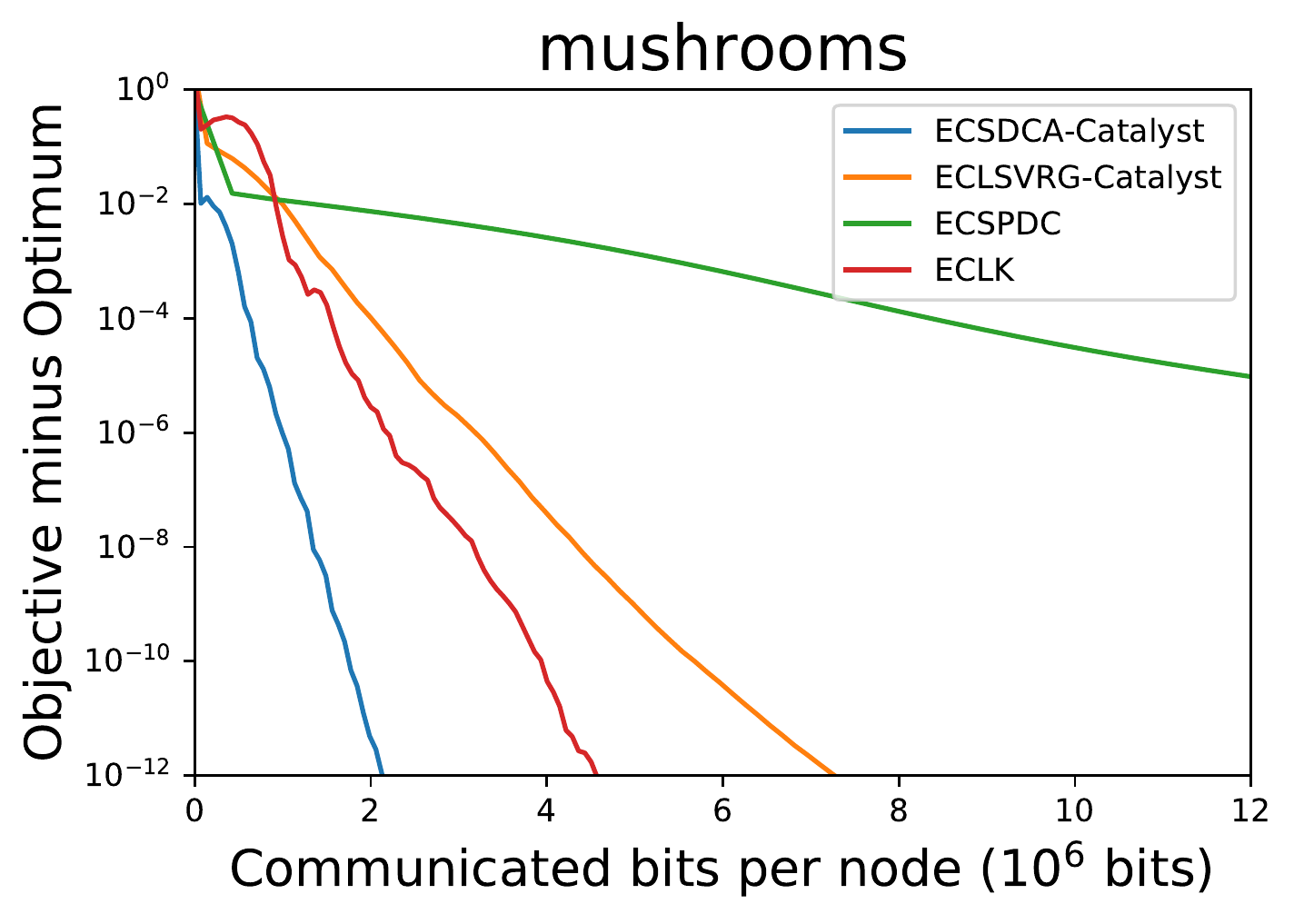}   
		\end{tabular}
	}
	\vspace{.3in}
	\caption{The Communication Complexity Performance of ECSDCA-Catalyst VS ECLSVRG-Catalyst VS ECSPDC VS ECLK for Top1 Compressor on \dataset{a9a}, \dataset{w6a}, and \dataset{mushrooms} Data Sets}
	\label{fig:acc}
\end{figure*}

\section{EXPERIMENTS}\label{sec:exp}

In this section, we implement our algorithms on the real world binary logistic regression tasks:
$$
x \mapsto \log \left(1+\exp(-y_i A_i^\top x)\right)+\tfrac{\lambda}{2}\Vert x\Vert^2,
$$
where $A_i,y_i$ are training sample pairs. We use the data sets: \dataset{a9a}, \dataset{w6a}, \dataset{phishing}, and \dataset{mushrooms} from LIBSVM Library \citep{chang2011libsvm}. More experiments can be found in the Appendix. 
\vskip 0.2mm
\noindent\textbf{Compressors.} In the experiments, we use Top1 and some contraction compressors transformed by unbiased ones such as random dithering ($s=\sqrt{d}$) and natural compression.  
\vskip 0.2mm
\noindent\textbf{Parameters.} We set $\lambda=1\times 10^{-5}$ and $n=20$. For all experiments, we use grid search to obtain the learning rate  $\{10^{-t},t = 0,1,2\cdots\}$. For ECSPDC, we use bisect method to obtain the argmax operator, $\theta$ is chosen by Theorem \ref{th:pf-SPDC}. For ECLSVRG and ECLK, we set $Q=Q_1$ and $p=\delta$. For Catalyst, we choose $\kappa$ by grid search $\{10^t\lambda: t\in \mathbb Z\}$. For the stopping criteria of the inner loop, a heuristic strategy was proposed for Catalyst in \citep{catalyst}, where the inner loop is constrained to perform at most $mn$ iterations. We employ this strategy similarly and the inner loop size is searched from $\{kd: k=1,2,5,10,100\}$, where $d$ is the dimension of data. 

\subsection{Effectiveness of TopK Compressor}

First, we demonstrate the effectiveness of TopK compressor compared with random dithering,  natural compression, and no compression. Figure \ref{fig:diff-a9a} shows that compression can improve the  performance with respect to the communication complexity in general, and TopK is specifically effective. 


\begin{figure*}[ht]
	\vspace{.1in}
	\centerline{\begin{tabular}{ccc}
			\includegraphics[width=0.3\linewidth]{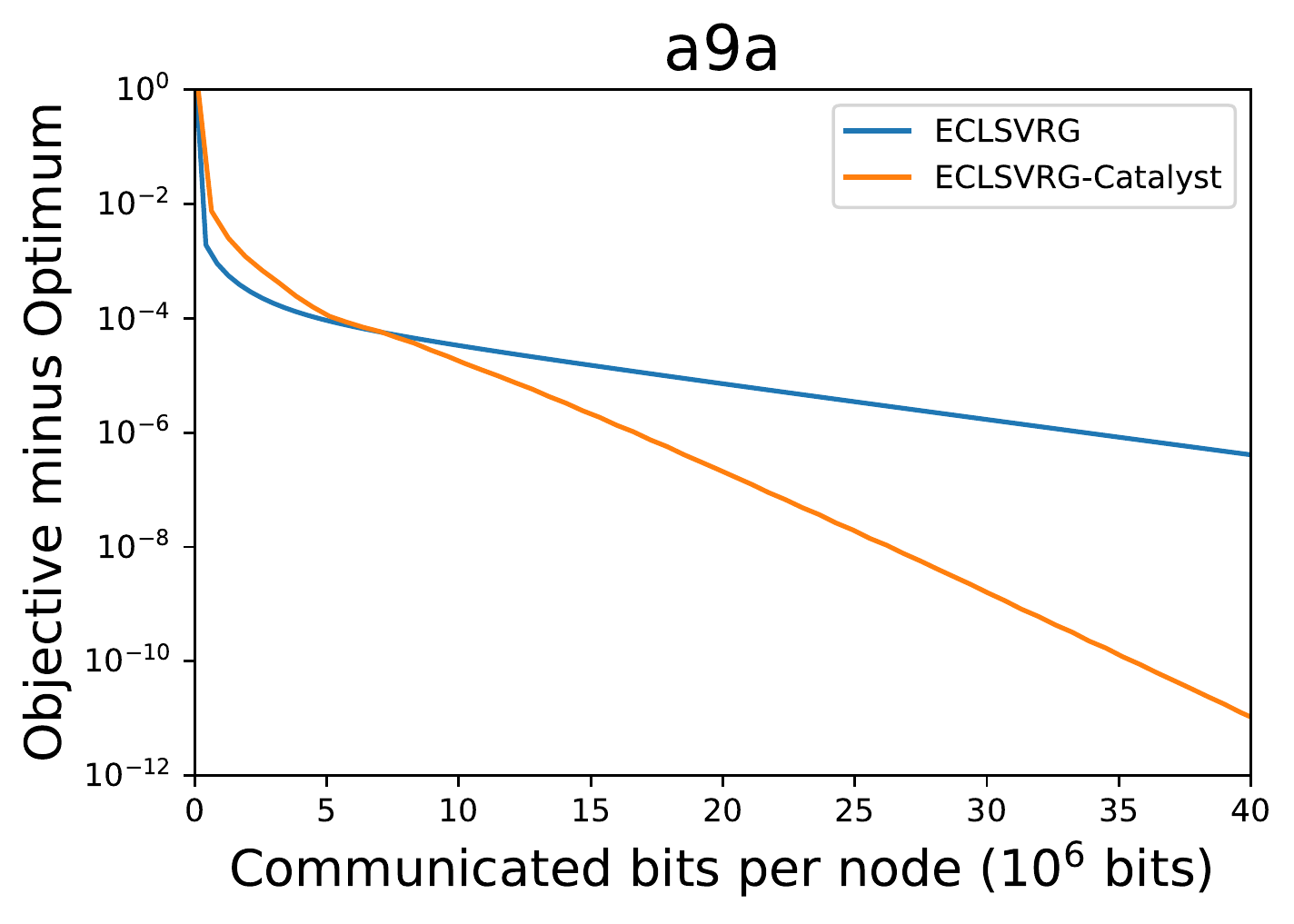}
			&\includegraphics[width=0.3\linewidth]{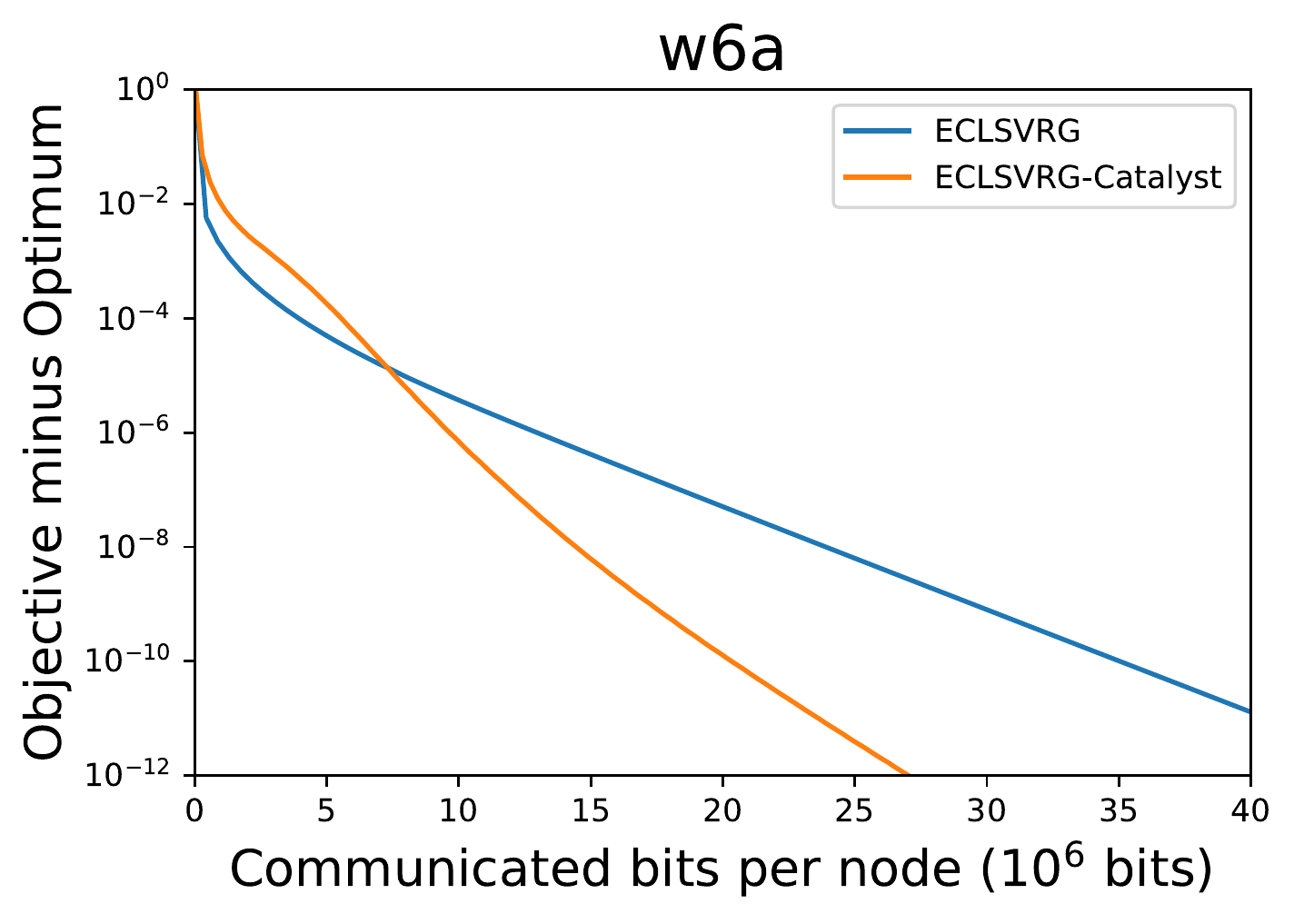}
			&\includegraphics[width=0.3\linewidth]{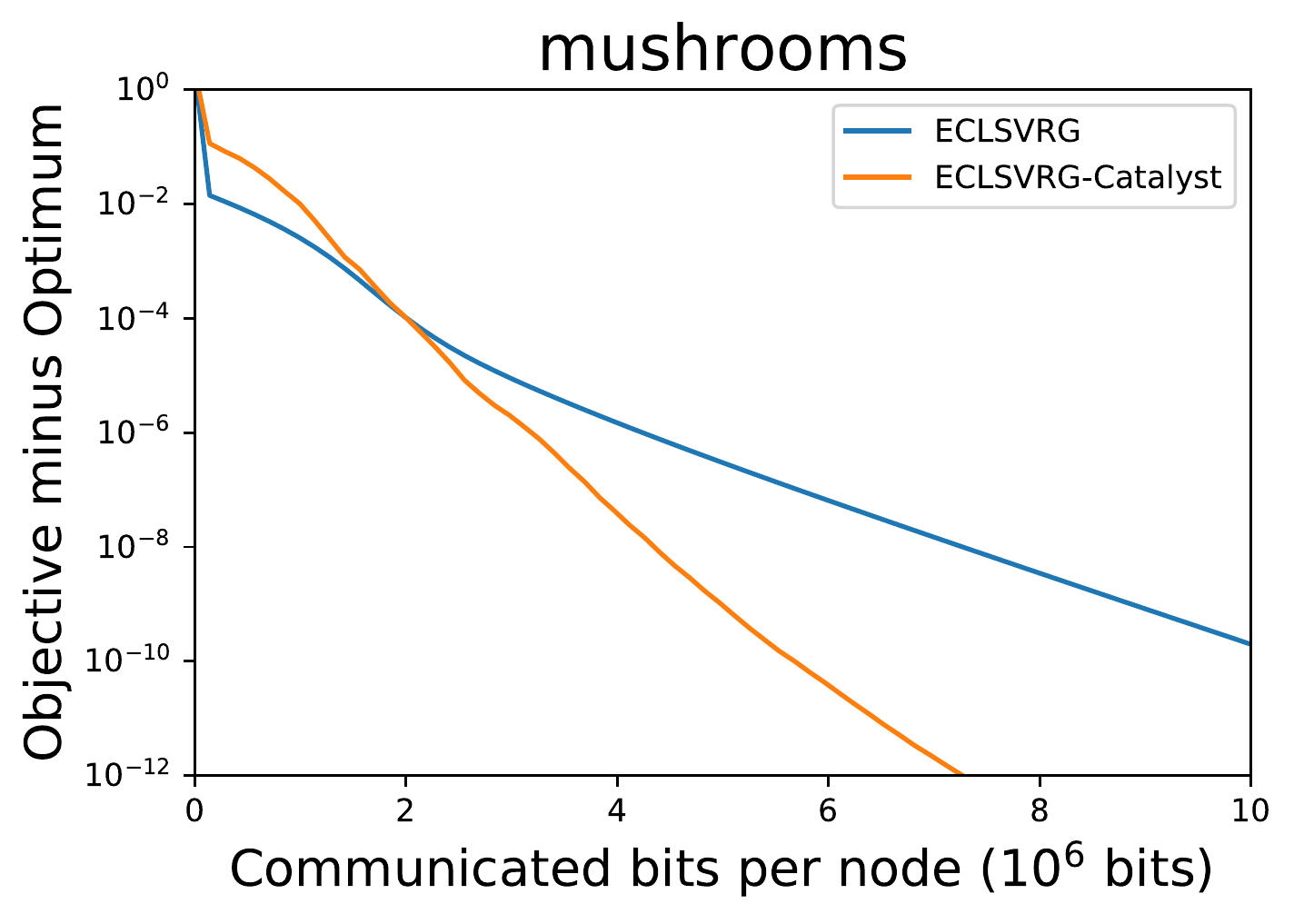}   
		\end{tabular}
	}
	\vspace{.3in}
	\caption{The Communication Complexity Performance of ECLSVRG VS ECLSVRG-Catalyst for Top1 Compressor on \dataset{a9a}, \dataset{w6a}, and \dataset{mushrooms} Data Sets}
	\label{fig:ECLSVRG-catalyst}
\end{figure*}

\begin{figure*}[ht]
	
	\vspace{.1in}
	
	\centerline{\begin{tabular}{ccc}
			\includegraphics[width=0.3\linewidth]{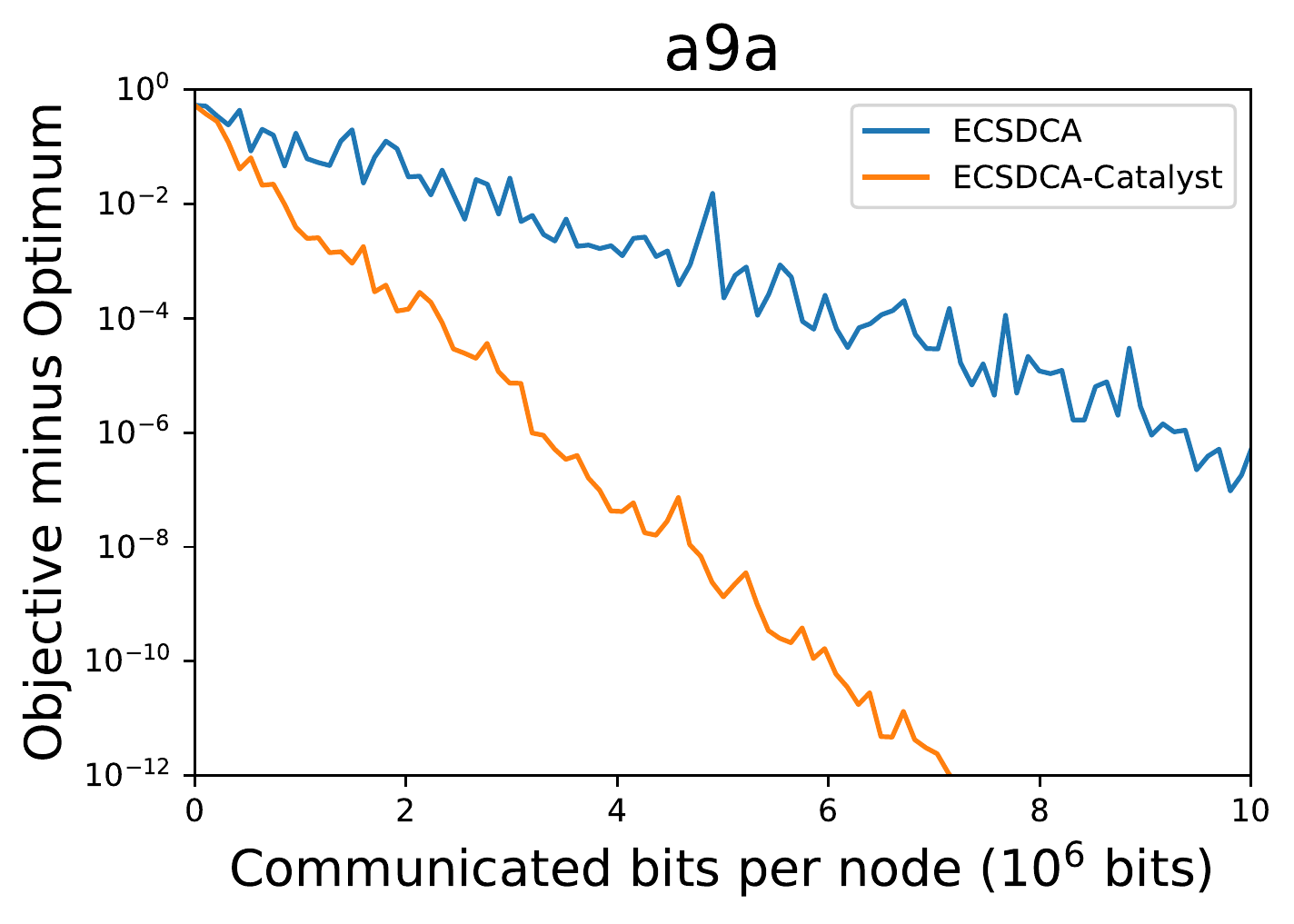}
			&\includegraphics[width=0.3\linewidth]{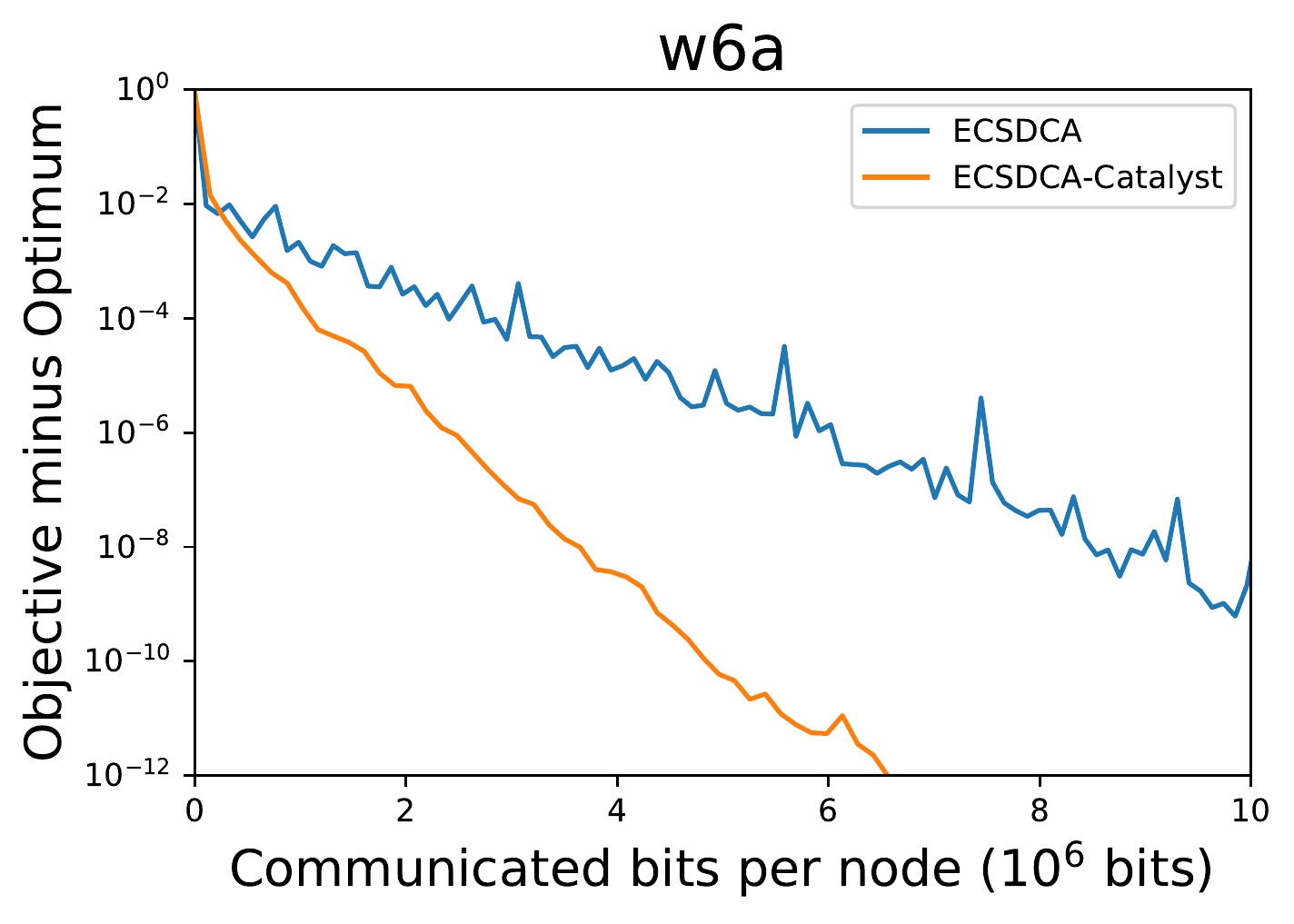}
			&\includegraphics[width=0.3\linewidth]{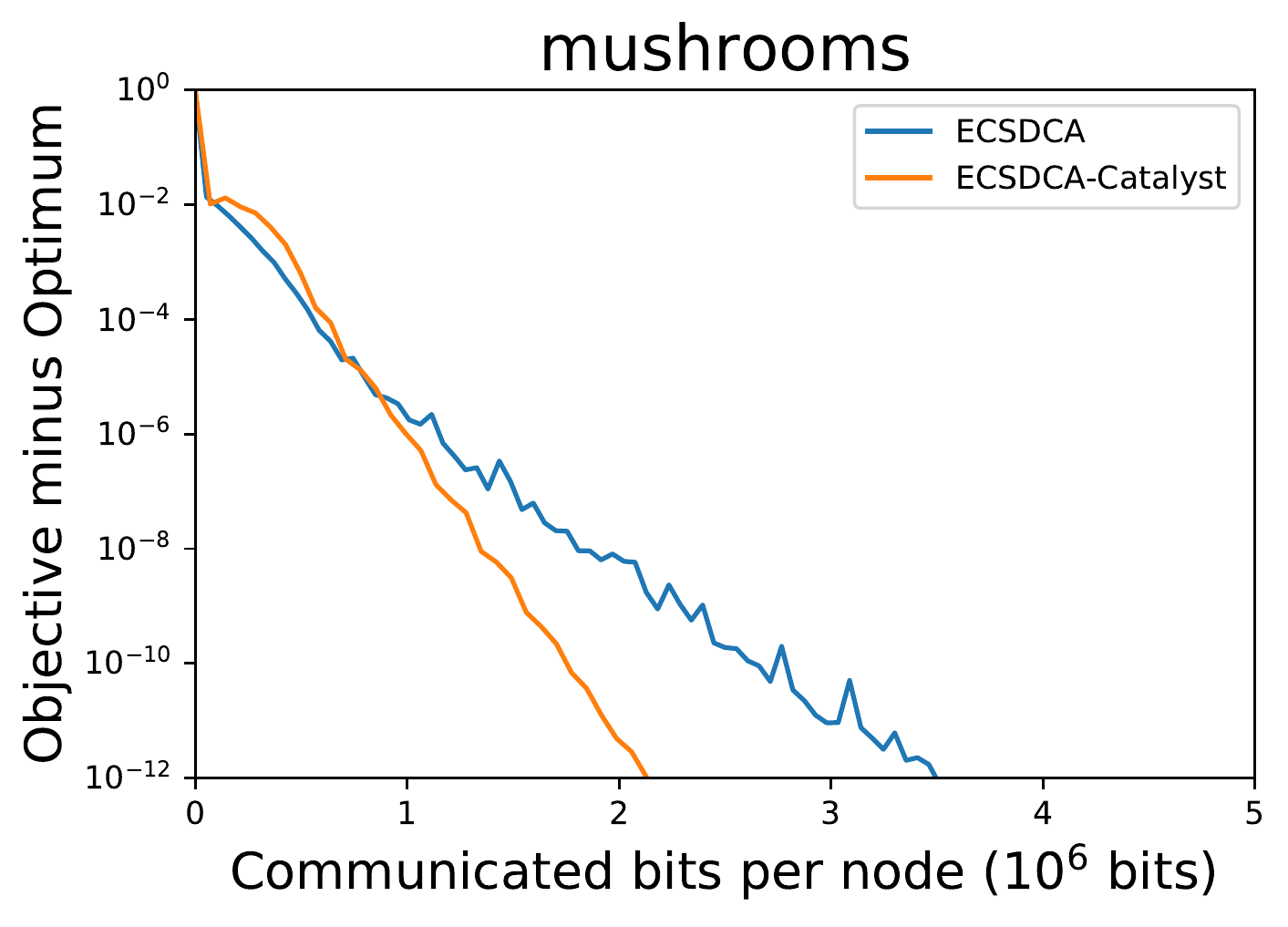}   
		\end{tabular}
	}
	\vspace{.3in}
	\caption{The Communication Complexity Performance of EC-SDCA VS ECSDCA-Catalyst for Top1 Compressor on \dataset{a9a}, \dataset{w6a}, and \dataset{mushrooms} Data Sets}
	\label{fig:ECSDCA-catalyst}
\end{figure*}

\subsection{Comparison of Different Accelerated Error Compensated Algorithms}

We compare Catalyst-based error-compensated algorithms and ECSPDC with ECLK, and also use the Top1 compressor. Figure \ref{fig:acc} shows that the performance of ECSDCA-Catalyst is the best for our tested data sets, which indicates the potential of the Catalyst-based error-compensated algorithm.

\subsection{Improvements from Catalyst Acceleration}

In this subsection, we compare Catalyst-based error-compensated algorithms with their baselines, namely, ECSDCA and ECLSVRG, where Top1 compressor is used. 
Figures \ref{fig:ECLSVRG-catalyst} and \ref{fig:ECSDCA-catalyst} show that Catalyst acceleration can indeed boost the speed of both ECSDCA and ECLSVRG with respect to the communication complexity significantly, which matches our theory.


\bibliography{references}
\bibliographystyle{plainnat}

\newpage
\appendix

\part*{Appendix}

\section{EXTRA EXPERIMENTS}

\subsection{Effectiveness of TopK Compressor}

\begin{figure}[ht]
	\vspace{.1in}
	\centerline{\begin{tabular}{ccc}
			\includegraphics[width=0.3\linewidth]{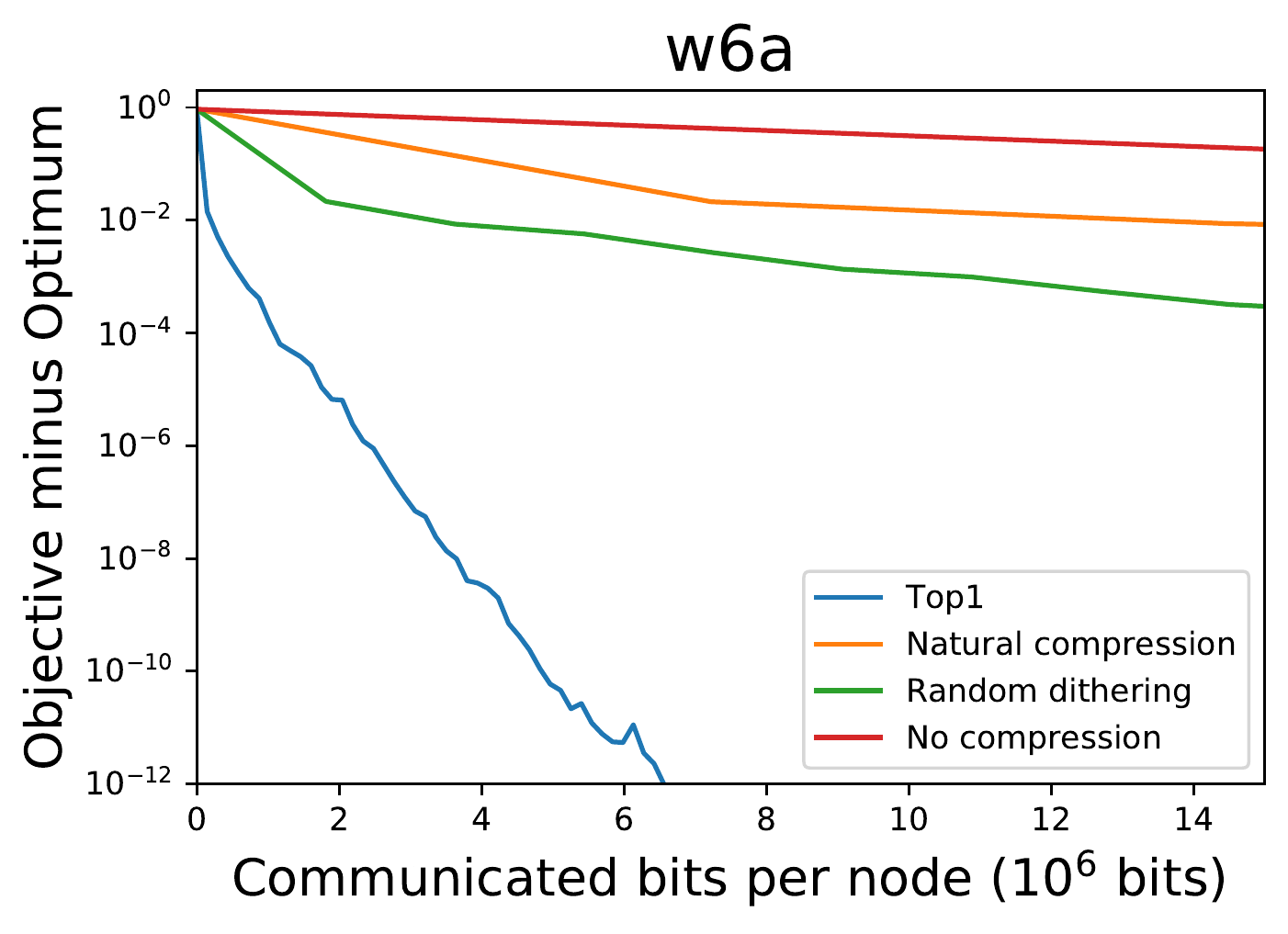}
			&\includegraphics[width=0.3\linewidth]{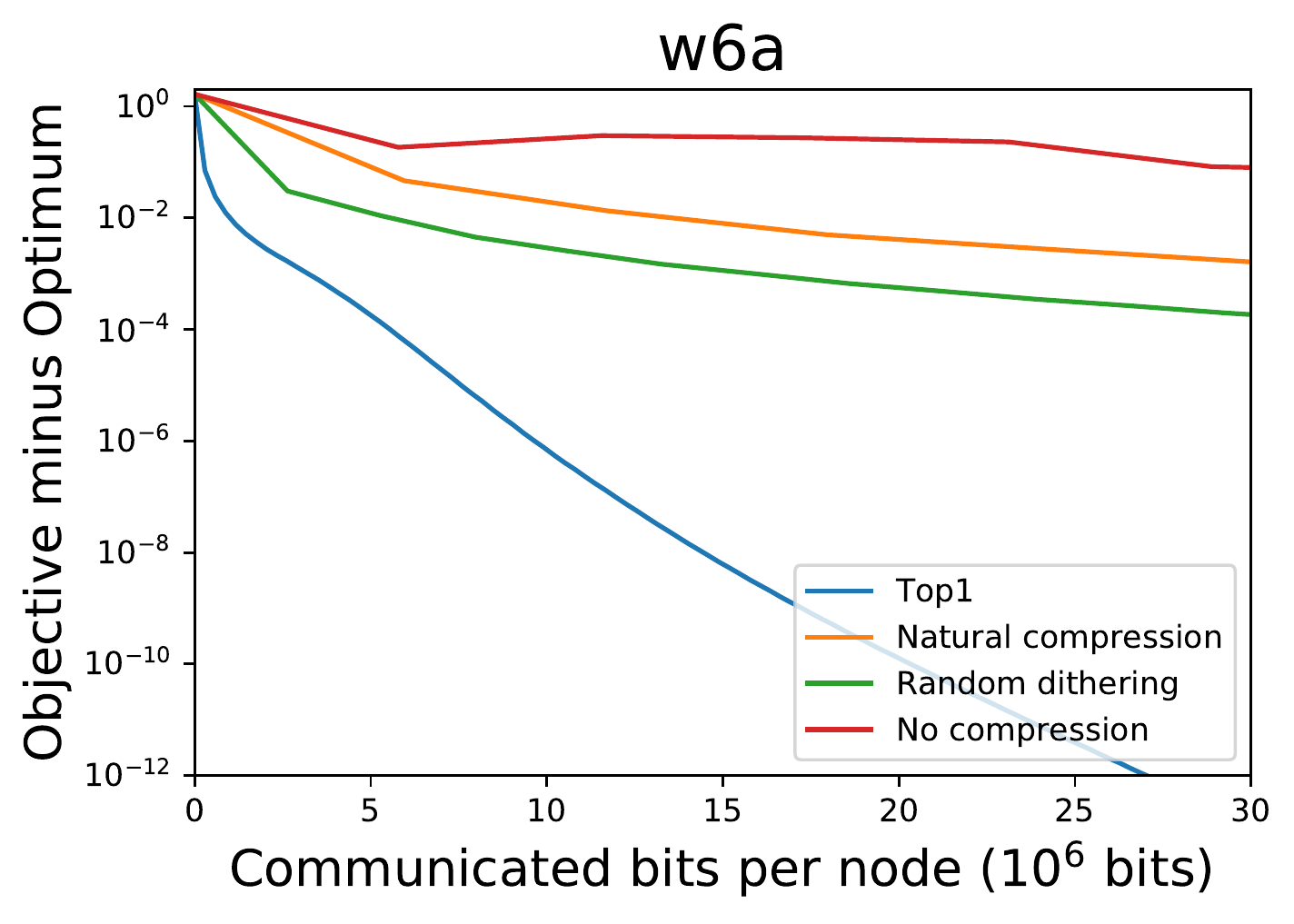}  
			&\includegraphics[width=0.3\linewidth]{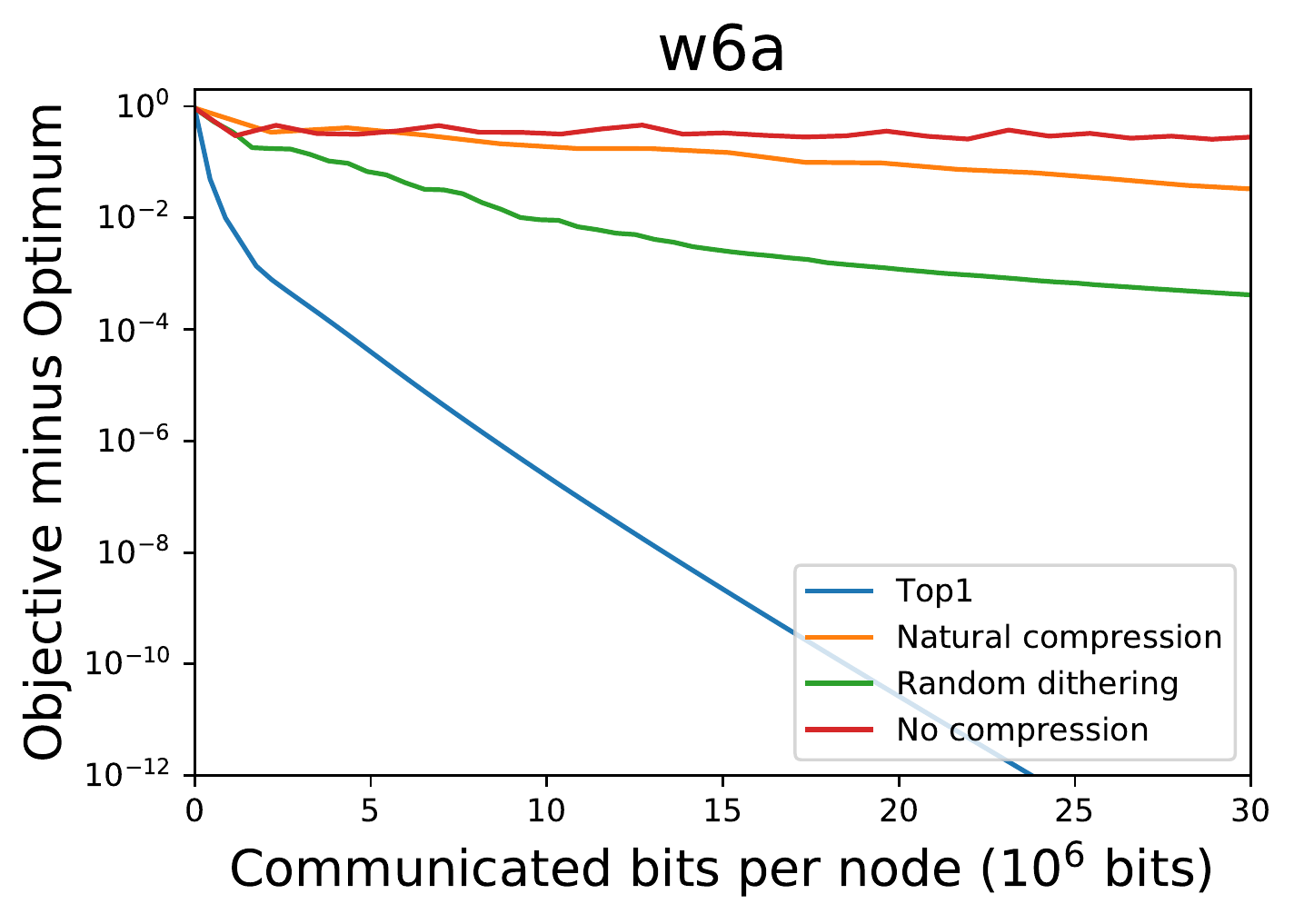} \\
		\end{tabular}
	}
	\vspace{.3in}
	\caption{The Communication Complexity Performance of ECSDCA-Catalyst, ECLSVRG-Catalyst, and ECSPDC Used with Compressors: Top1 VS Random Dithering VS Natural Compression VS No Compression on \dataset{w6a} Data Set}
	\label{fig:diff-w6a}
\end{figure}

\begin{figure}[ht]
	
	\vspace{.1in}
	\centerline{\begin{tabular}{ccc}
			\includegraphics[width=0.3\linewidth]{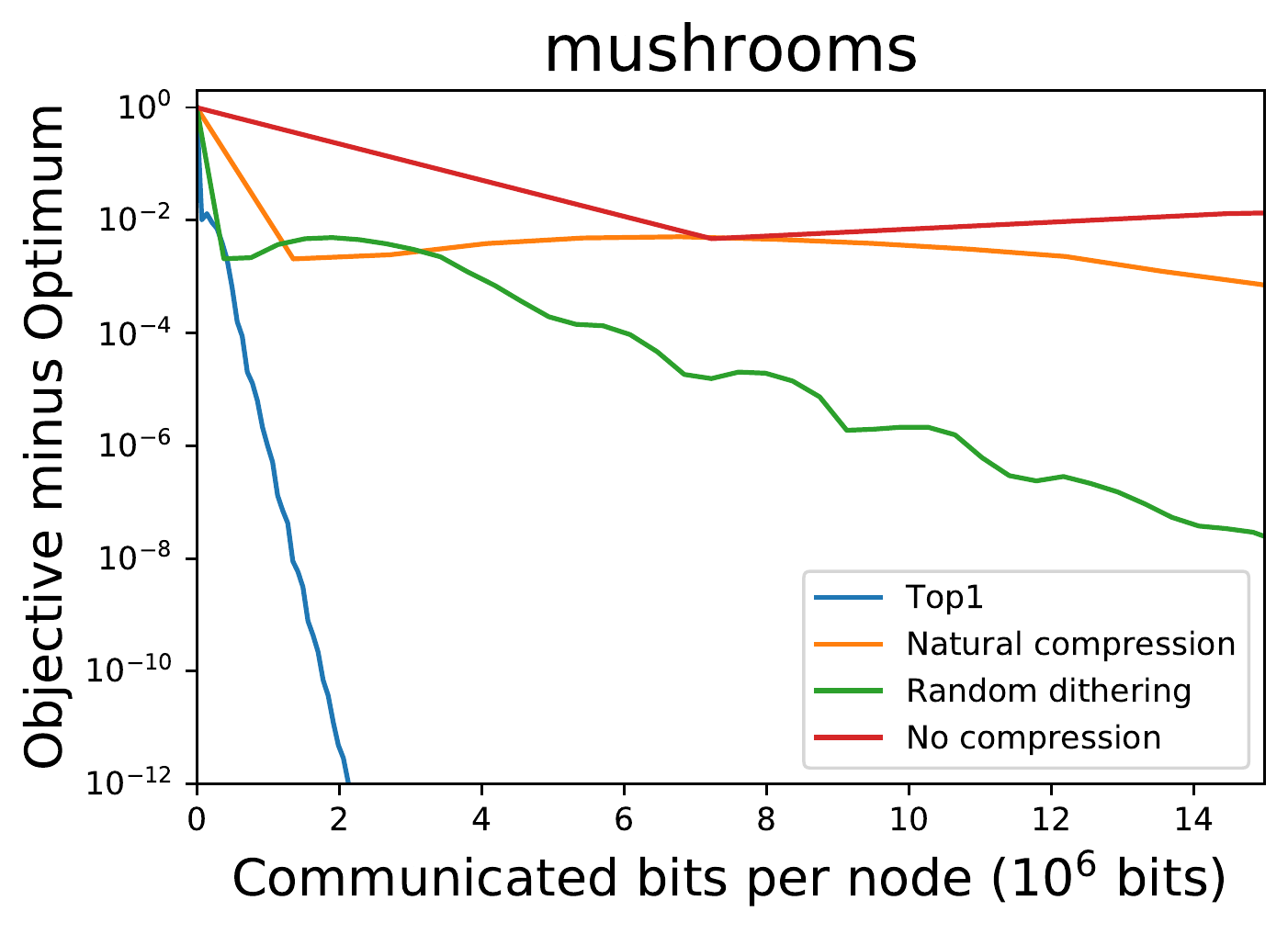}
			&\includegraphics[width=0.3\linewidth]{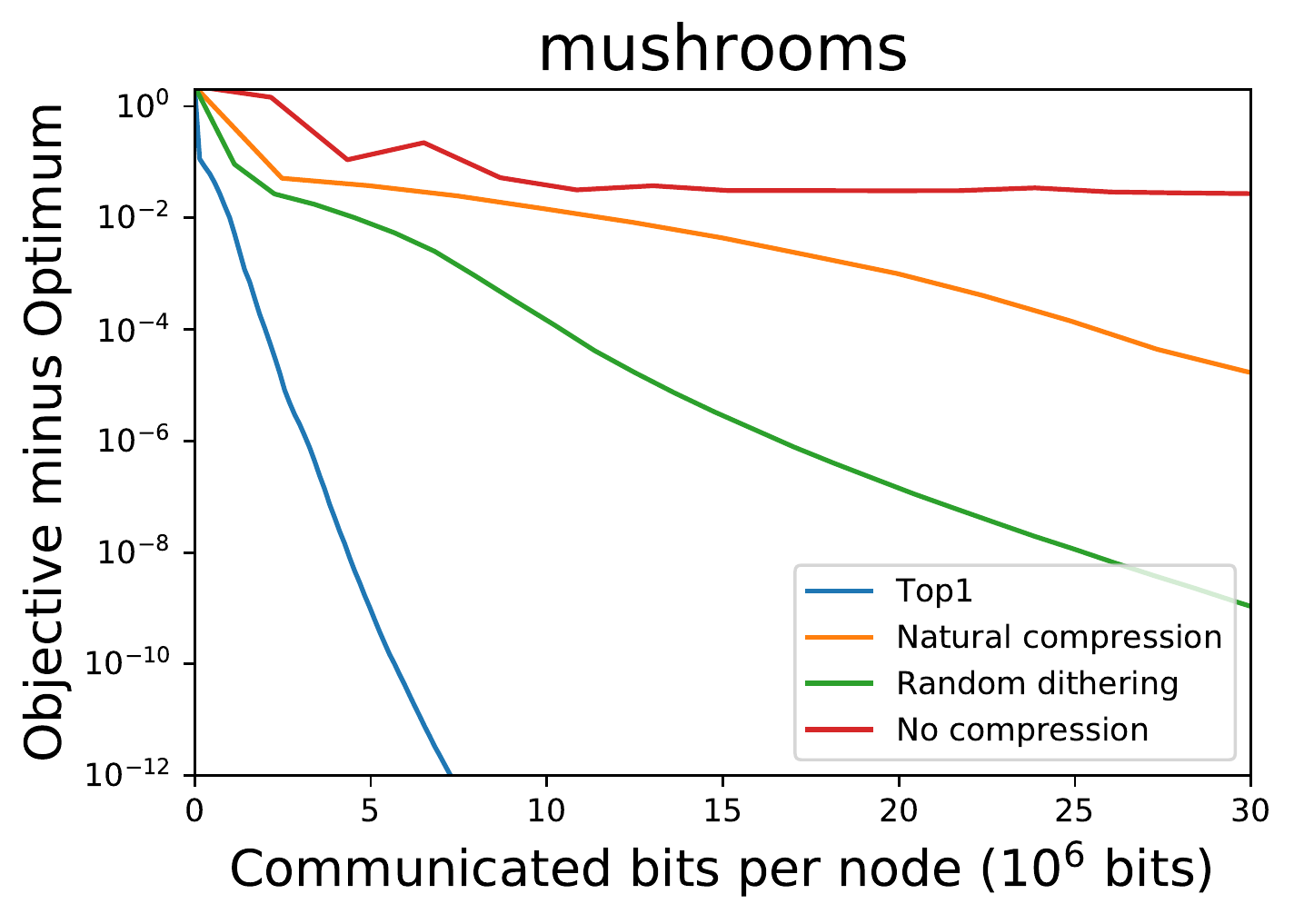}  
			&\includegraphics[width=0.3\linewidth]{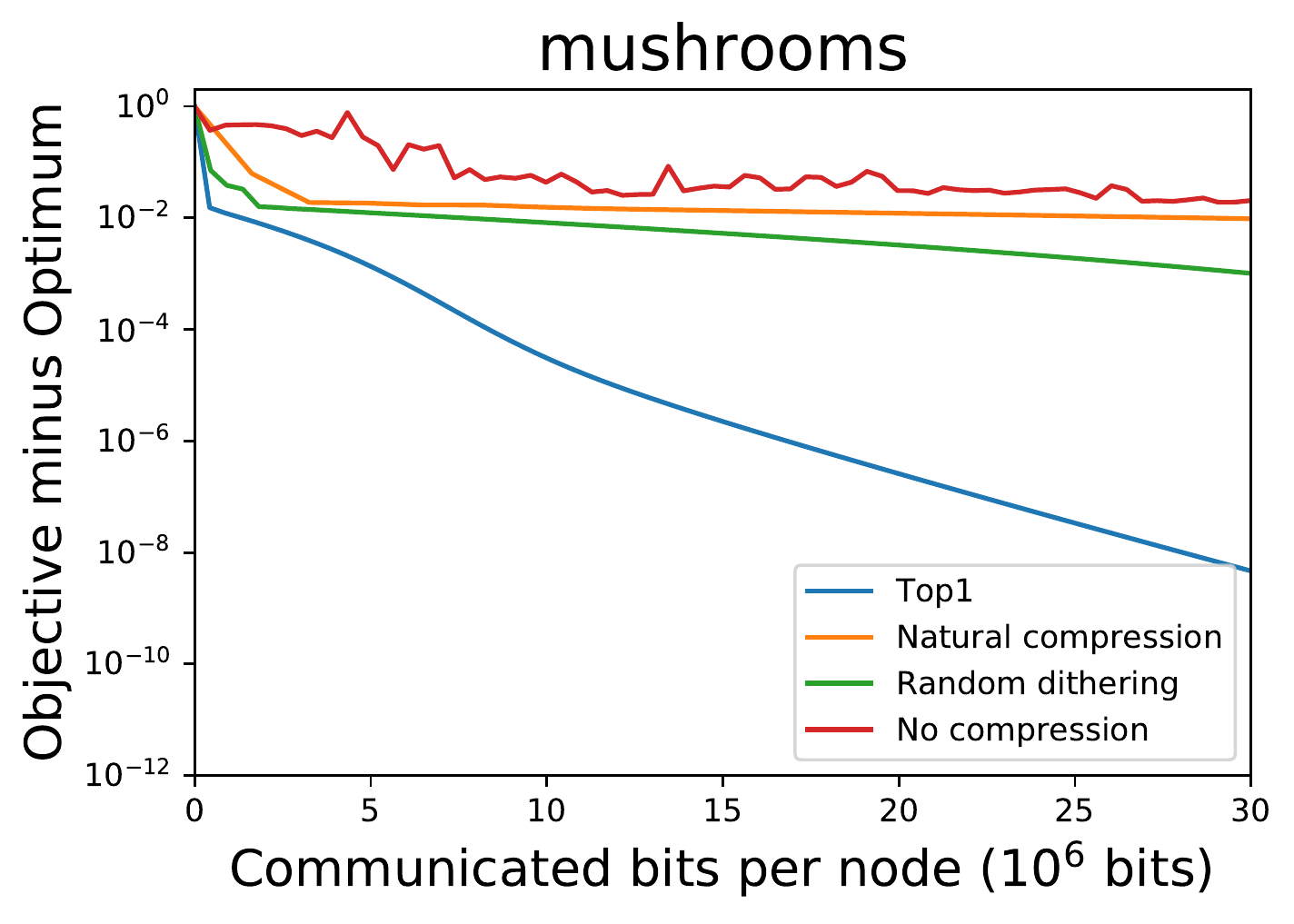} \\
		\end{tabular}
	}
	\vspace{.3in}
	\caption{The Communication Complexity Performance of ECSDCA-Catalyst, ECLSVRG-Catalyst, and ECSPDC Used with Compressors: Top1 VS Random Dithering VS Natural Compression VS No Compression on \dataset{mushrooms} Data Set}
	\label{fig:diff-mushrooms}
\end{figure}

\begin{figure}[ht]
	
	\vspace{.1in}
	\centerline{\begin{tabular}{ccc}
			\includegraphics[width=0.3\linewidth]{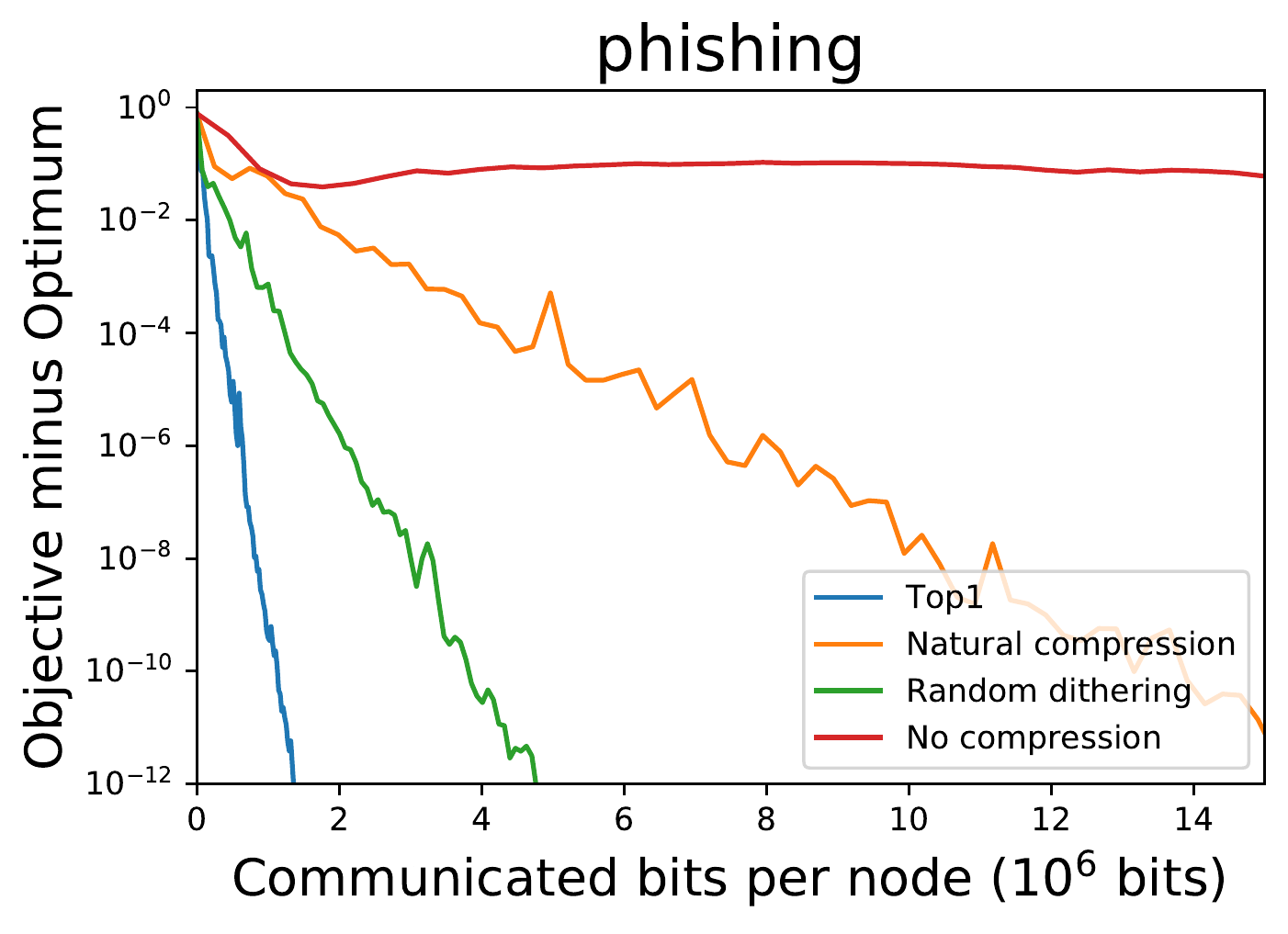}
			&\includegraphics[width=0.3\linewidth]{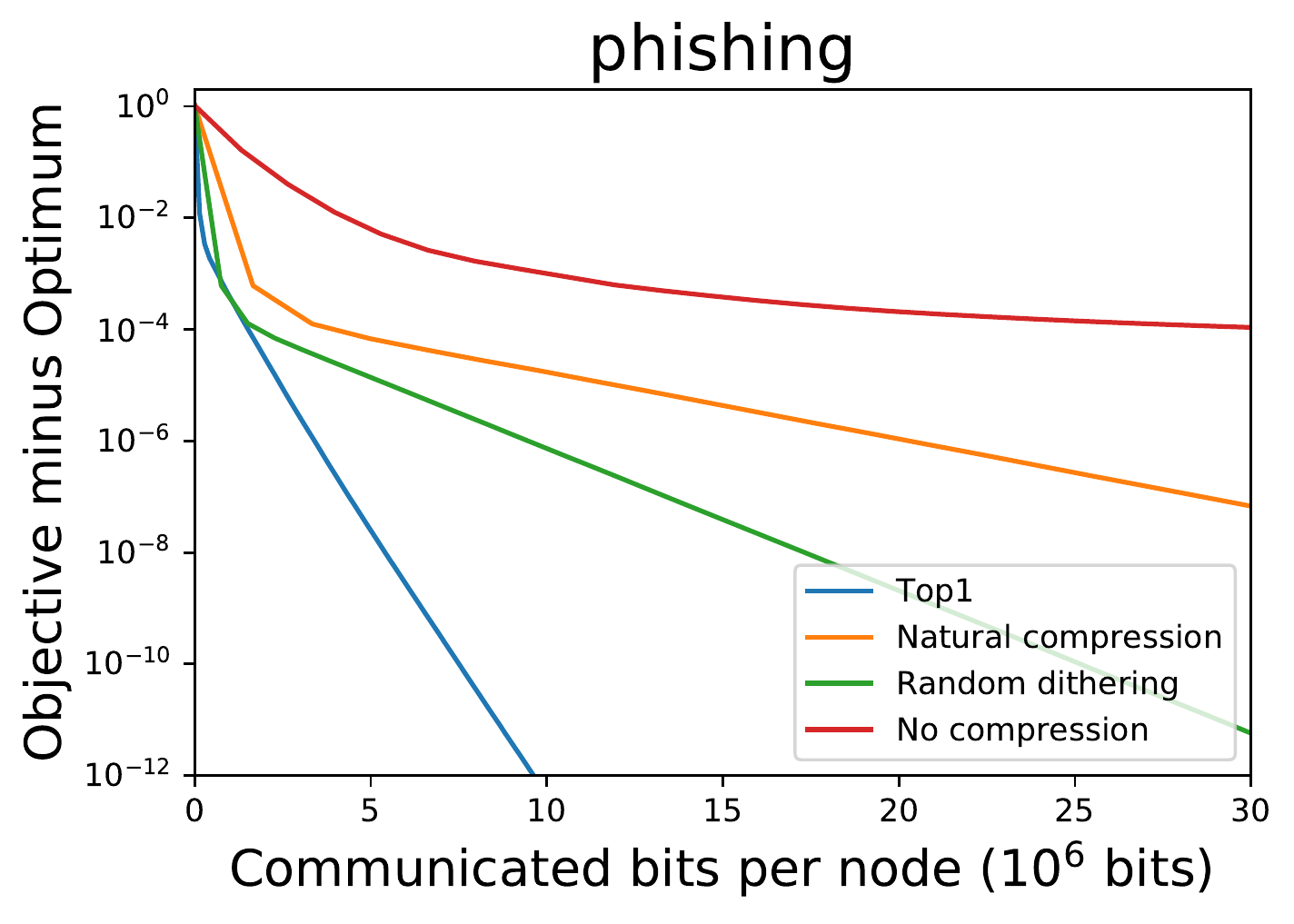}  
			&\includegraphics[width=0.3\linewidth]{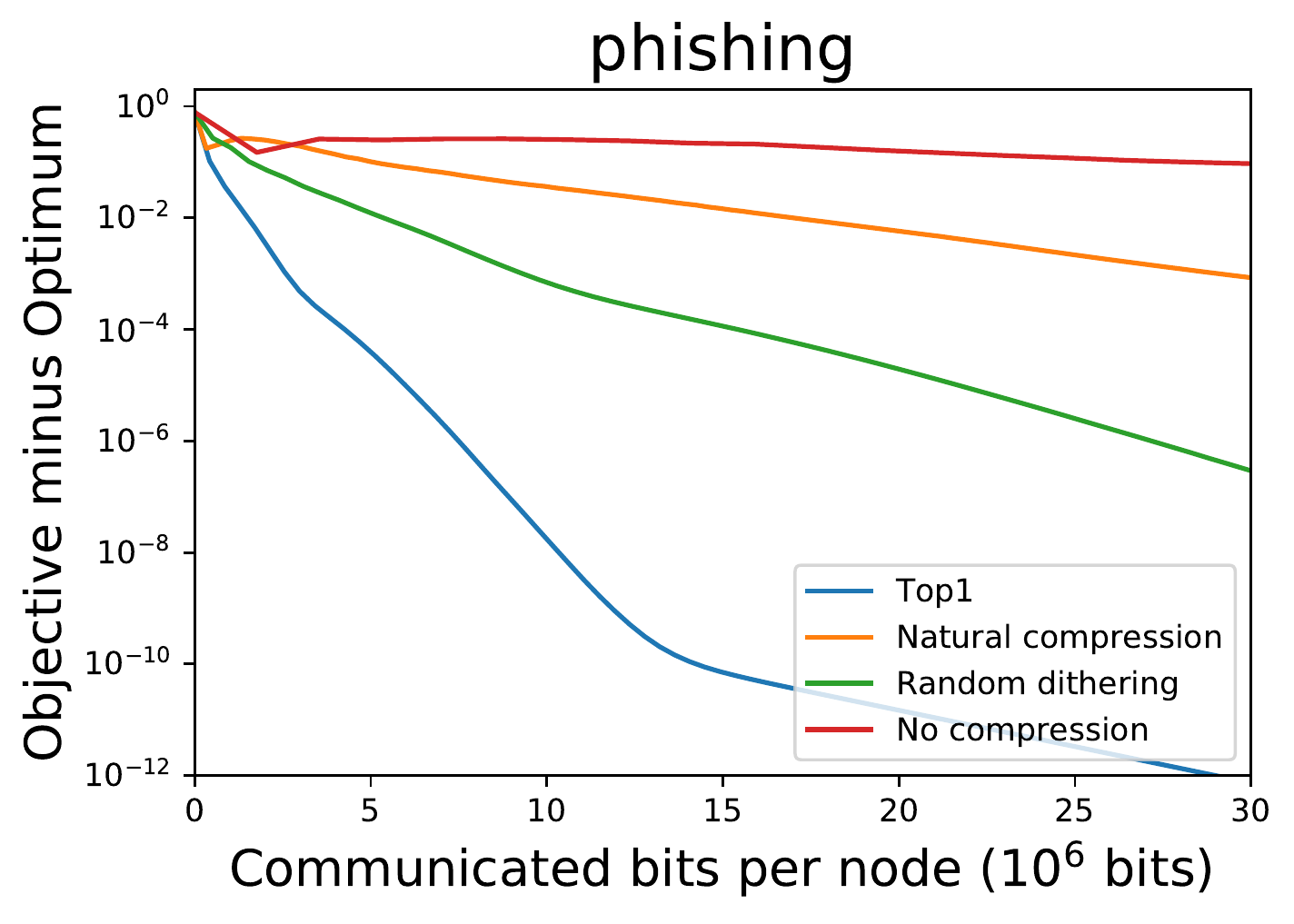} \\
		\end{tabular}
	}
	\vspace{.3in}
	\caption{The Communication Complexity Performance of ECSDCA-Catalyst, ECLSVRG-Catalyst, and ECSPDC Used with Compressors: Top1 VS Random Dithering VS Natural Compression VS No Compression on \dataset{phishing} Data Set.}
	\label{fig:diff-phishing}
\end{figure}

We demonstrate the effectiveness of TopK compressor compared with random dithering,  natural compression, and no compression. Figures \ref{fig:diff-w6a}, \ref{fig:diff-mushrooms}, and \ref{fig:diff-phishing} show that compression can improve the  performance with respect to the communication complexity in general, and TopK is specifically effective. 


\clearpage

\section{EC-LSVRG AND EC-SDCA ALGORITHMS}

In this section, we restate the two algorithms: EC-LSVRG and EC-SDCA in \citep{ecsdca}.

\begin{algorithm}[h!]
	\caption{Error compensated loopless SVRG (EC-LSVRG)}
	\label{alg:ec-lsvrg}
	\begin{algorithmic}
		\STATE {\bfseries Parameters:} stepsize $\eta >0$; probability $p \in (0, 1]$
		\STATE{\bfseries Initialization:}
		$x^0 = w^0 \in \R^d$; $e^0_\tau = 0 \in \R^d$; $u^0=1 \in \R$; $h^0_{\tau} \in \R^d$; $h^0 = \frac{1}{n} \sum_{\tau=1}^n h^0_\tau$
		\FOR{ $k = 0, 1, 2, \dots$} 
		\FOR{ $\tau = 1, \dots, n$} 
		\STATE Sample $i_k^\tau$ uniformly and independently in $[m]$ on each node 
		\STATE $g^k_{\tau} = \nabla f_{i_k^\tau}^{(\tau)}(x^k) - \nabla f_{i_k^\tau}^{(\tau)}(w^k) + \nabla f^{(\tau)}(w^k) - h^k_{\tau}$ 
		\STATE $y^k_{\tau} = Q(\eta g^k_{\tau} + e^k_{\tau})$, \quad $e^{k+1}_{\tau} = e^k_{\tau} + \eta g^k_{\tau} - y^k_{\tau}$
		\STATE $z^k_{\tau} = Q_1(\nabla f^{(\tau)}(w^k) - h^k_\tau)$, \quad $h^{k+1}_{\tau} = h^k_{\tau} + z^k_{\tau}$
		\STATE $u^{k+1}_\tau = 0$ for $\tau = 2, \dots, n$ 
		\STATE $
		u^{k+1}_1 = \left\{ \begin{array}{rl}
			1 & \mbox{ with probability $p$} \\
			0 &\mbox{ with probability $1-p$}
		\end{array} \right.
		$
		\STATE Send $y^k_{\tau}$, $z^k_{\tau}$, and $u^{k+1}_\tau$ to the other nodes 
		\STATE Receive $y^k_{\tau}$, $z^k_{\tau}$, and $u^{k+1}_\tau$  from the other nodes
		\STATE $y^k = \frac{1}{n} \sum_{\tau=1}^n y^k_{\tau}$, \quad $z^k = \frac{1}{n} \sum_{\tau=1}^n z^k_{\tau}$
		\STATE $u^{k+1} = \sum_{\tau=1}^n u^{k+1}_\tau$ 
		\STATE $x^{k+0.5} = x^k - (y^k + \eta h^k)$
		\STATE $x^{k+1} = \prox_{\eta \psi} (x^{k+0.5})$
		\STATE $
		w^{k+1} = \left\{ \begin{array}{rl}
			x^k & \mbox{ if $u^{k+1} = 1$} \\
			w^k &\mbox{ otherwise }
		\end{array} \right.
		$
		\STATE $h^{k+1} = h^k + z^k$
		\ENDFOR
		\ENDFOR
	\end{algorithmic}
\end{algorithm}

\begin{algorithm}[h!]
	\caption{Error compensated SDCA (EC-SDCA)}
	\label{alg:ec-sdca}
	\begin{algorithmic}
		\STATE {\bfseries Parameters:} $\theta>0$; $R_m \eqdef \max_{i, \tau} \|A_{i\tau}\|$; ${\bar R}^2 \eqdef \max_{\tau \in [n]} \{  \frac{1}{m}\lambda_{\rm max}(\sum_{i=1}^m A_{i\tau}A_{i\tau}^\top) \}$; $R^2 \eqdef \frac{1}{N} \lambda_{\rm max} (\sum_{\tau=1}^n \sum_{i=1}^m A_{i\tau} A_{i\tau}^\top)$; $p_{i\tau} = \frac{1}{m} \in \R$ for $i\in [m]$ and $\tau \in [n]$; positive constants $v_{i\tau} = R_m^2 + nR^2 \in \R$ for $i\in [m]$ and $\tau \in [n]$ 
		\STATE {\bfseries Initialization:}
		$\alpha^0 \in \R^{tN}$; $x^0 \in \R^d$; $u^0 = \frac{1}{\lambda N} \sum_{\tau=1}^n \sum_{i=1}^m A_{i\tau} \alpha^0_{i\tau} \in \R^d$; $e^0_\tau=0 \in \R^d$ for $\tau \in [n]$
		\FOR{ $k = 0, 1, 2, \dots$} 
		\FOR{ $\tau = 1, \dots, n$} 
		\STATE  $x^{k+1} = \nabla g^*(u^k)$ 
		\STATE $\alpha_{i\tau}^{k+1} = \alpha_{i\tau}^k$ for $i \in [m]$
		\STATE Sample $i_k^\tau$ uniformly and independently in $[m]$ on each node 
		\STATE$\Delta\alpha_{i^\tau_k \tau}^{k+1} = -\theta p_{i^\tau_k \tau}^{-1} \alpha_{i^\tau_k \tau}^{k} - \theta p_{i_k^\tau\tau}^{-1} \nabla \phi_{i^\tau_k \tau} (A_{i_k^\tau \tau}^\top x^{k+1}) $
		\STATE$\alpha_{i^\tau_{k} \tau}^{k+1} = \alpha_{i^\tau_{k} \tau}^{k} + \Delta\alpha_{i^\tau_{k} \tau}^{k+1}$
		\STATE $y_\tau^k = Q\left(  \frac{1}{\lambda m}A_{i_k^\tau \tau} \Delta\alpha_{i^\tau_{k} \tau}^{k+1} + e_\tau^k  \right)$
		\STATE $e_{\tau}^{k+1} = e_\tau^k + \frac{1}{\lambda m}A_{i_k^\tau \tau} \Delta\alpha_{i^\tau_{k} \tau}^{k+1} - y_\tau^k$
		\STATE Send $y^k_{\tau}$ to the other nodes 
		\STATE Receive $y^k_{\tau}$ from the other nodes
		\STATE$u^{k+1} = u^k + \frac{1}{n}\sum_{\tau=1}^n y_\tau^k$
		\ENDFOR
		\ENDFOR
	\end{algorithmic}
\end{algorithm}

\clearpage

\section{ERROR COMPENSATED SPDC}

In problem (\ref{primal-spdc}), we can replace each $\phi_{i\tau}(A_{i\tau}^\top x)$ by convex conjugation, i.e., 
$$
\phi_{i\tau}(A_{i\tau}^\top x) = \sup_{y\in \R^t} \{  \langle y, A_{i\tau}^\top x \rangle - \phi_{i\tau}^*(y)  \}, 
$$ 
where $\phi^*_{i\tau}$ is the conjugate function of $\phi_{i\tau}$. This leads to the following convex-concave saddle point problem 
\begin{align*}
	\min_{x\in \R^d} \max_{Y \in \R^{tN}}  f(x, Y)   \eqdef g(x)  +   \tfrac{1}{N} \sum_{\tau=1}^n \sum_{i=1}^m ( \langle y_{i\tau}, A_{i\tau}^\top x \rangle - \phi_{i\tau}^*(y_{i\tau})) ), 
\end{align*}
where 
$Y = (y_{11}^\top, ..., y_{m1}^\top, ..., y_{n1}^\top, ..., y_{mn}^\top)^\top \in \R^{tN}$ and $y_{i\tau} \in \R^t$.

\begin{algorithm}[h!]
	\caption{Error Compensated SPDC (ECSPDC)}
	\label{alg:ec-spdc}
	\begin{algorithmic}
		\STATE {\bfseries Parameters:} stepsize parameters $\sigma >0$; $\eta>0$; $\theta \in (0, 1)$ ; $R_m \eqdef \max_{i, \tau} \|A_{i\tau}\|$; ${\bar R}^2 \eqdef \max_{\tau \in [n]} \{  \tfrac{1}{m}\lambda_{\rm max}(\sum_{i=1}^m A_{i\tau}A_{i\tau}^\top) \}$; $R^2 \eqdef \tfrac{1}{N} \lambda_{\rm max} (\sum_{\tau=1}^n \sum_{i=1}^m A_{i\tau} A_{i\tau}^\top)$
		\STATE {\bfseries Initialization:}
		$x^0 = z^0 \in \R^d$; $e^0_\tau = 0 \in \R^d$; $y_{i\tau}^0 \in \R^t$; $u^0_\tau= \frac{1}{m}\sum_{i=1}^mA_{i\tau}y_{i\tau}^0$; $h_\tau^0\in \R^d$; $h^0 = \frac{1}{n} \sum_{\tau=1}^n h_\tau^0$
		\STATE {\bf for} {$k = 0, 1, 2, ...$} {\bf do}
		\STATE \quad  {\bf for} {$\tau = 1, ..., n$} {\bf do in parallel}
		\STATE \qquad Sample $i_k^\tau$ uniformly and independently in $[m]$ on each node 
		\STATE \qquad $
		y_{i\tau}^{k+1} = \left\{ \begin{array}{ll}
			\arg\max_{y\in \R^t}  \left\{  \langle y, A_{i\tau}^\top z^k \rangle - \phi^*_{i\tau} (y) - \frac{1}{2\sigma}\|y - y_{i\tau}^k\|^2  \right\}  & \mbox{ if $i = i_k^\tau $ } \\
			y_{i\tau}^k &\mbox{ if $i\neq i_k^\tau$}
		\end{array} \right.
		$
		\STATE \qquad $\Delta_\tau^k = Q(A_{i_k^\tau \tau} (y_{i_k^\tau \tau}^{k+1} - y_{i_k^\tau \tau}^k) + u_\tau^k - h_\tau^k + e_\tau^k)$
		\STATE \qquad $u_\tau^{k+1} = u_\tau^k + \frac{1}{m} A_{i_k^\tau \tau} (y_{i_k^\tau \tau}^{k+1} - y_{i_k^\tau \tau}^k)$, \quad $h_\tau^{k+1} = h_\tau^k + Q_1(u_\tau^k - h_\tau^k)$ 
		\STATE\qquad $e_\tau^{k+1} = e_\tau^k + A_{i_k^\tau \tau} (y_{i_k^\tau \tau}^{k+1} - y_{i_k^\tau \tau}^k) + u_\tau^k - h_\tau^k - \Delta_\tau^k$
		\STATE \qquad Send $\Delta_\tau^k$ and $Q_1(u_\tau^k - h^k_\tau )$ to the other nodes 
		\STATE \qquad Receive $\Delta_\tau^k$ and $Q_1(u_\tau^k - h^k_\tau )$ from the other nodes
		\STATE \qquad $\Delta^k = \frac{1}{n} \sum_{\tau=1}^n \Delta_\tau^k$ 
		\STATE \qquad $x^{k+1} = \arg\min_{x\in \R^d} \left\{  g(x) + \langle h^k + \Delta^k, x \rangle +  \tfrac{\|x-x^k\|^2}{2\eta}  \right\} $
		\STATE \qquad $h^{k+1} = h^k + \frac{1}{n} \sum_{\tau=1}^n Q_1(u_\tau^k- h^k_\tau )$, \quad $z^{k+1} = x^{k+1} + \theta(x^{k+1} - x^k)$
		\STATE \quad {\bf end for}
		\STATE {\bf end for}
	\end{algorithmic}
\end{algorithm}

\noindent {\bf Description of error-compensated SPDC (Algorithm \ref{alg:ec-spdc}).} In distributed SPDC, the search direction at the $k$-th iteration is 
$$
\frac{1}{n} \sum_{\tau=1}^n \left(  \frac{1}{m} \sum_{i=1}^m A_{i\tau}y_{i\tau}^k  +  A_{i_k^\tau \tau}(y_{i_k^\tau \tau}^{k+1} - y_{i_k^\tau \tau}^k)  \right), 
$$
where $i_k^\tau$ is sampled uniformly and independently in $[m]\eqdef \{1, 2, ..., m\}$ on each node. When $y_{i\tau}$ goes to the optimal solution, the term $y_{i_k^\tau \tau}^{k+1} - y_{i_k^\tau \tau}^k$ will go to zero, while another term $\frac{1}{m} \sum_{i=1}^m A_{i\tau}y_{i\tau}$ may not. Then in the presence of the compression error, the linear convergence rate could not be achieved by compressing this search direction directly. Hence, like ECLK, we introduce a vector $h_\tau^k$ to learn $u_\tau^k = \frac{1}{m} \sum_{i=1}^m A_{i\tau}y_{i\tau}$ iteratively. This learning scheme was first proposed in DIANA \citep{Samuel19} with the unbiased compressor. More precisely, we perform the following update on each node 
$$
h_{\tau}^{k+1} = h_\tau^k + Q_1(u_\tau^k - h_\tau^k), 
$$
where $Q_1$ is a contraction compressor. Now we apply the compression and error feedback mechanism to 
\begin{equation}\label{eq:searchd}
	u_\tau^k - h_\tau^k + A_{i_k^\tau \tau}(y_{i_k^\tau \tau}^{k+1} - y_{i_k^\tau \tau}^k), 
\end{equation} 
and add $h^k \eqdef \frac{1}{n} \sum_{\tau=1}^n h_\tau^k$ back after aggregation. We use $e_\tau^k$ to denote the compression error on each node, and add it to (\ref{eq:searchd}) before compression. After compression, $e_\tau^k$ is updated by the compression error at the current step: 
$$
e_\tau^{k+1} = e_\tau^k + u_\tau^k - h_\tau^k + A_{i_k^\tau \tau}(y_{i_k^\tau \tau}^{k+1} - y_{i_k^\tau \tau}^k) - Q(e_\tau^k + u_\tau^k - h_\tau^k + A_{i_k^\tau \tau}(y_{i_k^\tau \tau}^{k+1} - y_{i_k^\tau \tau}^k)), 
$$
where $Q$ is also a contraction compressor. The rest steps are the same as SPDC \citep{SPDC}. Next we introduce some useful variables. 

\vskip 2mm

Let $e^k \eqdef \tfrac{1}{n} \sum_{\tau=1}^n e_\tau^k$ and $u^k \eqdef \tfrac{1}{n} \sum_{\tau=1}^n u_\tau^k$ for $k\geq 0$. Define ${\tilde x}^k = x^k - \eta e^k$ for $k\geq 0$. We denote the optimal solution of the above saddle point problem as $(x^*, Y^*)$, where 
$$
Y^* = ((y^*_{11})^\top, ..., (y^*_{m1})^\top, ..., (y^*_{n1})^\top, ..., (y^*_{mn})^\top)^\top. 
$$ 
Now we are ready to construct some Lyapunov functions. 
For $k\geq 0$, define 
\begin{align*}
	\Phi_2^k & \eqdef \left(  \tfrac{1}{2\eta} + \tfrac{\lambda}{4}  \right)  \|{\tilde x}^{k} -x^*\|^2 + \left(  \tfrac{1}{4\sigma} + \tfrac{\gamma}{2}  \right) \tfrac{1}{n} \sum_{\tau=1}^n \sum_{i=1}^m  \|y_{i\tau}^{k} - y_{i\tau}^*\|^2 + \tfrac{3(\eta+ \lambda \eta^2)}{\delta n} \sum_{\tau=1}^n \|e_\tau^{k}\|^2 \\ 
	& \quad  +  f(x^{k}, Y^*)   - f(x^*, Y^*) + m\left(  f(x^*, Y^*) - f(x^*, Y^{k})   \right)  + \tfrac{42(1-\delta) (\eta+\lambda\eta^2)}{\delta^2 \delta_1 n} \sum_{\tau=1}^n \|h_\tau^{k} - u_\tau^{k} \|^2, 
\end{align*}
where $Y^k = ((y^k_{11})^\top, ..., (y^k_{m1})^\top,  ..., (y^k_{n1})^\top, ..., (y^k_{mn})^\top)^\top$, and 
\begin{align*}
	\Psi_2^k & \eqdef  \left(  \tfrac{1}{2\eta} + \tfrac{\lambda}{4}  \right) \|{\tilde x}^{k} -x^*\|^2 + \left(  \tfrac{1}{4\sigma} + \tfrac{\gamma}{2}  \right) \tfrac{1}{n} \sum_{\tau=1}^n \sum_{i=1}^m  \|y_{i\tau}^{k} - y_{i\tau}^*\|^2 + f(x^{k}, Y^*) - f(x^*, Y^*)   \\ 
	& \quad  + m\left(  f(x^*, Y^*) - f(x^*, Y^{k})   \right)  + \tfrac{3(\eta + \lambda \eta^2)}{\delta}\|e^{k}\|^2  + \tfrac{21(1-\delta)(\eta+\lambda \eta^2)}{\delta n^2} \sum_{\tau=1}^n  \|e_\tau^{k}\|^2 \\ 
	& \quad  + \tfrac{84(1-\delta)(\eta+\lambda\eta^2)}{5\delta^2 \delta_1}  \|h^{k} - u^{k}\|^2 + \tfrac{1512(1-\delta)(\eta+\lambda\eta^2)}{5\delta^2 \delta_1n^2} \sum_{\tau=1}^n\|h_\tau^{k} - u_\tau^{k}\|^2.  
\end{align*}

\noindent First, we introduce Lemma \ref{lm:ESO}, which is useful in the analysis of samplings in distributed systems.

\begin{lemma}\label{lm:ESO}
	Let $S = \{  (i^\tau, \tau) | \ i^\tau \mbox{ is chosen from $[m]$ uniformly and independently for all } \tau \in [n]  \}$. For any given $w_{i\tau} \in \R^t$ for $i\in[m]$ and $\tau\in [n]$, we have 
	$$
	\mathbb{E} \left[  \left\| \sum_{\tau=1}^n A_{i^\tau \tau} w_{i^\tau \tau} \right\|^2  \right] \leq \left( \frac{nR^2 + R_m^2}{m}  \right) \sum_{\tau=1}^n \sum_{i=1}^m \| w_{i\tau}\|^2. 
	$$
\end{lemma}

\noindent The following two lemmas show the evolution of the error terms $\sum_{\tau=1}^n \|e_\tau^k\|^2$ and $\|e^k\|^2$. The proofs are similar to that of Lemmas 3.4 and B.4 in \citep{eclk}, hence we omit them. 

\begin{lemma}\label{lm:etauk+1-SPDC}
	We have 
	\begin{align*}
		\frac{1}{n} \sum_{\tau=1}^n \mathbb{E}_k [\|e^{k+1}_\tau\|^2 ] & \leq \left(  1 - \frac{\delta}{2}  \right) \frac{1}{n} \sum_{\tau=1}^n \|e^k_\tau\|^2 + \frac{4(1-\delta)}{\delta n} \sum_{\tau=1}^n\|u_\tau^k - h_\tau^k\|^2 \\ 
		& \quad + \frac{(1-\delta)}{mn} \left(  \frac{4{\bar R}^2}{\delta} + R_m^2  \right) \sum_{\tau=1}^n \sum_{i=1}^m \|{\tilde y}^k_{i\tau} - y^k_{i\tau}\|^2. 
	\end{align*}
\end{lemma}

\begin{lemma}\label{lm:ek+1-SPDC} 
	Under Assumption \ref{as:expcompressor}, we have 
	\begin{align*}
		\mathbb{E}_k\|e^{k+1}\|^2 & \leq \left(  1 - \frac{\delta}{2}  \right) \|e^k\|^2 + \frac{2(1-\delta)\delta}{n^2} \sum_{\tau=1}^n \|e_\tau^k\|^2 + \frac{4(1-\delta)\delta}{n^2} \sum_{\tau=1}^n \|u_\tau^k - h_\tau^k\|^2 \\ 
		& \quad + \frac{4(1-\delta)}{\delta} \|u^k - h^k\|^2 + \frac{(1-\delta)}{mn} \left(  \frac{4R^2}{\delta} + \frac{5R_m^2}{n}  \right) \sum_{\tau=1}^n \sum_{i=1}^m \|{\tilde y}^k_{i\tau} - y^k_{i\tau}\|^2. 
	\end{align*}
	
\end{lemma}

\noindent We analyze the evolution of $\sum_{\tau=1}^n \|h_\tau^k - u_\tau^k\|^2$ and $\|h^k-u^k\|^2$ in the next two lemmas.

\begin{lemma}\label{lm:htauk+1-SPDC}
	We have 
	\begin{align*}
		\frac{1}{n} \sum_{\tau=1}^n \mathbb{E}_k [\|h_{\tau}^{k+1} - u_\tau^{k+1}\|^2] &\leq \left(  1 - \frac{\delta_1}{2}  \right) \frac{1}{n} \sum_{\tau=1}^n \|h_\tau^k - u_\tau^k\|^2 \\ 
		& \quad + \frac{1}{m^3 n} \left(  \frac{2(1-\delta_1) {\bar R}^2}{\delta_1} + R_m^2  \right) \sum_{\tau=1}^n \sum_{i=1}^m \|{\tilde y}^k_{i\tau} - y^k_{i\tau}\|^2. 
	\end{align*}
\end{lemma}

\begin{lemma}\label{lm:hk+1-SPDC}
	Under Assumption \ref{as:expcompressor}, we have 
	\begin{align*}
		\mathbb{E}_k \|h^{k+1} - u^{k+1} \|^2 &\leq  \left(  1 - \delta_1 \right) \|h^k - u^k\|^2 + \frac{\delta_1}{n^2} \sum_{\tau=1}^n\|h_\tau^k - u_\tau^k\|^2  \\ 
		& + \frac{1}{m^3 n} \left(  \frac{R^2}{\delta_1} + \frac{R_m^2}{n}  \right) \sum_{\tau=1}^n \sum_{i=1}^m \|{\tilde y}^k_{i\tau} - y^k_{i\tau}\|^2. 
	\end{align*}
\end{lemma}

\noindent The dual problem of problem (\ref{primal-spdc}) is 
\begin{align}
	\max_{Y \in \mathbb{R}^{tN}} & D(Y) \eqdef \min_{x \in \R^d} f(x, Y) = -\tfrac{1}{N} \sum_{\tau=1}^n \sum_{i=1}^m \phi^*_{i\tau} (y_{i\tau}) - g^*\left(  - \tfrac{1}{N} AY  \right), \label{dual-spdc}
\end{align}
where $g^*$ is the conjugate function of $g$ and 
\begin{equation}\label{eq:defA}
	A = [ A_{11}, ..., A_{m1}, ..., A_{n1}, ..., A_{mn} ] \in \R^{d \times tN}.
\end{equation}

\noindent Recall that $R^2 = \tfrac{1}{N} \lambda_{\rm max} (\sum_{\tau=1}^n \sum_{i=1}^m A_{i\tau} A_{i\tau}^\top) = \tfrac{1}{N} \|A\|^2$. We have the following lemma. 

\begin{lemma}\label{lm:PD-ecSPDC}[Lemma 3 in \citealp{SPDC}]
	Let Assumption \ref{as:ecSPDC} hold. Then for any point $(x, Y) \in {\rm dom}(f(x, Y))$, we have 
	$$
	P(x) \leq f(x, Y^*) + \tfrac{R^2}{ 2\gamma} \|x - x^*\|^2, \quad {\rm and} \quad	D(Y) \geq f(x^*, Y) - \tfrac{R^2}{2\lambda N} \|Y-Y^*\|^2. 
	$$
\end{lemma}

\begin{theorem}\label{th:pf-SPDC}
	Let Assumption \ref{as:contracQQ1-ecSPDC} and Assumption \ref{as:ecSPDC} hold. Set $\sigma = \tfrac{1}{2\cR_1} \sqrt{\tfrac{m\lambda}{\gamma}}$, $\eta = \tfrac{1}{2\cR_1}\sqrt{\tfrac{\gamma}{m\lambda}}$, and $\theta = 1 - \min\left\{  \tfrac{1}{m + 4\cR_1 \sqrt{m/(\lambda \gamma)}}, \tfrac{\delta}{6}, \tfrac{\delta_1}{6}  \right\}$, where $\cR_1 >0$ will be chosen later. 
	(i) Let $\cR_1^2  = \cR_2^2 \eqdef 2R^2 + \tfrac{2R_m^2}{n} + \tfrac{3(1-\delta)}{4}  \left(  \tfrac{14{\bar R}^2}{\delta^2}  + \tfrac{7R_m^2}{2\delta}   + \tfrac{84(1-\delta_1){\bar R}^2}{\delta^2 \delta_1^2 m^2} + \tfrac{42R_m^2}{\delta^2 \delta_1m^2} \right)$. Assume $\tfrac{{\cR}_2^2}{\lambda \gamma} \geq 1$. Then $$\mathbb{E}[\Phi_2^k] \leq  \epsilon \left(  \Phi_2^0 +  \tfrac{1}{4\sigma n} \sum_{\tau=1}^n \sum_{i=1}^m \|y_{i\tau}^{0} - y_{i\tau}^*\|^2 \right),$$ as long as 
	$
	k \geq {\cal O} \left(  \left(  \tfrac{1}{\delta} + \tfrac{1}{\delta_1} + m + \cR_2 \sqrt{\tfrac{m}{\lambda \gamma}}   \right) \log \tfrac{1}{\epsilon}  \right). 
	$
	In particular, if $\tfrac{1}{m} \leq {\cal O}(\delta_1)$, then the iteration complexity becomes 
	$$
	k \geq {\cal O} \left(  \left(  \left( R + \tfrac{R_m}{\sqrt{n}}  + \tfrac{\sqrt{(1-\delta)}{\bar R}}{\delta}  +  \tfrac{\sqrt{(1-\delta)}R_m}{\sqrt{\delta}}   \right) \sqrt{\tfrac{m}{\lambda \gamma}}  \tfrac{1}{\delta} + m \right) \log \tfrac{1}{\epsilon}  \right). 
	$$
	(ii)Let $\cR_1^2 = \cR_3^2$, where 
	\begin{align*}\cR_3^2 & \eqdef 2R^2 + \tfrac{2R_m^2}{n} + \tfrac{21(1-\delta)}{4} \left(  \tfrac{2R^2}{\delta^2} + \tfrac{11R_m^2}{2\delta n} + \tfrac{12(1-\delta){\bar R}^2}{\delta^2 n}  \tfrac{12R^2}{5\delta^2 \delta_1^2 m^2}  + \tfrac{228R_m^2}{5\delta^2 \delta_1m^2 n} + \tfrac{432(1-\delta_1){\bar R}^2}{5\delta^2 \delta_1^2m^2 n} \right). 
	\end{align*}
	Let Assumption \ref{as:expcompressor} hold and assume $\tfrac{{\cR}_3^2}{\lambda \gamma} \geq 1$. Then $
	\mathbb{E}[\Psi_2^k] \leq  \epsilon \left(  \Psi_2^0 +  \tfrac{1}{4\sigma n} \sum_{\tau=1}^n \sum_{i=1}^m \|y_{i\tau}^{0} - y_{i\tau}^*\|^2  \right)
	$ as long as $k \geq {\cal O} \left(  \left(  \tfrac{1}{\delta} + \tfrac{1}{\delta_1} + m + \cR_3 \sqrt{\tfrac{m}{\lambda \gamma}}   \right) \log \tfrac{1}{\epsilon}  \right)$.  
	If $\tfrac{1}{m} \leq {\cal O}(\delta_1)$, then the iteration complexity becomes 
	\begin{align*}
		k & \geq {\cal O} \left(  \left(  \left( R + \tfrac{R_m}{\sqrt{n}}  + \tfrac{\sqrt{(1-\delta)}{R}}{\delta}  +  \tfrac{\sqrt{(1-\delta)}R_m}{\sqrt{\delta n}}   \right) \sqrt{\tfrac{m}{\lambda \gamma}}  +  \tfrac{1}{\delta} + m \right) \log \tfrac{1}{\epsilon}  \right). 
	\end{align*}
\end{theorem}

From Lemma \ref{lm:PD-ecSPDC}, same as Corollary 4 in \citep{SPDC}, we can bound the primal-dual gap in the following theorem. 
\begin{theorem}
	Let Assumption \ref{as:ecSPDC} hold. Then we have 
	\begin{align*}
		P(x^k) - D(Y^k)   &\leq \left(  1 + \tfrac{R^2}{\lambda \gamma}  \right) \left(   f(x^{k}, Y^*) - f(x^*, Y^*)  + m\left(  f(x^*, Y^*) - f(x^*, Y^{k})   \right)    \right). 
	\end{align*}
\end{theorem}
Since $ f(x^{k}, Y^*) - f(x^*, Y^*) + m (  f(x^*, Y^*) - f(x^*, Y^{k}))$ is bounded by $\Phi_2^k$ or $\Psi_2^k$, the iteration complexity of the primal-dual gap can be deduced easily from Theorem~\ref{th:pf-SPDC}. Hence, we omit it.

\newpage

\section{PROOFS OF LEMMA \ref{lm:ESO}, LEMMA \ref{lm:htauk+1-SPDC}, LEMMA \ref{lm:hk+1-SPDC}, AND THEOREM \ref{th:pf-SPDC}} 

\subsection{Proof of Lemma \ref{lm:ESO}}

Let $W = (w_{11}^\top, ..., w_{m1}^\top, w_{12}^\top, ..., w_{m2}^\top, ..., w_{n1}^\top, ..., w_{mn}^\top)^\top \in \R^{tN}$. We have 
\begin{align*}
	& \quad \mathbb{E} \left[  \left\| \sum_{\tau=1}^n A_{i^\tau \tau} w_{i^\tau \tau} \right\|^2  \right] \\ 
	& = \mathbb{E} \left[ \sum_{\tau_1 \neq \tau_2} \langle A_{i^{\tau_1} \tau_1} w_{i^{\tau_1} \tau_1}, A_{i^{\tau_2} \tau_2} w_{i^{\tau_2} \tau_2} \rangle  \right] + \mathbb{E} \left[  \sum_{\tau=1}^n \| A_{i^\tau \tau} w_{i^\tau \tau}\|^2 \right] \\ 
	& = \sum_{\tau_1 \neq \tau_2} \left\langle \frac{1}{m} \sum_{i=1}^m A_{i \tau_1} w_{i \tau_1}, \frac{1}{m} \sum_{j=1}^m A_{j \tau_2} w_{j \tau_2} \right\rangle + \frac{1}{m} \sum_{\tau=1}^n \sum_{i=1}^m \|A_{i\tau} w_{i\tau}\|^2 \\ 
	& = \frac{1}{m^2} \sum_{\tau_1 \neq \tau_2} \sum_{i, j=1}^m \langle  A_{i \tau_1} w_{i \tau_1}, A_{j \tau_2} w_{j \tau_2}  \rangle + \frac{1}{m} \sum_{\tau=1}^n \sum_{i=1}^m \|A_{i\tau} w_{i\tau}\|^2 \\ 
	& =  \frac{1}{m^2} \sum_{\tau_1, \tau_2=1}^n  \sum_{i, j=1}^m \langle  A_{i \tau_1} w_{i \tau_1}, A_{j \tau_2} w_{j \tau_2}  \rangle - \frac{1}{m^2} \sum_{\tau=1}^n \sum_{i, j=1}^m \langle  A_{i \tau} w_{i \tau}, A_{j \tau} w_{j \tau}  \rangle  + \frac{1}{m} \sum_{\tau=1}^n \sum_{i=1}^m \|A_{i\tau} w_{i\tau}\|^2 \\ 
	& = \frac{1}{m^2} \|AW\|^2 - \frac{1}{m^2} \sum_{\tau=1}^n \left\| \sum_{i=1}^m A_{i\tau}w_{i\tau} \right\|^2 + \frac{1}{m} \sum_{\tau=1}^n \sum_{i=1}^m \|A_{i\tau} w_{i\tau}\|^2  \\ 
	& \leq  \frac{1}{m^2} \|AW\|^2 + \frac{1}{m} \sum_{\tau=1}^n \sum_{i=1}^m \|A_{i\tau} w_{i\tau}\|^2 \\ 
	& \leq \frac{NR^2}{m^2} \sum_{\tau=1}^n \sum_{i=1}^m \| w_{i\tau}\|^2 + \frac{1}{m} R_m^2  \sum_{\tau=1}^n \sum_{i=1}^m \| w_{i\tau}\|^2 \\ 
	& =  \left( \frac{nR^2}{m} +  \frac{1}{m} R_m^2  \right) \sum_{\tau=1}^n \sum_{i=1}^m \| w_{i\tau}\|^2, 
\end{align*}
where in the second equality, we use the fact that $i^{\tau_1}$ is indpendent of $i^{\tau_2}$ for $\tau_1 \neq \tau_2$, and in the last inequality, we use $\|A_{i\tau}\| \leq \max_{i, \tau} \|A_{i\tau}\| = R_m$ and $$\frac{1}{N}\lambda_{\rm max}(A A^\top)  = \frac{1}{N} \lambda_{\rm max} (\sum_{\tau=1}^n \sum_{i=1}^m A_{i\tau} A_{i\tau}^\top) = R^2. $$

\subsection{Proof of Lemma \ref{lm:htauk+1-SPDC}} 

First, we have 
\begin{align}
	& \quad  \mathbb{E}_k \|h_\tau^{k+1} - u_\tau^{k+1}\|^2 \nonumber \\ 
	& = \mathbb{E}_k \left\|h_\tau^{k+1} - u_\tau^k - \mathbb{E}_k[u_\tau^{k+1} - u_\tau^k] + \mathbb{E}_k[u_\tau^{k+1} - u_\tau^k] - (u_\tau^{k+1} - u_\tau^k) \right\|^2 \nonumber  \\ 
	& = \mathbb{E}_k\left\|h_\tau^{k+1} - u_\tau^k - \mathbb{E}_k[u_\tau^{k+1} - u_\tau^k] \right\|^2 + \mathbb{E}_k\left\| \mathbb{E}_k[u_\tau^{k+1} - u_\tau^k] - (u_\tau^{k+1} - u_\tau^k) \right\|^2 \nonumber \\ 
	& = \mathbb{E}_k\left\|h_\tau^{k+1} - u_\tau^k - \mathbb{E}_k[u_\tau^{k+1} - u_\tau^k] \right\|^2 + \mathbb{E}_k\|u_\tau^{k+1} - u_\tau^k\|^2 - \left\|  \mathbb{E}_k[u_\tau^{k+1} - u_\tau^k] \right\|^2 \nonumber \\ 
	& \leq (1+\beta) \mathbb{E}_k\|h_\tau^{k+1} - u_\tau^k\|^2 + \left(  1 + \frac{1}{\beta} - 1 \right) \left\|  \mathbb{E}_k[u_\tau^{k+1} - u_\tau^k] \right\|^2 + \mathbb{E}_k\|u_\tau^{k+1} - u_\tau^k\|^2  \label{eq:htauk+1} \\ 
	& \leq (1-\delta_1) (1+\beta) \|h_\tau^k - u_\tau^k\|^2 + \frac{1}{\beta} \left\|  \mathbb{E}_k[u_\tau^{k+1} - u_\tau^k] \right\|^2 + \mathbb{E}_k\|u_\tau^{k+1} - u_\tau^k\|^2  \nonumber \\ 
	& = \left(  1 - \frac{\delta_1}{2}  \right) \|h_\tau^k - u_\tau^k\|^2 + \frac{2(1-\delta_1)}{\delta_1}  \left\|  \mathbb{E}_k[u_\tau^{k+1} - u_\tau^k] \right\|^2 + \mathbb{E}_k\|u_\tau^{k+1} - u_\tau^k\|^2, \nonumber
\end{align}
where we use Young's inequality for any $\beta>0$ in the first inequality, in the second inequality we use the contraction property of $Q_1$, in the last equality we choose $\beta = \frac{\delta_1}{2(1-\delta_1)}$ when $\delta_1<1$. When $\delta_1=1$, it is easy to see that the above inequality also holds.

\noindent Since $u_\tau^{k+1} - u_\tau^k = \frac{1}{m} A_{i_k^\tau \tau} (y_{i_k^\tau \tau}^{k+1} - y_{i_k^\tau \tau}^k)$, we have 
\begin{align*}
	\mathbb{E}_k\|u_\tau^{k+1} - u_\tau^k\|^2 & \leq \frac{R_m^2}{m^2} \mathbb{E}_k \left\| y_{i_k^\tau \tau}^{k+1} - y_{i_k^\tau \tau}^k \right\|^2 \\ 
	& = \frac{R_m^2}{m^3} \sum_{i=1}^m \|{\tilde y}_{i\tau}^k - y_{i\tau}^k\|^2. 
\end{align*}

\noindent From $ \mathbb{E}_k[u_\tau^{k+1} - u_\tau^k] = \frac{1}{m^2} \sum_{i=1}^m A_{i\tau}({\tilde y}^k_{i\tau} - y^k_{i\tau})$, we can get 
\begin{align*}
	\left\|  \mathbb{E}_k[u_\tau^{k+1} - u_\tau^k] \right\|^2 & = \frac{1}{m^4} \left\| \sum_{i=1}^m A_{i\tau}({\tilde y}^k_{i\tau} - y^k_{i\tau}) \right\|^2 \\ 
	& \leq \frac{1}{m^4} \left\|[A_{1\tau}, ..., A_{m\tau}] \right\|^2 \cdot \sum_{i=1}^m \|{\tilde y}^k_{i\tau} - y^k_{i\tau}\|^2 \\ 
	& = \frac{1}{m^4} \lambda_{\rm max} \left(  \sum_{i=1}^m A_{i\tau} A_{i\tau}^\top  \right) \sum_{i=1}^m \|{\tilde y}^k_{i\tau} - y^k_{i\tau}\|^2  \\ 
	& \leq \frac{{\bar R}^2}{m^3} \sum_{i=1}^m \|{\tilde y}^k_{i\tau} - y^k_{i\tau}\|^2. 
\end{align*}

\noindent Combining the above three inequalities, we arrive at 
\begin{align*}
	\mathbb{E}_k \|h_\tau^{k+1} - u_\tau^{k+1}\|^2 &\leq  \left(  1 - \frac{\delta_1}{2}  \right) \|h_\tau^k - u_\tau^k\|^2 + \frac{1}{m^3} \left(  \frac{2(1-\delta_1) {\bar R}^2}{\delta_1} + R_m^2  \right) \sum_{i=1}^m \|{\tilde y}^k_{i\tau} - y^k_{i\tau}\|^2. 
\end{align*}

\noindent Summing up the above inequality from $\tau=1$ to $n$ and dividing both sides of the resulting inequality by $n$ , we can get the result.

\subsection{Proof of Lemma \ref{lm:hk+1-SPDC}}

First, same as (\ref{eq:htauk+1}), we can obtain 
$$
\mathbb{E}_k \|h^{k+1} - u^{k+1} \|^2 \leq (1+\beta) \|h^{k+1} - u^k\|^2 + \frac{1}{\beta} \left\| \mathbb{E}_k [u^{k+1} - u^k] \right\|^2 + \mathbb{E}_k \|u^{k+1} - u^k\|^2, 
$$
for any $\beta>0$. 

\noindent Under Assumption \ref{as:expcompressor}, same as the analysis of $\mathbb{E}_k\|h^{k+1} - \nabla f(w^k)\|^2$ in Lemma B.5 of \citep{eclk}, we can get 
$$
\mathbb{E}_k\|h^{k+1} - u^k\|^2 \leq (1-\delta_1)^2 \|h^k - u^k\|^2 + \frac{(1-\delta_1)\delta_1}{n^2} \sum_{\tau=1}^n\|h_\tau^k - u_\tau^k\|^2 . 
$$

\noindent  Combining the above two inequalities yields that 
\begin{align*}
	\mathbb{E}_k \|h^{k+1} - u^{k+1} \|^2 & \leq (1+\beta)(1-\delta_1)^2 \|h^k - u^k\|^2 + \frac{(1+\beta)(1-\delta_1)\delta_1}{n^2} \sum_{\tau=1}^n\|h_\tau^k - u_\tau^k\|^2 \\ 
	& \quad + \frac{1}{\beta} \left\| \mathbb{E}_k [u^{k+1} - u^k] \right\|^2 + \mathbb{E}_k \|u^{k+1} - u^k\|^2\\ 
	& = \left(  1 - \delta_1 \right) \|h^k - u^k\|^2 + \frac{\delta_1}{n^2} \sum_{\tau=1}^n\|h_\tau^k - u_\tau^k\|^2 \\ 
	& \quad + \frac{(1-\delta_1)}{\delta_1} \left\| \mathbb{E}_k [u^{k+1} - u^k] \right\|^2 + \mathbb{E}_k \|u^{k+1} - u^k\|^2, 
\end{align*}
where we choose $\beta = \frac{\delta_1}{1-\delta_1}$ when $\delta_1<1$. When $\delta_1=1$, $h^{k+1} = u^k$, thus the above inequality also holds. 

\noindent  Since $u_\tau^{k+1} - u_\tau^k = \frac{1}{m} A_{i_k^\tau \tau} (y_{i_k^\tau \tau}^{k+1} - y_{i_k^\tau \tau}^k)$, we have 
\begin{align*}
	\mathbb{E}_k\|u^{k+1} - u^k\|^2 & = \mathbb{E}_k \left\|  \frac{1}{mn} \sum_{\tau=1}^n A_{i_k^\tau \tau} (y_{i_k^\tau \tau}^{k+1} - y_{i_k^\tau \tau}^k)  \right\|^2 \\ 
	& = \mathbb{E}_k \left\|  \frac{1}{mn} \sum_{\tau=1}^n A_{i_k^\tau \tau} ({\tilde y}_{i_k^\tau \tau}^{k} - y_{i_k^\tau \tau}^k)  \right\|^2 \\ 
	& \overset{Lemma~\ref{lm:ESO}}{\leq} \frac{nR^2 + R_m^2}{m^3n^2} \sum_{\tau=1}^n \sum_{i=1}^m \|{\tilde y}^k_{i\tau} - y^k_{i\tau}\|^2. 
\end{align*}

\noindent From $\mathbb{E}_k[u^{k+1} - u^k] = \frac{1}{m^2 n} \sum_{\tau=1}^n \sum_{i=1}^m A_{i\tau}({\tilde y}^k_{i\tau} - y^k_{i\tau})$, we can obtain 
\begin{align*}
	\left\| \mathbb{E}_k [u^{k+1} - u^k] \right\|^2 & \leq \frac{1}{m^4 n^2} \|A\|^2 \sum_{\tau=1}^n \sum_{i=1}^m \|{\tilde y}^k_{i\tau} - y^k_{i\tau}\|^2  \\ 
	& = \frac{R^2}{m^3n} \sum_{\tau=1}^n \sum_{i=1}^m \|{\tilde y}^k_{i\tau} - y^k_{i\tau}\|^2, 
\end{align*}
where $A$ is defined in (\ref{eq:defA}). 
Combining the above three inequalities, we can get the result.

\newpage 

\subsection{Proof of Theorem \ref{th:pf-SPDC}}

(i) Let ${\tilde y}^k_{i\tau} \eqdef	\arg\max_{y\in \R^t} \left\{  \langle y, A_{i\tau}^\top z^k \rangle - \phi^*_{i\tau} (y) - \frac{1}{2\sigma}\|y - y_{i\tau}^k\|^2  \right\}$. Then similar to (47) in \citep{SPDC}, we can get 
\begin{align}
	& \quad \left(  \frac{1}{2\sigma} + \frac{(m-1)\gamma}{2m}  \right) \frac{1}{n} \sum_{\tau=1}^n \sum_{i=1}^m \|y_{i\tau}^k - y_{i\tau}^*\|^2  \nonumber \\
	& \geq \left(  \frac{1}{2\sigma} + \frac{\gamma}{2}  \right) \frac{1}{n} \sum_{\tau=1}^n \sum_{i=1}^m \mathbb{E}_k \|y_{i\tau}^{k+1} - y_{i\tau}^*\|^2 + \frac{1}{2\sigma n} \sum_{\tau=1}^n \sum_{i=1}^m \mathbb{E}_k \|y_{i\tau}^{k+1} - y_{i\tau}^k\|^2  \nonumber  \\ 
	& \quad + \frac{1}{n} \mathbb{E}_k\left[  \sum_{\tau=1}^n (\phi^*_{i_k^\tau \tau}(y^{k+1}_{i_k^\tau \tau}) - \phi^*_{i_k^\tau \tau}(y^{k}_{i_k^\tau \tau}) ) \right] \nonumber \\ 
	& \quad + \frac{1}{N} \sum_{\tau=1}^n \sum_{i=1}^m (\phi^*_{i\tau}(y^k_{i\tau}) - \phi^*_{i\tau}(y^*_{i\tau})) - \mathbb{E}_k \langle u^k-u^* + m(u^{k+1} - u^k), z^k \rangle, \label{eq:yk+1}
\end{align}
where $u^k = \frac{1}{N} \sum_{\tau=1}^n \sum_{i=1}^m A_{i\tau} y^k_{i\tau}$ and $u^* = \frac{1}{N} \sum_{\tau=1}^n \sum_{i=1}^m A_{i\tau} y^*_{i\tau}$. 

\noindent  From the update rule of $x^{k+1}$ and optimality condition, we have 
$$
\partial g(x^{k+1}) + h^k+\Delta^k + \frac{1}{\eta} (x^{k+1} - x^k) = 0. 
$$
For $h^k+\Delta^k$, we have 
\begin{align*}
	h^k + \Delta^k & = \frac{1}{n} \sum_{\tau=1}^n (h_\tau^k + \Delta_\tau^k) \\ 
	& = \frac{1}{n} \sum_{\tau=1}^n (h_\tau^k + e_\tau^k + A_{i_k^\tau \tau} (y_{i_k^\tau \tau}^{k+1} - y_{i_k^\tau \tau}^k) + u_\tau^k - h_\tau^k - e_\tau^{k+1} ) \\ 
	& = \frac{1}{n} \sum_{\tau=1}^n (e_\tau^k - e_\tau^{k+1} + m(u_\tau^{k+1} - u_\tau^k)) \\ 
	& = e^k - e^{k+1} + u^k + m(u^{k+1} - u^k). 
\end{align*}

\noindent By using the above two equalities, we can obtain 
\begin{align*}
	{\tilde x}^{k+1} & = x^{k+1} - \eta e^{k+1} \\
	& = x^k - \eta(\partial g(x^{k+1}) + h^k + \Delta^k) - \eta e^{k+1} \\ 
	& = x^k - \eta (\partial g(x^{k+1}) + e^k - e^{k+1} + u^k + m(u^{k+1} - u^k) - e^{k+1}) \\ 
	& = {\tilde x}^k - \eta (\partial g(x^{k+1}) + u^k + m(u^{k+1} - u^k) ). 
\end{align*}

\noindent Then similar to Lemma B.3 in \citep{eclk}, we can get 
\begin{align*}
	\langle u^k + m(u^{k+1} - u^k), x^* - x^{k+1} \rangle &\geq \left(  \frac{1}{2\eta} + \frac{\lambda}{4}  \right) \|{\tilde x}^{k+1} -x^*\|^2 - \frac{1}{2\eta} \|{\tilde x}^k - x^*\|^2 - \frac{\eta}{2}\|e^k\|^2 \\ 
	& \quad - \left(  \frac{\eta}{2} + \frac{\lambda \eta^2}{2}  \right)\|e^{k+1}\|^2 + \frac{1}{4\eta} \|x^{k+1} - x^k\|^2 + g(x^{k+1}) - g(x^*). 
\end{align*}

\noindent Rearranging terms and taking conditional expectation, we arrive at 
\begin{align}
	\frac{1}{2\eta} \|{\tilde x}^k - x^*\|^2 &\geq  \left(  \frac{1}{2\eta} + \frac{\lambda}{4}  \right) \mathbb{E}_k \|{\tilde x}^{k+1} -x^*\|^2 + \frac{1}{4\eta} \mathbb{E}_k \|x^{k+1} - x^k\|^2 + \mathbb{E}_k[ g(x^{k+1}) - g(x^*)] \nonumber \\ 
	& \quad + \mathbb{E}_k \langle u^k + m(u^{k+1} - u^k), x^{k+1} - x^* \rangle - \frac{\eta}{2}\|e^k\|^2 - \left(  \frac{\eta}{2} + \frac{\lambda \eta^2}{2}  \right) \mathbb{E}_k \|e^{k+1}\|^2. \label{eq:xk+1}
\end{align}

\noindent Similar to (51) in \citep{SPDC}, we can get 
\begin{align*}
	&\quad f(x^{k+1}, Y^*) - f(x^*, Y^*) + m\left(  f(x^*, Y^*) - f(x^*, Y^{k+1})   \right) - (m-1) \left(  f(x^*, Y^*) - f(x^*, Y^k)  \right) \\ 
	& = \frac{1}{N} \sum_{\tau=1}^n \sum_{i=1}^m \left(  \phi^*_{i\tau} (y^k_{i\tau}) - \phi^*_{i\tau} (y^*_{i\tau})  \right) + \frac{1}{n} \sum_{\tau=1}^n \left(\phi^*_{i_k^\tau \tau}(y^{k+1}_{i_k^\tau \tau}) - \phi^*_{i_k^\tau \tau}(y^{k}_{i_k^\tau \tau}) \right) + g(x^{k+1}) - g(x^*) \\ 
	& \quad + \langle u^*, x^{k+1} \rangle - \langle u^k, x^* \rangle + m \langle u^k - u^{k+1}, x^* \rangle. 
\end{align*}

\noindent Combining (\ref{eq:yk+1}), (\ref{eq:xk+1}), and the above equality after taking conditional expectation, we can obtain 
\begin{align}
	& \quad \frac{1}{2\eta} \|{\tilde x}^k - x^*\|^2 + \left(  \frac{1}{2\sigma} + \frac{(m-1)\gamma}{2m}  \right) \frac{1}{n} \sum_{\tau=1}^n \sum_{i=1}^m \|y_{i\tau}^k - y_{i\tau}^*\|^2 + (m-1) \left(  f(x^*, Y^*) - f(x^*, Y^k)  \right) \nonumber \\ 
	& \geq  \left(  \frac{1}{2\eta} + \frac{\lambda}{4}  \right) \mathbb{E}_k \|{\tilde x}^{k+1} -x^*\|^2 + \left(  \frac{1}{2\sigma} + \frac{\gamma}{2}  \right) \frac{1}{n} \sum_{\tau=1}^n \sum_{i=1}^m \mathbb{E}_k \|y_{i\tau}^{k+1} - y_{i\tau}^*\|^2 +  \frac{1}{4\eta} \mathbb{E}_k \|x^{k+1} - x^k\|^2 \nonumber \\ 
	& \quad + \frac{1}{2\sigma m n} \sum_{\tau=1}^n \sum_{i=1}^m \|{\tilde y}_{i\tau}^{k} - y_{i\tau}^k\|^2 + \mathbb{E}_k\left[ f(x^{k+1}, Y^*) - f(x^*, Y^*) + m\left(  f(x^*, Y^*) - f(x^*, Y^{k+1})   \right)   \right] \nonumber \\ 
	& \quad + \mathbb{E}_k \langle u^k-u^* + m(u^{k+1} - u^k), x^{k+1} - z^k \rangle  - \frac{\eta}{2}\|e^k\|^2 - \left(  \frac{\eta}{2} + \frac{\lambda \eta^2}{2}  \right) \mathbb{E}_k \|e^{k+1}\|^2, \label{eq:u1}
\end{align}
where we also use the fact that $\sum_{\tau=1}^n \sum_{i=1}^m \mathbb{E}_k \|y_{i\tau}^{k+1} - y_{i\tau}^k\|^2 = \frac{1}{m} \sum_{\tau=1}^n \sum_{i=1}^m \|{\tilde y}_{i\tau}^{k} - y_{i\tau}^k\|^2$.

\noindent Next we estimate the term $\mathbb{E}_k \langle u^k-u^* + m(u^{k+1} - u^k), x^{k+1} - z^k \rangle $. We have 
\begin{align}
	& \quad \langle u^k-u^* + m(u^{k+1} - u^k), x^{k+1} - z^k \rangle \nonumber \\ 
	& = \frac{1}{N} \langle AY^k - AY^* + mA(Y^{k+1} - Y^k), x^{k+1} - x^k - \theta(x^k - x^{k-1}) \rangle \nonumber \\ 
	& = \frac{1}{N} \langle AY^{k+1} - AY^*, x^{k+1} - x^k \rangle - \frac{\theta}{N} \langle AY^k - AY^*, x^k-x^{k-1} \rangle \nonumber \\ 
	& \quad + \frac{m-1}{N} \langle AY^{k+1} - AY^k, x^{k+1}-x^k\rangle + \frac{m\theta}{N} \langle AY^{k+1} - AY^k, x^k-x^{k-1} \rangle, \label{eq:u2}
\end{align}
where we define $x^{-1} \eqdef x^0$ to guarantee $z^0=x^0$. By Cauchy-Schwarz inequality, we have 
\begin{align*}
	|\langle AY^{k+1} - AY^k, x^{k+1}-x^k\rangle | & \leq \frac{n}{8\eta} \|x^{k+1}-x^k\|^2 + \frac{2\eta}{n} \|A(Y^{k+1} - Y^k)\|^2. 
\end{align*}

\noindent By Lemma \ref{lm:ESO}, we further have 
\begin{align*}
	\mathbb{E}_k  \|A(Y^{k+1} - Y^k)\|^2 & = \mathbb{E}_k \left\|  \sum_{\tau=1}^n A_{i_k^\tau \tau}(y^{k+1}_{i_k^\tau \tau} - y^k_{i_k^\tau \tau}) \right\|^2 \\ 
	& = \mathbb{E}_k \left\|  \sum_{\tau=1}^n A_{i_k^\tau \tau}({\tilde y}^k_{i_k^\tau \tau} - y^k_{i_k^\tau \tau}) \right\|^2 \\ 
	& \overset{Lemma~\ref{lm:ESO}}{\leq} \frac{nR^2 + R_m^2}{m} \sum_{\tau=1}^n \sum_{i=1}^m \|{\tilde y}^k_{i\tau} - y^k_{i\tau}\|^2. 
\end{align*}

\noindent Thus, we arrive at 
$$
\mathbb{E}_k \langle AY^{k+1} - AY^k, x^{k+1}-x^k\rangle \geq - \frac{n}{8\eta} \mathbb{E}_k \|x^{k+1}-x^k\|^2 - \frac{2\eta (nR^2 + R_m^2)}{N} \sum_{\tau=1}^n \sum_{i=1}^m \|{\tilde y}^k_{i\tau} - y^k_{i\tau}\|^2. 
$$

\noindent Simiarly, we can get 
$$
\mathbb{E}_k \langle AY^{k+1} - AY^k, x^{k}-x^{k-1} \rangle \geq - \frac{n}{8\eta}  \|x^k-x^{k-1}\|^2 - \frac{2\eta (nR^2 + R_m^2)}{N} \sum_{\tau=1}^n \sum_{i=1}^m \|{\tilde y}^k_{i\tau} - y^k_{i\tau}\|^2. 
$$

\noindent Combining the above two inequalities with (\ref{eq:u1}) and (\ref{eq:u2}), we have 
\begin{align}
	&\quad \frac{1}{2\eta} \|{\tilde x}^k - x^*\|^2 + \left(  \frac{1}{2\sigma} + \frac{(m-1)\gamma}{2m}  \right) \frac{1}{n} \sum_{\tau=1}^n \sum_{i=1}^m \|y_{i\tau}^k - y_{i\tau}^*\|^2 + (m-1) \left(  f(x^*, Y^*) - f(x^*, Y^k)  \right) \nonumber \\
	& \quad + \theta \left(  f(x^k, Y^*) - f(x^*, Y^*)  \right) + \frac{\theta}{8\eta} \|x^k - x^{k-1}\|^2 + \frac{\theta}{N} \langle AY^k - AY^*, x^k-x^{k-1} \rangle \nonumber \\ 
	& \geq \left(  \frac{1}{2\eta} + \frac{\lambda}{4}  \right) \mathbb{E}_k \|{\tilde x}^{k+1} -x^*\|^2 + \left(  \frac{1}{2\sigma} + \frac{\gamma}{2}  \right) \frac{1}{n} \sum_{\tau=1}^n \sum_{i=1}^m \mathbb{E}_k \|y_{i\tau}^{k+1} - y_{i\tau}^*\|^2 +  \frac{1}{8\eta} \mathbb{E}_k \|x^{k+1} - x^k\|^2 \nonumber \\ 
	& \quad + \mathbb{E}_k\left[ f(x^{k+1}, Y^*) - f(x^*, Y^*) + m\left(  f(x^*, Y^*) - f(x^*, Y^{k+1})   \right)   \right] + \frac{1}{N} \mathbb{E}_k \langle AY^{k+1} - AY^*, x^{k+1} - x^k \rangle  \nonumber \\ 
	& \quad + \frac{1}{2N} \left(  \frac{1}{\sigma} - \frac{8\eta (nR^2 + R_m^2)}{n}  \right)  \sum_{\tau=1}^n \sum_{i=1}^m \|{\tilde y}^k_{i\tau} - y^k_{i\tau}\|^2  - \frac{\eta}{2}\|e^k\|^2 - \left(  \frac{\eta}{2} + \frac{\lambda \eta^2}{2}  \right) \mathbb{E}_k \|e^{k+1}\|^2, \label{eq:11-SPDC}
\end{align}
where we add the nonnegative term $\theta \left(  f(x^k, Y^*) - f(x^*, Y^*)  \right)$ to the left-hand side of the above inequality.

\noindent From Lemma \ref{lm:etauk+1-SPDC}, we have 
\begin{align}
	& \quad \frac{3(\eta+ \lambda \eta^2)}{\delta n} \sum_{\tau=1}^n \mathbb{E}_k\|e_\tau^{k+1}\|^2 + \frac{\eta}{2n}\sum_{\tau=1}^n \|e_\tau^{k}\|^2 + \frac{\eta+\lambda \eta^2}{2n}\sum_{\tau=1}^n \mathbb{E}_k\|e_\tau^{k+1}\|^2 \nonumber \\ 
	& \leq \frac{\eta + \lambda \eta^2}{n} \left( \left(  \frac{3}{\delta} + \frac{1}{2}  \right) \left(  1 - \frac{\delta}{2}  \right) + \frac{1}{2} \right) \sum_{\tau=1}^n \|e_\tau^{k}\|^2 + \frac{4(1-\delta) (\eta+\lambda\eta^2)}{\delta n}  \left(  \frac{3}{\delta} + \frac{1}{2}  \right) \sum_{\tau=1}^n \|h_\tau^k - u_\tau^k\|^2 \nonumber \\ 
	& \quad + \frac{(1-\delta)(\eta+\lambda\eta^2)}{mn}  \left(  \frac{3}{\delta} + \frac{1}{2}  \right) \left(  \frac{4{\bar R}^2}{\delta} + R_m^2  \right) \sum_{\tau=1}^n \sum_{i=1}^m \|{\tilde y}^k_{i\tau} - y^k_{i\tau}\|^2 \nonumber \\ 
	& \leq \frac{3(\eta+ \lambda \eta^2)}{\delta n} \left(  1 - \frac{\delta}{6}  \right) \sum_{\tau=1}^n \|e_\tau^{k}\|^2 +  \frac{14(1-\delta) (\eta+\lambda\eta^2)}{\delta^2 n} \sum_{\tau=1}^n \|h_\tau^k - u_\tau^k\|^2 \nonumber \\ 
	& \quad + \frac{(1-\delta)(\eta+\lambda\eta^2)}{mn} \left(   \frac{14{\bar R}^2}{\delta^2}  + \frac{7R_m^2}{2\delta} \right) \sum_{\tau=1}^n \sum_{i=1}^m \|{\tilde y}^k_{i\tau} - y^k_{i\tau}\|^2. \label{eq:22-SPDC}
\end{align}

\noindent From Lemma \ref{lm:htauk+1-SPDC}, we have 
\begin{align}
	& \quad \frac{42(1-\delta) (\eta+\lambda\eta^2)}{\delta^2 \delta_1 n} \sum_{\tau=1}^n \mathbb{E}_k \|h_\tau^{k+1} - u_\tau^{k+1} \|^2  +  \frac{14(1-\delta) (\eta+\lambda\eta^2)}{\delta^2 n} \sum_{\tau=1}^n \|h_\tau^k - u_\tau^k\|^2 \nonumber \\ 
	& \leq \frac{42(1-\delta) (\eta+\lambda\eta^2)}{\delta^2 \delta_1 n} \left(  1 - \frac{\delta_1}{2} + \frac{\delta_1}{3}  \right)  \sum_{\tau=1}^n \|h_\tau^k - u_\tau^k\|^2 \nonumber \\ 
	& \quad +  \frac{42(1-\delta) (\eta+\lambda\eta^2)}{\delta^2 \delta_1 m^3 n}  \left(  \frac{2(1-\delta_1) {\bar R}^2}{\delta_1} + R_m^2  \right) \sum_{\tau=1}^n \sum_{i=1}^m \|{\tilde y}^k_{i\tau} - y^k_{i\tau}\|^2. \label{eq:33-SPDC}
\end{align}

\noindent From (\ref{eq:11-SPDC}), (\ref{eq:22-SPDC}), (\ref{eq:33-SPDC}), and the fact that $\|e^k\|^2 \leq \frac{1}{n} \sum_{\tau=1}^n \|e_\tau^k\|^2$, we arrive at 
\begin{align}
	& \quad \left(  \frac{1}{2\eta} + \frac{\lambda}{4}  \right) \mathbb{E}_k \|{\tilde x}^{k+1} -x^*\|^2 + \left(  \frac{1}{2\sigma} + \frac{\gamma}{2}  \right) \frac{1}{n} \sum_{\tau=1}^n \sum_{i=1}^m \mathbb{E}_k \|y_{i\tau}^{k+1} - y_{i\tau}^*\|^2 +  \frac{1}{8\eta} \mathbb{E}_k \|x^{k+1} - x^k\|^2 \nonumber \\ 
	& \qquad + \mathbb{E}_k\left[ f(x^{k+1}, Y^*) - f(x^*, Y^*) + m\left(  f(x^*, Y^*) - f(x^*, Y^{k+1})   \right)   \right] + \frac{1}{N} \mathbb{E}_k \langle AY^{k+1} - AY^*, x^{k+1} - x^k \rangle \nonumber \\ 
	& \qquad + \frac{3(\eta+ \lambda \eta^2)}{\delta n} \sum_{\tau=1}^n \mathbb{E}_k\|e_\tau^{k+1}\|^2 + \frac{42(1-\delta) (\eta+\lambda\eta^2)}{\delta^2 \delta_1 n} \sum_{\tau=1}^n \mathbb{E}_k \|h_\tau^{k+1} - u_\tau^{k+1} \|^2 \nonumber \\ 
	& \leq \frac{1}{2\eta} \|{\tilde x}^k - x^*\|^2 + \left(  \frac{1}{2\sigma} + \frac{(m-1)\gamma}{2m}  \right) \frac{1}{n} \sum_{\tau=1}^n \sum_{i=1}^m \|y_{i\tau}^k - y_{i\tau}^*\|^2 + (m-1) \left(  f(x^*, Y^*) - f(x^*, Y^k)  \right) \nonumber  \\
	& \quad + \theta \left(  f(x^k, Y^*) - f(x^*, Y^*)  \right) + \frac{\theta}{8\eta} \|x^k - x^{k-1}\|^2 + \frac{\theta}{N} \langle AY^k - AY^*, x^k-x^{k-1} \rangle  \nonumber \\ 
	& \quad + \frac{3(\eta+ \lambda \eta^2)}{\delta n} \left(  1 - \frac{\delta}{6}  \right) \sum_{\tau=1}^n \|e_\tau^{k}\|^2 + \frac{42(1-\delta) (\eta+\lambda\eta^2)}{\delta^2 \delta_1 n} \left(  1 - \frac{\delta_1}{6}  \right)  \sum_{\tau=1}^n \|h_\tau^k - u_\tau^k\|^2  \nonumber \\ 
	& \quad - \frac{1}{mn} \left(  \frac{1}{2 \sigma} - \frac{4\eta (nR^2 + R_m^2)}{n} - (1-\delta) (\eta+\lambda\eta^2) \left(  \frac{14{\bar R}^2}{\delta^2}  + \frac{7R_m^2}{2\delta}   + \frac{84(1-\delta_1){\bar R}^2}{\delta^2 \delta_1^2 m^2} + \frac{42R_m^2}{\delta^2 \delta_1m^2} \right)   \right)  \nonumber \\ 
	& \quad \cdot \sum_{\tau=1}^n \sum_{i=1}^m \|{\tilde y}^k_{i\tau} - y^k_{i\tau}\|^2. \label{eq:ev-SPDC}
\end{align}

\noindent Define 
\begin{align*}
	\Phi_1^k & \eqdef \left(  \frac{1}{2\eta} + \frac{\lambda}{4}  \right) \|{\tilde x}^{k} -x^*\|^2 + \left(  \frac{1}{2\sigma} + \frac{\gamma}{2}  \right) \frac{1}{n} \sum_{\tau=1}^n \sum_{i=1}^m \|y_{i\tau}^{k} - y_{i\tau}^*\|^2 +  \frac{1}{8\eta} \|x^{k} - x^{k-1}\|^2 \\ 
	& \quad +  f(x^{k}, Y^*) - f(x^*, Y^*) + m\left(  f(x^*, Y^*) - f(x^*, Y^{k})   \right)  + \frac{1}{N}  \langle AY^{k} - AY^*, x^{k} - x^{k-1} \rangle \\ 
	& \quad + \frac{3(\eta+ \lambda \eta^2)}{\delta n} \sum_{\tau=1}^n  \|e_\tau^{k}\|^2 + \frac{42(1-\delta) (\eta+\lambda\eta^2)}{\delta^2 \delta_1 n} \sum_{\tau=1}^n  \|h_\tau^{k} - u_\tau^{k} \|^2, 
\end{align*}
for $k\geq 0$, where $x^{-1} = x^0$. Assume $\frac{{\cR}_2^2}{\lambda \gamma} \geq 1$. Then $\lambda \eta = \frac{1}{2\cR_2} \sqrt{\frac{\lambda\gamma}{m}} \leq \frac{1}{2}$, and thus 
\begin{align*}
	& \quad \frac{4\eta (nR^2 + R_m^2)}{n} + (1-\delta) (\eta+\lambda\eta^2) \left(  \frac{14{\bar R}^2}{\delta^2}  + \frac{7R_m^2}{2\delta}   + \frac{84(1-\delta_1){\bar R}^2}{\delta^2 \delta_1^2 m^2} + \frac{42R_m^2}{\delta^2 \delta_1m^2} \right) \\ 
	& \leq \frac{4\eta (nR^2 + R_m^2)}{n} + \frac{3(1-\delta)\eta}{2} \left(  \frac{14{\bar R}^2}{\delta^2}  + \frac{7R_m^2}{2\delta}   + \frac{84(1-\delta_1){\bar R}^2}{\delta^2 \delta_1^2 m^2} + \frac{42R_m^2}{\delta^2 \delta_1m^2} \right) \\ 
	& = 2\eta \cR_2^2  = \frac{1}{2\sigma}. 
\end{align*}

\noindent From (\ref{eq:ev-SPDC}), the above inequality, and the definition of $\Phi_1^k$, we can get 
$$
\mathbb{E}_k[\Phi_1^{k+1}] \leq \theta \Phi_1^k, 
$$
where we use $\left( 1 - \frac{\delta}{6} \right) \leq \theta$, $\left( 1 - \frac{\delta_1}{6} \right) \leq \theta$, 
$$
\frac{m-1}{m} = 1 - \frac{1}{m} \leq \theta,  \quad \left. \frac{1}{2\eta}\right/ \left(  \frac{1}{2\eta} + \frac{\lambda}{4}  \right) = 1 - \frac{1}{1 + 4\cR_2\sqrt{m/(\lambda\gamma)}} \leq \theta,  
$$
and 
$$
\left.  \left(  \frac{1}{2\sigma} + \frac{(m-1)\gamma}{2m}  \right) \right/ \left(  \frac{1}{2\sigma} + \frac{\gamma}{2}  \right) = 1 - \frac{1}{m + m/(\gamma\sigma)} = 1 - \frac{1}{m + 2\cR_2\sqrt{m/(\lambda\gamma)}} \leq \theta. 
$$

\noindent By the tower property, we further have $\mathbb{E}[\Phi_1^{k+1}] \leq \theta \mathbb{E} [\Phi_1^k]$. Apply this relation recursively, we can obtain 
\begin{equation}\label{eq:phi1k-SPDC}
	\mathbb{E} [\Phi_1^k] \leq \theta^k \Phi_1^0. 
\end{equation}

\noindent From the definition of $\Phi_2^k$, we know that 
$$
\Phi_1^k = \Phi_2^k + \frac{1}{4\sigma n} \sum_{\tau=1}^n \sum_{i=1}^m \|y_{i\tau}^{k} - y_{i\tau}^*\|^2 + \frac{1}{8\eta} \|x^{k} - x^{k-1}\|^2 + \frac{1}{N}  \langle AY^{k} - AY^*, x^{k} - x^{k-1} \rangle. 
$$

\noindent From Young's inequality, we have 
\begin{align}
	\frac{1}{N} \left| \langle AY^{k} - AY^*, x^{k} - x^{k-1} \rangle \right| & \leq \frac{\|x^k - x^{k-1}\|^2}{8\eta} + \frac{\|A\|^2 \|Y^k-Y^*\|^2}{N^2/(2\eta)} \nonumber \\ 
	& = \frac{\|x^k - x^{k-1}\|^2}{8\eta}  + \frac{2\eta R^2}{N} \sum_{\tau=1}^n \sum_{i=1}^m \|y_{i\tau}^{k} - y_{i\tau}^*\|^2  \nonumber \\ 
	& \leq \frac{\|x^k - x^{k-1}\|^2}{8\eta}  + \frac{\eta \cR_2^2}{n} \sum_{\tau=1}^n \sum_{i=1}^m \|y_{i\tau}^{k} - y_{i\tau}^*\|^2  \nonumber \\ 
	& = \frac{\|x^k - x^{k-1}\|^2}{8\eta} + \frac{1}{4\sigma n} \sum_{\tau=1}^n \sum_{i=1}^m \|y_{i\tau}^{k} - y_{i\tau}^*\|^2, \label{eq:last-SPDC}
\end{align}
which indicates that $\Phi_2^k \leq \Phi_1^k$ for $k\geq 0$. Therefore, from (\ref{eq:phi1k-SPDC}) we have 
$$
\mathbb{E}[\Phi_2^k] \leq \theta^k \Phi_1^0 = \theta^k \left(  \Phi_2^0 +  \frac{1}{4\sigma n} \sum_{\tau=1}^n \sum_{i=1}^m \|y_{i\tau}^{0} - y_{i\tau}^*\|^2 \right). 
$$

\noindent Finally, from $\frac{R_m^2}{m} \leq {\bar R}^2$, we can get the results.

\noindent (ii) Under Assumption \ref{as:expcompressor}, from Lemma \ref{lm:etauk+1-SPDC} and Lemma \ref{lm:ek+1-SPDC}, we have 
\begin{align*}
	& \quad \frac{3(\eta + \lambda \eta^2)}{\delta} \mathbb{E}_k \|e^{k+1}\|^2 + \frac{21(1-\delta)(\eta+\lambda \eta^2)}{\delta n^2} \sum_{\tau=1}^n \mathbb{E}_k \|e_\tau^{k+1}\|^2 + \frac{\eta}{2}\|e^k\|^2 + \frac{\eta+\lambda \eta^2}{2} \mathbb{E}_k \|e^{k+1}\|^2 \\ 
	& \leq \frac{3(\eta + \lambda \eta^2)}{\delta} \left(  1 - \frac{\delta}{6}  \right) \|e^k\|^2 +  \frac{21(1-\delta)(\eta+\lambda \eta^2)}{\delta n^2} \left(  1 - \frac{\delta}{6}  \right) \sum_{\tau=1}^n \|e_\tau^k\|^2 \\ 
	& \quad + \frac{14(1-\delta)(\eta + \lambda \eta^2)}{\delta^2} \|h^k - u^k\|^2 + \frac{84(1-\delta)(\eta+\lambda\eta^2)}{\delta^2 n^2} \sum_{\tau=1}^n \|h_\tau^k - u_\tau^k\|^2 \\ 
	& \quad + \frac{7(1-\delta)(\eta+\lambda\eta^2)}{\delta mn} \left(  \frac{2R^2}{\delta} + \frac{11R_m^2}{2n} + \frac{12(1-\delta){\bar R}^2}{\delta n}  \right) \sum_{\tau=1}^n \sum_{i=1}^m \|{\tilde y}^k_{i\tau} - y^k_{i\tau}\|^2. 
\end{align*}

\noindent From Lemma \ref{lm:htauk+1-SPDC} and Lemma \ref{lm:hk+1-SPDC}, we have 
\begin{align*}
	& \quad \frac{84(1-\delta)(\eta+\lambda\eta^2)}{5\delta^2 \delta_1} \mathbb{E}_k \|h^{k+1} - u^{k+1}\|^2 + \frac{1512(1-\delta)(\eta+\lambda\eta^2)}{5\delta^2 \delta_1n^2} \sum_{\tau=1}^n \mathbb{E}_k\|h_\tau^{k+1} - u_\tau^{k+1}\|^2 \\
	& \qquad + \frac{14(1-\delta)(\eta + \lambda \eta^2)}{\delta^2} \|h^k - u^k\|^2 + \frac{84(1-\delta)(\eta+\lambda\eta^2)}{\delta^2 n^2} \sum_{\tau=1}^n \|h_\tau^k - u_\tau^k\|^2 \\ 
	& \leq \frac{84(1-\delta)(\eta+\lambda\eta^2)}{5\delta^2 \delta_1} \left(  1 - \frac{\delta_1}{6}  \right) \|h^k - u^k\|^2 + \frac{1512(1-\delta)(\eta+\lambda\eta^2)}{5\delta^2 \delta_1n^2} \left(  1- \frac{\delta_1}{6}  \right) \sum_{\tau=1}^n \|h_\tau^k - u_\tau^k\|^2 \\ 
	& \quad + \frac{84(1-\delta)(\eta+\lambda\eta^2)}{5\delta^2 \delta_1m^3 n} \left(  \frac{R^2}{\delta_1} + \frac{19R_m^2}{n} + \frac{36(1-\delta_1){\bar R}^2}{\delta_1 n}  \right) \sum_{\tau=1}^n \sum_{i=1}^m \|{\tilde y}^k_{i\tau} - y^k_{i\tau}\|^2. 
\end{align*}

\noindent Combining the above two inequalities and (\ref{eq:11-SPDC}), we arrive at 
\begin{align}
	& \quad \left(  \frac{1}{2\eta} + \frac{\lambda}{4}  \right) \mathbb{E}_k \|{\tilde x}^{k+1} -x^*\|^2 + \left(  \frac{1}{2\sigma} + \frac{\gamma}{2}  \right) \frac{1}{n} \sum_{\tau=1}^n \sum_{i=1}^m \mathbb{E}_k \|y_{i\tau}^{k+1} - y_{i\tau}^*\|^2 +  \frac{1}{8\eta} \mathbb{E}_k \|x^{k+1} - x^k\|^2 \nonumber \\ 
	& \qquad + \mathbb{E}_k\left[ f(x^{k+1}, Y^*) - f(x^*, Y^*) + m\left(  f(x^*, Y^*) - f(x^*, Y^{k+1})   \right)   \right] \nonumber \\ 
	& \qquad + \frac{3(\eta + \lambda \eta^2)}{\delta} \mathbb{E}_k\|e^{k+1}\|^2 + \frac{21(1-\delta)(\eta+\lambda \eta^2)}{\delta n^2} \sum_{\tau=1}^n \mathbb{E}_k \|e_\tau^{k+1}\|^2 + \frac{1}{N} \mathbb{E}_k \langle AY^{k+1} - AY^*, x^{k+1} - x^k \rangle  \nonumber \\
	& \qquad + \frac{84(1-\delta)(\eta+\lambda\eta^2)}{5\delta^2 \delta_1} \mathbb{E}_k \|h^{k+1} - u^{k+1}\|^2 + \frac{1512(1-\delta)(\eta+\lambda\eta^2)}{5\delta^2 \delta_1n^2} \sum_{\tau=1}^n \mathbb{E}_k\|h_\tau^{k+1} - u_\tau^{k+1}\|^2 \nonumber \\ 
	& \leq \frac{1}{2\eta} \|{\tilde x}^k - x^*\|^2 + \left(  \frac{1}{2\sigma} + \frac{(m-1)\gamma}{2m}  \right) \frac{1}{n} \sum_{\tau=1}^n \sum_{i=1}^m \|y_{i\tau}^k - y_{i\tau}^*\|^2 + (m-1) \left(  f(x^*, Y^*) - f(x^*, Y^k)  \right) \nonumber  \\
	& \quad + \theta \left(  f(x^k, Y^*) - f(x^*, Y^*)  \right) + \frac{\theta}{8\eta} \|x^k - x^{k-1}\|^2 + \frac{\theta}{N} \langle AY^k - AY^*, x^k-x^{k-1} \rangle  \nonumber \\ 
	& \quad + \frac{3(\eta + \lambda \eta^2)}{\delta} \left(  1 - \frac{\delta}{6}  \right) \|e^k\|^2 +  \frac{21(1-\delta)(\eta+\lambda \eta^2)}{\delta n^2} \left(  1 - \frac{\delta}{6}  \right) \sum_{\tau=1}^n \|e_\tau^k\|^2  \nonumber \\ 
	& \quad + \frac{84(1-\delta)(\eta+\lambda\eta^2)}{5\delta^2 \delta_1} \left(  1 - \frac{\delta_1}{6}  \right) \|h^k - u^k\|^2 + \frac{1512(1-\delta)(\eta+\lambda\eta^2)}{5\delta^2 \delta_1n^2} \left(  1- \frac{\delta_1}{6}  \right) \sum_{\tau=1}^n \|h_\tau^k - u_\tau^k\|^2  \nonumber \\
	& \quad - \frac{1}{mn}  \sum_{\tau=1}^n \sum_{i=1}^m \|{\tilde y}^k_{i\tau} - y^k_{i\tau}\|^2 \cdot \left(  \frac{1}{2 \sigma} - \frac{4\eta (nR^2 + R_m^2)}{n} - 7(1-\delta) (\eta+\lambda\eta^2) \right.  \nonumber \\
	& \quad \cdot \left. \left(  \frac{2R^2}{\delta^2} + \frac{11R_m^2}{2\delta n} + \frac{12(1-\delta){\bar R}^2}{\delta^2 n} + \frac{12R^2}{5\delta^2 \delta_1^2 m^2}  + \frac{228R_m^2}{5\delta^2 \delta_1m^2 n} + \frac{432(1-\delta_1){\bar R}^2}{5\delta^2 \delta_1^2m^2 n} \right)   \right) . \label{eq:ev-SPDC-2}
\end{align}

\noindent  Define 
\begin{align*}
	\Psi_1^k & \eqdef  \left(  \frac{1}{2\eta} + \frac{\lambda}{4}  \right) \|{\tilde x}^{k} -x^*\|^2 + \left(  \frac{1}{2\sigma} + \frac{\gamma}{2}  \right) \frac{1}{n} \sum_{\tau=1}^n \sum_{i=1}^m  \|y_{i\tau}^{k} - y_{i\tau}^*\|^2 +  \frac{1}{8\eta} \|x^{k} - x^{k-1}\|^2  \\ 
	& \qquad + f(x^{k}, Y^*) - f(x^*, Y^*) + m\left(  f(x^*, Y^*) - f(x^*, Y^{k})   \right)   + \frac{1}{N} \mathbb{E}_k \langle AY^{k} - AY^*, x^{k} - x^{k-1} \rangle  \\ 
	& \qquad + \frac{3(\eta + \lambda \eta^2)}{\delta}\|e^{k}\|^2 + \frac{21(1-\delta)(\eta+\lambda \eta^2)}{\delta n^2} \sum_{\tau=1}^n  \|e_\tau^{k}\|^2  \\
	& \qquad + \frac{84(1-\delta)(\eta+\lambda\eta^2)}{5\delta^2 \delta_1}  \|h^{k} - u^{k}\|^2 + \frac{1512(1-\delta)(\eta+\lambda\eta^2)}{5\delta^2 \delta_1n^2} \sum_{\tau=1}^n\|h_\tau^{k} - u_\tau^{k}\|^2, 
\end{align*}
for $k\geq 0$, where $x^{-1} = x^0$.

\noindent Assume $\frac{{\cR}_3^2}{\lambda \gamma} \geq 1$. Then $\lambda \eta = \frac{1}{2\cR_3} \sqrt{\frac{\lambda\gamma}{m}} \leq \frac{1}{2}$, and thus 
\begin{align*}
	& \quad \frac{4\eta (nR^2 + R_m^2)}{n} + 7(1-\delta) (\eta+\lambda\eta^2) \\ 
	& \quad \cdot  \left(  \frac{2R^2}{\delta^2} + \frac{11R_m^2}{2\delta n} + \frac{12(1-\delta){\bar R}^2}{\delta^2 n} + \frac{12R^2}{5\delta^2 \delta_1^2 m^2}  + \frac{228R_m^2}{5\delta^2 \delta_1m^2 n} + \frac{432(1-\delta_1){\bar R}^2}{5\delta^2 \delta_1^2m^2 n} \right) \\ 
	& \leq 2\eta \cR_3^2 = \frac{1}{2\sigma}. 
\end{align*}

\noindent From (\ref{eq:ev-SPDC-2}), the above inequality, and the definition of $\Psi_1^k$, we can get 
$$
\mathbb{E}_k[\Psi_1^{k+1}] \leq \theta \Psi_1^k. 
$$
\noindent By the tower property, we further have $\mathbb{E}[\Psi_1^{k+1}] \leq \theta \mathbb{E} [\Psi_1^k]$. Apply this relation recursively, we can obtain 
\begin{equation}\label{eq:psi1k-SPDC}
	\mathbb{E} [\Psi_1^k] \leq \theta^k \Psi_1^0. 
\end{equation}

\noindent From the definition of $\Psi_2^k$, we know 
$$
\Psi_1^k = \Psi_2^k + \frac{1}{4\sigma n} \sum_{\tau=1}^n \sum_{i=1}^m \|y_{i\tau}^{k} - y_{i\tau}^*\|^2 + \frac{1}{8\eta} \|x^{k} - x^{k-1}\|^2 + \frac{1}{N}  \langle AY^{k} - AY^*, x^{k} - x^{k-1} \rangle. 
$$
\noindent From (\ref{eq:last-SPDC}), we have $\Psi_2^k \leq \Psi_1^k$ for $k\geq 0$. Thus, from (\ref{eq:psi1k-SPDC}) we can obtain 
$$
\mathbb{E}[\Psi_2^k] \leq \theta^k \Psi_1^0 = \theta^k \left(  \Psi_2^0 +  \frac{1}{4\sigma n} \sum_{\tau=1}^n \sum_{i=1}^m \|y_{i\tau}^{0} - y_{i\tau}^*\|^2  \right). 
$$

\noindent At the end, from $\frac{R_m^2}{m} \leq {\bar R}^2$ and $\frac{{\bar R}^2}{n} \leq R^2$, we can get the results.

\newpage 

\section{PROOFS FOR EC-LSVRG + CATALYST}

\subsection{Proof of Lemma \ref{lm:eclsvrg-1}}

First, from Theorem 2.10 in \citep{ecsdca} and the initialization rules of $h_{\tau, (k)}^0$ and $e_{\tau, (k)}^0$, we have 
$$
\mathbb{E}[G_k({\bar x}_{(k)}^K) - G_k^*] \leq \frac{9(\lambda+\kappa)\|x_{(k)}^0-x^*_{(k)}\|^2 + 2(G_k(x_{(k)}^0) - G_k^*)}{1 - (1 - {\tilde \theta})^{K+1}} (1 - {\tilde \theta})^K, 
$$
where we denote ${\tilde \theta} \eqdef \min\{  \frac{(\lambda+\kappa)\eta}{2}, \frac{\delta}{4}, \frac{\delta_1}{4}, \frac{p}{4}  \}$. Since $G_k$ is $(\lambda+\kappa)$-strongly convex, we have 
$$
G_k(x) - G_k^* \geq \frac{\lambda+\kappa}{2}\|x-x^*_{(k)}\|^2, 
$$
for any $x\in \R^d$, which indicates that 
$$
\mathbb{E}[G_k({\bar x}_{(k)}^K) - G_k^*] \leq \frac{18(G_k(x_{(k)}^0) - G_k^*)}{1 - (1 - {\tilde \theta})^{K+1}} (1 - {\tilde \theta})^K. 
$$
Noticing that $\ln(1-a) + a \leq 0$ for any $a \in (0, 1)$, we have $(1-{\tilde \theta})^{\frac{1}{\tilde \theta}} \leq \frac{1}{e}<0.37$. Now we first let $K\geq \frac{1}{\tilde \theta}$, then we have $(1-{\tilde \theta})^{K+1} \leq 0.37$, which yields 
$$
\mathbb{E}[G_k({\bar x}_{(k)}^K) - G_k^*] \leq {30(G_k(x_{(k)}^0) - G_k^*)} (1 - {\tilde \theta})^K. 
$$
Then similar to the proof of Lemma C.1 in \citep{catalyst}, but choosing $T_0$ as 
$$
T_0 = \max \left\{  \frac{1}{\tilde \theta}, \frac{1}{\tilde \theta} \log \left(  \frac{1}{1-e^{-{\tilde \theta}}} \frac{30(G_k(x_{(k)}^0) - G_k^*)}{\epsilon_k}  \right) \right\} 
$$
instead, we can obtain 
$$
\mathbb{E}[T_k] \leq \max \left\{  \frac{1}{\tilde \theta},  \frac{1}{\tilde \theta} \log \left(  \frac{60(G_k(x_{(k)}^0) - G_k^*)}{{\tilde \theta} \epsilon_k}  \right)  \right\} + 1. 
$$
Within the above inequality, similar to the proof of Proposition 3.2 in \citep{catalyst}, we can get $\mathbb{E}[T_k] \leq {\tilde {\cal O}} (1/{\tilde \theta})$, where the notation ${\tilde {\cal O}}$ hides some constants and some logorithmic dependencies in $\lambda$, $\kappa$, and ${\tilde \theta}$. At last, by the stepsize rule in Theorem 2.10 in \citep{ecsdca}, we can obtain the result.

\subsection{Proof of Lemma \ref{lm:eclsvrg-1-re}}

From the inequality above (27) in the proof of Theorem 2.10 in \citep{ecsdca}, we have 
$$
\mathbb{E}[G_k(x_{(k)}^K) - G_k^*] \leq \frac{18}{\eta} \mathbb{E}[\Phi_{3,(k)}^K] \leq \frac{18}{\eta} (1 - {\tilde \theta}_1)^K \Phi_{3,(k)}^0, 
$$
where we denote ${\tilde \theta}_1 \eqdef \min\left\{  \frac{(\lambda+\kappa)\eta}{2}, \frac{\delta}{4}, \frac{\delta_1}{4}, \frac{p}{4}  \right\}$. Thus, we can get 
$$
\mathbb{E}[\Phi_{3, (k)}^K + G_k(x_{(k)}^K) - G_k^*] \leq \left( 1 + \frac{18}{\eta}  \right)  (1 - {\tilde \theta}_1)^K \Phi_{3,(k)}^0. 
$$
Then from Lemma C.1 in \citep{catalyst}, we know 
$$
\mathbb{E}[T_k] \leq {\tilde {\cal O}} \left(  \frac{1}{{\tilde \theta}_1} \log \left(  \frac{(1+18/\eta) \Phi_{3,(k)}^0}{\epsilon_k}  \right)  \right) = {\tilde {\cal O}} \left(  \frac{1}{{\tilde \theta}_1} \log \left(  \frac{\Phi_{3,(k)}^0}{\epsilon_k}  \right)  \right). 
$$
Next we will show that $\log \left(  \frac{\Phi_{3,(k)}^0}{\epsilon_k}  \right) \leq {\tilde {\cal O}}(1)$, which concludes the proof. From the definition of $\Phi_{3, (k)}^K$ and the initialization rule at each outer iteration, we have 
\begin{align*}
	\Phi_{3, (k)}^0 & = \|x^{k-1} - e_{(k)}^0 - x_{(k)}^*\|^2 + \tfrac{12(L_f+\kappa)\eta}{n \delta} \sum_{\tau=1}^n \|e_{\tau, (k))}^0\|^2 + \tilde \eta (G_k(x^{k-1}) - G_k^*) \\ 
	& \quad + \tfrac{192(1-\delta)(L_f+\kappa)\eta^3}{\delta^2\delta_1 n} \sum_{\tau=1}^n \|h_{\tau, (k)}^{0} - \nabla f^{\tau} (x^{k-1}) - \nabla \psi(x^{k-1}) - \kappa (x^{k-1} - y^{k-1})\|^2, 
\end{align*}
where we denote ${\tilde \eta} = \tfrac{4}{3p}  \left(   \tfrac{48(1-\delta) (L_f+\kappa) \eta^3}{\delta} \left(  \tfrac{4({\bar L}+\kappa)}{\delta} + L + \kappa  + \tfrac{16({\bar L}+\kappa)p}{\delta \delta_1} \left(  1 + \tfrac{2p}{\delta_1}  \right) \right) + \tfrac{4(L+\kappa)\eta^2}{n}  \right)$. We estimate each term in the above equality respectively. Since $G_k$ is $(\lambda+\kappa)$-strongly convex, we have 
$$
\|x^{k-1} - e_{(k)}^0 - x_{(k)}^*\|^2 \leq 2\|x^{k-1}-x_{(k)}^*\|^2 + 2\|e_{(k)}^0\|^2 \leq \frac{4}{\lambda+\kappa} (G_k(x^{k-1})-G_k^*) + 2\|e_{(k)}^0\|^2. 
$$
For the third term, define $G_k^{(\tau)}(x) \eqdef f^{(\tau)}(x) + \psi(x) + \frac{\kappa}{2}\|x-y^{k-1}\|^2$ for simplicity. If $h_{\tau, (k)}^0 = h_{\tau, (k-1)}^{T_{k-1}}$, then we have 
\begin{align*}
	& \quad \frac{1-\delta}{n} \sum_{\tau=1}^n \|h_{\tau, (k)}^0 - \nabla f^{\tau} (x^{k-1}) - \nabla \psi(x^{k-1}) - \kappa (x^{k-1} - y^{k-1})\|^2 \\
	& = \frac{1-\delta}{n} \sum_{\tau=1}^n \|h_{\tau, (k-1)}^{T_{k-1}} - \nabla f^{\tau} (x^{k-1}) - \nabla \psi(x^{k-1}) - \kappa (x^{k-1} - y^{k-1})\|^2 \\ 
	& = \quad \frac{1-\delta}{n} \sum_{\tau=1}^n \|h_{\tau, (k-1)}^{T_{k-1}} - \nabla G_k^{\tau} (x^{k-1})\|^2 \\ 
	& \leq \frac{3(1-\delta)}{n} \sum_{\tau=1}^n \left(\|h_{\tau, (k-1)}^{T_{k-1}} - \nabla G_{k-1}^{(\tau)}(w_{(k-1)}^{T_{k-1}})\|^2 \right. \\ 
	& \quad \left.  + \| \nabla G_{k-1}^{(\tau)}(w_{(k-1)}^{T_{k-1}}) -  \nabla G_{k-1}^{(\tau)}(x^{k-1})  \|^2 + \|\nabla G_{k-1}^{(\tau)}(x^{k-1}) - \nabla G_{k}^{(\tau)}(x^{k-1})\|^2 \right).
\end{align*}

\noindent Since $\Phi_{3, (k-1)}^{T_{k-1}} + G_{k-1}({x}_{(k-1)}^{T_{k-1}}) - G_{k-1}^* \leq \epsilon_{k-1}$, we have 
$$
\frac{3(1-\delta)}{n} \sum_{\tau=1}^n \|h_{\tau, (k-1)}^{T_{k-1}} - \nabla G_{k-1}^{(\tau)}(w_{(k-1)}^{T_{k-1}})\|^2 \leq \frac{\delta^2\delta_1}{192(L_f+\kappa)\eta^3} \Phi_{3, (k-1)}^{T_{k-1}} \leq \frac{\delta^2\delta_1}{192(L_f+\kappa)\eta^3} \epsilon_{k-1}. 
$$
\noindent From the smoothness of $G_{k-1}^{(\tau)}$, we have 
\begin{align*}
	& \quad \frac{1}{n} \sum_{\tau=1}^n \| \nabla G_{k-1}^{(\tau)}(w_{(k-1)}^{T_{k-1}}) -  \nabla G_{k-1}^{(\tau)}(x^{k-1})  \|^2  \\ 
	& \leq \frac{2}{n} \sum_{\tau=1}^n \|\nabla G_{k-1}^{(\tau)}(w_{(k-1)}^{T_{k-1}}) - \nabla G_{k-1}^{(\tau)}(x_{(k-1)}^*) \|^2 + \frac{2}{n} \sum_{\tau=1}^n \|\nabla G_{k-1}^{(\tau)}(x^{k-1}) - \nabla G_{k-1}^{(\tau)}(x_{(k-1)}^*) \|^2 \\ 
	& \leq 4({\bar L} + \kappa) \left(  G_{k-1}(w_{(k-1)}^{T_{k-1}}) - G_{k-1}^* + G_{k-1}(x^{k-1}) - G_{k-1}^*  \right) \\ 
	& \leq \frac{4({\bar L}+\kappa)}{{\tilde \eta}} \Phi_{3, (k-1)}^{T_{k-1}} +  4({\bar L} + \kappa) (G_{k-1}(x^{k-1}) - G_{k-1}^*) \\ 
	& \leq  4({\bar L} + \kappa) \left(  \frac{1}{\tilde \eta} + 1  \right) \epsilon_{k-1}. 
\end{align*}

\noindent For $\frac{1}{n} \sum_{\tau=1}^n \|\nabla G_{k-1}^{(\tau)}(x^{k-1}) - \nabla G_{k}^{(\tau)}(x^{k-1})\|^2$, we have 
$$
\frac{1}{n} \sum_{\tau=1}^n \|\nabla G_{k-1}^{(\tau)}(x^{k-1}) - \nabla G_{k}^{(\tau)}(x^{k-1})\|^2 = \kappa^2 \|y^{k-1} - y^{k-2}\|^2. 
$$

\noindent For the third term, if $h_{\tau, (k)}^0 = h_{\tau, (k-1)}^{T_{k-1}} + \kappa(y^{k-2} - y^{k-1})$, then we have 
\begin{align*}
	& \quad \frac{1-\delta}{n} \sum_{\tau=1}^n \|h_{\tau, (k)}^0 - \nabla f^{\tau} (x^{k-1}) - \nabla \psi(x^{k-1}) - \kappa (x^{k-1} - y^{k-1})\|^2 \\
	& = \frac{1-\delta}{n} \sum_{\tau=1}^n \|h_{\tau, (k-1)}^{T_{k-1}} - \nabla f^{\tau} (x^{k-1}) - \nabla \psi(x^{k-1}) - \kappa (x^{k-1} - y^{k-2})\|^2 \\ 
	& = \quad \frac{1-\delta}{n} \sum_{\tau=1}^n \|h_{\tau, (k-1)}^{T_{k-1}} - \nabla G_{k-1}^{\tau} (x^{k-1})\|^2 \\ 
	& \leq \frac{3(1-\delta)}{n} \sum_{\tau=1}^n \left(\|h_{\tau, (k-1)}^{T_{k-1}} - \nabla G_{k-1}^{(\tau)}(w_{(k-1)}^{T_{k-1}})\|^2   + \| \nabla G_{k-1}^{(\tau)}(w_{(k-1)}^{T_{k-1}}) -  \nabla G_{k-1}^{(\tau)}(x^{k-1})  \|^2 \right).
\end{align*}

\noindent If $e_{\tau, (k)}^0=0$, then $\|e_{(k)}^0\|^2 = \|e_{\tau, (k)}^0\|^2=0$. If $e_{\tau, (k)}^0 = e_{\tau, (k-1)}^{T_{k-1}}$, then 
$$
\|e_{\tau, (k)}^0\|^2 \leq \frac{1}{n} \sum_{\tau=1}^n \|e_{\tau, (k)}^0\|^2 \leq \frac{\delta}{12(L_f + \kappa)\eta} \Phi_{3, (k-1)}^{T_{k-1}} \leq \frac{\delta}{12(L_f + \kappa)\eta} \epsilon_{k-1}. 
$$

\noindent Moreover, for any $a\geq 1$ and $b\geq 1$, we have $\log(a+b) \leq \log(2\max\{a, b\}) \leq \log (a) + \log (b) + 1$. Using the above estimations, we conclude that 
$$
\log\left(  \frac{\Phi_{3, (k)}^0}{\epsilon_k}  \right) \leq {\tilde {\cal O}} \left(  \frac{G_k(x^{k-1}) - G_k^*}{\epsilon_k}  \right) + {\tilde {\cal O}} \left(  \frac{\|y^{k-1}-y^{k-2}\|^2}{\epsilon_k}  \right) + {\cal O}(1) 
$$

\noindent Finally, from Lemmas B.1 and B.2 in \citep{catalyst}, and similar to the proof of Proposition 3.2 in \citep{catalyst}, we can get $\log \left(  \frac{\Phi_{3,(k)}^0}{\epsilon_k}  \right) \leq {\tilde {\cal O}}(1)$.

\newpage

\section{PROOFS FOR EC-SDCA + CATALYST}

\subsection{Proof of Lemma \ref{lm:ecsdca-1}}

In the proof, we borrow the methodology of \citep{accSDCA} for an accelerated proximal stochastic dual coordinate ascent method. First, from Theorem 3.5 in \citep{ecsdca}, we have 
$$
\mathbb{E}[\Psi_{3,(k)}^K] \leq \left( 1 - \min\left\{ \theta, \frac{\delta}{4} \right\} \right)^K \epsilon_{D,(k)}^0, 
$$
for $K\geq1$, and $\theta$ is chosen as in (9) in \citep{ecsdca}. Moreover, from the proof of Theorem 3.5 in \citep{ecsdca}, it is easy to verify that for $K\geq 1$ 
$$
\mathbb{E}[\epsilon_{P,(k)}^K] = \mathbb{E}[G_k(x_{(k)}^K) - G_k^*] \leq \frac{1}{\theta}  \left( 1 - \min\left\{ \theta, \frac{\delta}{4} \right\} \right)^K \epsilon_{D,(k)}^0. 
$$
\noindent Combining the above two inequalities, we get 
\begin{equation}\label{eq:psi3+2g}
	\mathbb{E} [\sqrt{4n+\delta mn} \Psi^K_{3,(k)} + 2(G_k(x_{(k)}^{K+1}) - G_k^*)] \leq \left(  \sqrt{4n+\delta mn} + \frac{2}{\theta} \right) \left( 1 - \min\left\{ \theta, \frac{\delta}{4} \right\} \right)^K \epsilon_{D,(k)}^0, 
\end{equation}
for $K\geq 1$. Then from (\ref{eq:psi3+2g}) and similar to Lemma C.1 in \citep{catalyst}, we know $\mathbb{E}[T_1] \leq {\tilde {\cal O}} \left( \frac{1}{\delta} + m + \frac{{\cal M}_1}{\lambda + \kappa} \right)$. Next we consider $k>1$ case. 

\vskip 2mm 

\noindent Define ${\tilde \lambda} = \lambda + \kappa$, ${\tilde \xi}(x) = \frac{\lambda}{{\tilde \lambda}} \xi(x) + \frac{\kappa}{2{\tilde \lambda}}\|x\|^2$, and $\xi_k(x) = \frac{\lambda}{{\tilde \lambda}}\xi(x) + \frac{\kappa}{2{\tilde \lambda}}\|x-y^{k-1}\|^2$. Denote the objective function of the dual problem of minimizing $G_k$ as $D_k(\alpha)$ and $v(\alpha) = \frac{1}{{\tilde \lambda}N} \sum_{\tau=1}^n \sum_{i=1}^m A_{i\tau}\alpha_{i\tau}$. Then $G_k(x) = \frac{1}{N}\sum_{\tau=1}^n \sum_{i=1}^m \phi_{i\tau}(A_{i\tau}^\top x) + {\tilde \lambda}\xi_k(x)$ and $D_k(\alpha) = -\frac{1}{N}\sum_{\tau=1}^n\sum_{i=1}^m \phi^*_{i\tau}(-\alpha_{i\tau}) - {\tilde \lambda}\xi_k^*(v(\alpha))$. Noticing that 
\begin{equation}\label{eq:xikstar}
	\xi_k^*(u) = \max_{x} \langle x, u+\frac{\kappa}{{\tilde \lambda}} y^{k-1} \rangle - {\tilde \xi}(x) - \frac{\kappa}{2{\tilde \lambda}}\|y^{k-1}\|^2 = {\tilde \xi}^*\left( u+\frac{\kappa}{{\tilde \lambda}} y^{k-1} \right) - \frac{\kappa}{2{\tilde \lambda}}\|y^{k-1}\|^2, 
\end{equation}
we can also write $D_k(\alpha) = -\frac{1}{N}\sum_{\tau=1}^n\sum_{i=1}^m \phi^*_{i\tau}(-\alpha_{i\tau}) - {\tilde \lambda} {\tilde \xi}^*\left( v(\alpha)+\frac{\kappa}{{\tilde \lambda}} y^{k-1} \right) + \frac{\kappa}{2}\|y^{k-1}\|^2$. Then we have 
\begin{align}
	& \quad -D_k(\alpha^{k-1}) + D_{k-1}(\alpha^{k-1}) \nonumber \\ 
	&= {\tilde \lambda}{\tilde \xi}^*\left( v(\alpha^{k-1}) + \frac{\kappa}{{\tilde \lambda}}y^{k-1} \right) - {\tilde \lambda}{\tilde \xi}^*\left( v(\alpha^{k-1}) + \frac{\kappa}{{\tilde \lambda}}y^{k-2} \right) + \frac{\kappa}{2}\|y^{k-2}\|^2 - \frac{\kappa}{2}\|y^{k-1}\|^2 \nonumber \\ 
	& \leq \kappa \langle \nabla {\tilde \xi}^*\left( v(\alpha^{k-1}) + \frac{\kappa}{{\tilde \lambda}}y^{k-2} \right), y^{k-1}-y^{k-2} \rangle + \frac{\kappa^2}{2{\tilde \lambda}} \|y^{k-1}-y^{k-2}\|^2 + \frac{\kappa}{2}\|y^{k-2}\|^2 - \frac{\kappa}{2}\|y^{k-1}\|^2 \nonumber \\
	& = \kappa \langle \nabla { \xi}_{k-1}^*\left( v(\alpha^{k-1}) \right), y^{k-1}-y^{k-2} \rangle + \frac{\kappa^2}{2{\tilde \lambda}} \|y^{k-1}-y^{k-2}\|^2 + \frac{\kappa}{2}\|y^{k-2}\|^2 - \frac{\kappa}{2}\|y^{k-1}\|^2, \label{eq:diffDk}
\end{align}
where we use ${\tilde \xi}^*$ is 1-smooth in the first inequality and (\ref{eq:xikstar}) in the last equality. From (33) in \citep{ecsdca}, we know ${\tilde u_{(k-1)}}^{T_{k-1}} = v(\alpha^{T_{k-1}}_{(k-1)})= v(\alpha^{k-1}) =  u_{(k-1)}^{T_{k-1}} + e_{(k-1)}^{T_{k-1}}$. Moreover, from the update of EC-SDCA, we know the output at the $k$-th outer iteration $x^{k-1} = x_{(k-1)}^{T_{k-1}+1} = \nabla {\xi}^*_{k-1}( u_{(k-1)}^{T_{k-1}})$. Thus, for $\langle \nabla { \xi}_{k-1}^*\left( v(\alpha^{k-1}) \right), y^{k-1}-y^{k-2} \rangle$, we have 
\begin{align*}
	& \quad  \langle \nabla {\xi}_{k-1}^*\left( v(\alpha^{k-1}) \right), y^{k-1}-y^{k-2} \rangle - \langle x^{k-1}, y^{k-1}-y^{k-2} \rangle \\
	& =  \langle \nabla {\xi}_{k-1}^*\left({\tilde u_{(k-1)}}^{T_{k-1}} \right) - \nabla {\tilde \xi}^*_{k-1}( u_{(k-1)}^{T_{k-1}}), y^{k-1}-y^{k-2} \rangle \\ 
	& \leq \frac{1}{2} \|\nabla {\xi}_{k-1}^*\left({\tilde u_{(k-1)}}^{T_{k-1}} \right) - \nabla {\xi}^*_{k-1}( u_{(k-1)}^{T_{k-1}})\|^2 + \frac{1}{2}\|y^{k-1} - y^{k-2}\|^2 \\
	& \leq \frac{1}{2} \|e_{(k-1)}^{T_{k-1}}\|^2 + \frac{1}{2}\|y^{k-1} - y^{k-2}\|^2,
\end{align*}
where we use ${\xi}_{k-1}^*$ is 1-smooth in the last inequality. From (\ref{eq:diffDk}), the above inequality, and ${\tilde \lambda}>\kappa$, we arrive at 
\begin{align*}
	& \quad -D_k(\alpha^{k-1}) + D_{k-1}(\alpha^{k-1}) \\ 
	& \leq  \kappa \langle x^{k-1}, y^{k-1}-y^{k-2} \rangle + \frac{\kappa}{2}\|e_{(k-1)}^{T_{k-1}}\|^2 + \kappa \|y^{k-1} - y^{k-2}\|^2 + \frac{\kappa}{2}\|y^{k-2}\|^2 - \frac{\kappa}{2}\|y^{k-1}\|^2. 
\end{align*}

\noindent Furthermore, we have
$$
G_k(x^{k-1}) = G_{k-1}(x^{k-1}) + \frac{\kappa}{2}\|y^{k-1}\|^2 - \frac{\kappa}{2}\|y^{k-2}\|^2 + \kappa \langle x^{k-1}, y^{k-2}-y^{k-1} \rangle.
$$

\noindent Combining the above two inequalities, we obtain 
$$
G_k(x^{k-1}) - D_k(\alpha^{k-1}) \leq G_{k-1}(x^{k-1}) - D_{k-1}(\alpha^{k-1}) +  \frac{\kappa}{2}\|e_{(k-1)}^{T_{k-1}}\|^2 + \kappa \|y^{k-1} - y^{k-2}\|^2. 
$$

\noindent Since ${\bar R}^2 \leq nR^2$ and $R_m^2\leq mnR^2$, we have 
$$
\Psi_{3,(k)}^K \geq G_k^* - D_k(\alpha_{(k)}^K)  + \frac{2\rho}{\delta}\|e_{(k)}^K\|^2 \geq G_k^* -  D_k(\alpha_{(k)}^K) + \frac{\lambda+\kappa}{\sqrt{2n+\delta mn}}\|e_{(k)}^K\|^2, 
$$
which implies that 
$$
G_{k-1}^* - D_{k-1}(\alpha^{k-1}) \leq \frac{\sqrt{4n+\delta mn}}{2} \Psi_{3,(k-1)}^{T_{k-1}} \leq \frac{1}{2}\epsilon_{k-1}, 
$$
and 
$$
\frac{\kappa}{2}\|e_{(k-1)}^{T_{k-1}}\|^2 \leq \frac{\sqrt{2n+\delta mn}}{2} \cdot \frac{\lambda+\kappa}{\sqrt{2n+\delta mn}} \|e_{(k-1)}^{T_{k-1}}\|^2 \leq  \frac{\sqrt{2n+\delta mn}}{2} \Psi_{3, (k-1)}^{T_{k-1}} \leq \frac{\epsilon_{k-1}}{2}. 
$$

\noindent Moreover, it is easy to see that $G_{k-1}(x^{k-1}) - G_{k-1}^* = G_{k-1}(x_{(k-1)}^{T_{k-1}+1}) - G_{k-1}^* \leq \frac{1}{2}\epsilon_{k-1}$. Therefore, 
$$
\epsilon_{D,(k)}^0 \leq G_k(x^{k-1}) - D_k(\alpha^{k-1}) \leq \frac{3}{2}\epsilon_{k-1} + \kappa \|y^{k-1} - y^{k-2}\|^2. 
$$
\noindent Then similar to the proofs of Proposition 3.2 and Lemma C.1 in \citep{catalyst}, we can get $\mathbb{E}[T_k] \leq {\tilde {\cal O}} \left( \frac{1}{\delta} + m + \frac{{\cal M}_1}{\lambda + \kappa} + \frac{1}{\delta} \sqrt{\frac{(1-\delta)({\bar R}^2 + \delta R_m^2)}{(\lambda+\kappa)\gamma}}  \right)$.

\subsection{Proof of Lemma \ref{lm:ecsdca-1-re}}

From the equality above (33) in the proof of Theorem 3.3 in \citep{ecsdca}, we know 
$$
{\tilde u}_{(k)}^{K+1} -  {\tilde u}_{(k)}^{K} = \frac{1}{(\lambda+\kappa)N} \sum_{\tau=1}^n A_{i_K^\tau \tau} \Delta \alpha_{i_K^\tau \tau, (k)}^{K+1},
$$
where ${\tilde u}_{(k)}^K = u_{(k)}^K + e_{(k)}^K$. Hence, as long as 
\begin{equation}\label{eq:uk0}
	{\tilde u}_{(k)}^{0} = \frac{1}{(\lambda+\kappa)N} \sum_{\tau=1}^n\sum_{i=1}^m A_{i\tau} \alpha_{i\tau,(k)}^0,
\end{equation}
we will have ${\tilde u}_{(k)}^{K} = \frac{1}{(\lambda+\kappa)N} \sum_{\tau=1}^n\sum_{i=1}^m A_{i\tau} \alpha_{i\tau,(k)}^K$ for all $K\geq 0$. From the initialization rule at each outer iteration, it is easy to see that (\ref{eq:uk0}) is satisfied for all $k\geq 1$. Therefore, 
\begin{equation}\label{eq:ukK}
	{\tilde u}_{(k)}^{K} = \frac{1}{(\lambda+\kappa)N} \sum_{\tau=1}^n\sum_{i=1}^m A_{i\tau} \alpha_{i\tau,(k)}^K,
\end{equation}
for all $k\geq 1$ and $K\geq 0$. Moreover, from the proofs for EC-SDCA, it is easy to verify that the convergence results still hold as long as (\ref{eq:ukK}) is satisfied, and it does not matter whether $u_{(k)}^0$ equal to $\frac{1}{(\lambda+\kappa)N} \sum_{\tau=1}^n\sum_{i=1}^m A_{i\tau} \alpha_{i\tau,(k)}^0$ or not. Then similar to the proof of Lemma \ref{lm:ecsdca-1}, we have 
\begin{align*}
	& \quad \mathbb{E} [\sqrt{4n+\delta mn} \Psi^K_{3,(k)} + 2(G_k(x_{(k)}^{K+1}) - G_k^*)] \\ 
	& \leq \left(  \sqrt{4n+\delta mn} + \frac{2}{\theta} \right) \left( 1 - \min\left\{ \theta, \frac{\delta}{4} \right\} \right)^K \left(\epsilon_{D,(k)}^0 + \frac{2(\rho+\theta(\lambda+\kappa))}{\delta n}\sum_{\tau=1}^n \|e_{\tau, (k)}^0\|^2 \right),
\end{align*}
and 
$$
\epsilon_{D,(k)}^0 \leq \frac{3}{2}\epsilon_{k-1} + \kappa \|y^{k-1} - y^{k-2}\|^2. 
$$
Moreover, from the initialization of $e_{\tau, (k)}^0$, we have 
$$
\frac{2(\rho+\theta(\lambda+\kappa))}{\delta n}\sum_{\tau=1}^n \|e_{\tau, (k)}^0\|^2 \leq \Phi_{3, (k-1)}^{T_{k-1}} \leq \frac{1}{2}\epsilon_{k-1}. 
$$
Then similar to the proofs of Proposition 3.2 and Lemma C.1 in \citep{catalyst}, we can get the result.

\clearpage 

\section{TABLES}

\begin{table*}[h]
	\caption{Communication Complexity Results for EC-LSVRG + Catalyst in the Smooth Case ($r_Q$ represents the communication cost of the compressed vector $Q(x)$ for $x\in \R^d$. The common Assumptions \ref{as:primal} and \ref{as:contracQQ1-ecSPDC} are omitted. $A_1 \eqdef \delta L_f + \nicefrac{\delta L}{n} + \sqrt{1-\delta} (\sqrt{L_f {\bar L}} + \sqrt{\delta L_f L})$ and $A_2 \eqdef \delta L_f + \nicefrac{\delta L}{n} + \sqrt{1-\delta}L_f$.) }
	\label{tab:ECLSVRG+C}
	\begin{center}
		{\footnotesize
			\begin{tabular}{|c|c|}
				\hline
				Assumptions & \begin{tabular}{c} Communication complexity 	\end{tabular}\\
				\hline 
				\begin{tabular}{c} $A_1 \geq \lambda$ \\ $\kappa=A_1-\lambda$ \end{tabular} & 
				${\tilde {\cal O}} \left(  \tfrac{r_Q}{\sqrt{\lambda}} \left(  \sqrt{\tfrac{L_f}{\delta}}  + \sqrt{\tfrac{L}{\delta n}}  + \tfrac{(\sqrt{1-\delta} (\sqrt{L_f {\bar L}} + \sqrt{\delta L_f L}))^{\tfrac{1}{2}}}{\delta}  + \tfrac{\sqrt{1-\delta} (\sqrt{\bar L} + \sqrt{\delta L})}{\delta}  \right)   \log \tfrac{1}{\epsilon} \right)$  \\
				\hline
				\begin{tabular}{c} $A_1 < \lambda$ \\ $\kappa=0$ \end{tabular}& 
				${\tilde {\cal O}} \left(  r_Q \left(  \tfrac{1}{\delta} + \tfrac{L_f}{\lambda} + \tfrac{L}{n \lambda} + \tfrac{\sqrt{1-\delta} (\sqrt{L_f {\bar L}} + \sqrt{\delta L_f L})}{\delta \lambda}  + \tfrac{\sqrt{1-\delta} (\sqrt{\bar L} + \sqrt{\delta L})}{\delta \sqrt{\lambda}}  \right) \log \tfrac{1}{\epsilon} \right)$  \\
				\hline
				\begin{tabular}{c} Assumption \ref{as:expcompressor} \\ $A_2 \geq \lambda$ \\ $\kappa=A_2-\lambda$ \end{tabular} & 
				${\tilde {\cal O}} \left(  \tfrac{r_Q}{\sqrt{\lambda}} \left(  \sqrt{\tfrac{L_f}{\delta}}  + \sqrt{\tfrac{L}{\delta n}}  + \tfrac{(1-\delta)^{\tfrac{1}{4}} \sqrt{L_f}}{\delta} \right)   \log \tfrac{1}{\epsilon} \right)$  \\
				\hline
				\begin{tabular}{c} Assumption \ref{as:expcompressor} \\ $A_2 < \lambda$ \\ $\kappa=0$ \end{tabular}& 
				${\tilde {\cal O}} \left(  r_Q \left(  \tfrac{1}{\delta} + \tfrac{L_f}{\lambda} + \tfrac{L}{n \lambda} + \tfrac{\sqrt{1-\delta} L_f}{\delta \lambda}   \right) \log \tfrac{1}{\epsilon} \right)$  \\
				\hline
			\end{tabular}
		}
	\end{center} 
\end{table*}

\begin{table*}[h]
	\caption{Communication Complexity Results for EC-SDCA + Catalyst ($r_Q$ represents the communication cost of the compressed vector $Q(x)$ for $x\in \R^d$. The common Assumptions \ref{as:contracQQ1-ecSPDC} and \ref{as:ecSPDC}  are omitted. $A_3 \eqdef (\frac{R_m^2}{n \gamma} + \frac{R^2}{\gamma} + \frac{\sqrt{1-\delta}R{\bar R}}{\delta \gamma} + \frac{\sqrt{1-\delta}RR_m}{\sqrt{\delta}\gamma}) \left/ (\tfrac{1}{\delta} + m) \right.$ and $A_4 \eqdef (\frac{R_m^2}{n \gamma} + \frac{R^2}{\gamma} + \frac{\sqrt{1-\delta}R^2}{\delta \gamma} ) \left/  (\tfrac{1}{\delta} + m) \right.$, where $R, {\bar R}, R_m$ are defined in Algorithm \ref{alg:ec-spdc}.) }
	\label{tab:ECSDCA+C}
	\begin{center}
		{\footnotesize
			\begin{tabular}{|c|c|}
				\hline
				Assumptions & \begin{tabular}{c} Communication complexity 	\end{tabular}\\
				\hline 
				\begin{tabular}{c} $A_3 \geq \lambda$ \\ $\kappa=A_3-\lambda$ \end{tabular} & 
				${\tilde {\cal O}} \left(  \tfrac{r_Q \log \tfrac{1}{\epsilon} }{\sqrt{\lambda}} \left(  \sqrt{\tfrac{1+\delta m}{\delta}} \sqrt{\frac{R_m^2}{n \gamma} + \frac{R^2}{\gamma} + \frac{\sqrt{1-\delta}R{\bar R}}{\delta \gamma} + \frac{\sqrt{1-\delta}RR_m}{\sqrt{\delta}\gamma}}  + \sqrt{\frac{(1-\delta)({\bar R}^2 + \delta R_m^2)}{ \delta^2\gamma}} \right)   \right)$  \\
				\hline
				\begin{tabular}{c} $A_3 < \lambda$ \\ $\kappa=0$ \end{tabular}& 
				${\tilde {\cal O}} \left(  r_Q \left(  \tfrac{1}{\delta} + m  + \frac{R_m^2}{\lambda n \gamma} + \frac{R^2}{\lambda \gamma} + \frac{\sqrt{1-\delta}R{\bar R}}{\delta \lambda \gamma} + \frac{\sqrt{1-\delta}RR_m}{\sqrt{\delta} \lambda \gamma} + \sqrt{\frac{(1-\delta)({\bar R}^2 + \delta R_m^2)}{\delta^2 \lambda \gamma}} \right) \log \tfrac{1}{\epsilon} \right)$  \\
				\hline
				\begin{tabular}{c} Assumption \ref{as:expcompressor} \\ $A_4\geq \lambda$ \\ $\kappa=A_4-\lambda$ \end{tabular} & 
				${\tilde {\cal O}} \left(  \tfrac{r_Q}{\sqrt{\lambda}} \left(  \sqrt{\tfrac{1}{\delta} + m} \sqrt{\frac{R_m^2}{n \gamma} + \frac{R^2}{\gamma} + \frac{\sqrt{1-\delta}R^2}{\delta \gamma}} \right)   \log \tfrac{1}{\epsilon} \right)$  \\
				\hline
				\begin{tabular}{c} Assumption \ref{as:expcompressor} \\ $A_4 < \lambda$ \\ $\kappa=0$ \end{tabular}& 
				${\tilde {\cal O}} \left(  r_Q \left(  \tfrac{1}{\delta} + m + \frac{R_m^2}{n \lambda \gamma} + \frac{R^2}{\lambda \gamma} + \frac{\sqrt{1-\delta}R^2}{\delta \lambda \gamma}   \right) \log \tfrac{1}{\epsilon} \right)$  \\
				\hline
			\end{tabular}
		}
	\end{center} 
\end{table*}

\end{document}